\documentclass[a4paper, 11pt, final, BCOR=5mm, pdftex, pagesize]{scrartcl} 
\usepackage[a4paper,left=3cm,right=3cm,top=2.9cm,bottom=3cm]{geometry}
\usepackage[utf8]{inputenc} 
\usepackage{fancyhdr}
\usepackage{csquotes} 
\usepackage{verbatim}
\usepackage{underscore}
\usepackage{enumitem,tabularx}
\usepackage{color}
\usepackage{fixltx2e} 
\usepackage{lmodern} 
\usepackage[T1]{fontenc} 
\usepackage{setspace}
\usepackage{authblk}

\usepackage{amsmath}
\usepackage{amsfonts}
\usepackage{amssymb}
\usepackage{mathrsfs}
\usepackage{mathtools}
\usepackage{dsfont} 
\usepackage{bigints} 
\usepackage{amsthm}
\usepackage{accents}

\usepackage{float} 
\usepackage[margin=0.8cm]{caption} 
\usepackage{subcaption}
\usepackage{graphicx}
\usepackage{pgfplots}
\pgfplotsset{compat=1.8}
\usepackage[draft=false,babel,tracking=true,kerning=true,spacing=true]{microtype} 
\usepackage{tabularx} 

\definecolor{CMblue}{RGB}{0, 153, 153}
\definecolor{CMgreen}{RGB}{76, 153, 0}
\definecolor{CMorange}{RGB}{237, 125, 49}
\definecolor{CMlightblue}{RGB}{204, 255, 255}

\usepackage{chngcntr}
\usepackage[section]{placeins} 
\numberwithin{equation}{section}

\usepackage[noadjust]{cite}
\usepackage{hyperref}

\usepackage{tikz}
\usetikzlibrary{shapes.geometric, patterns}
\pgfplotsset{compat=1.8}
\usepackage{setspace,tabularx} 
\usepackage[draft=false,babel,tracking=true,kerning=true,spacing=true]{microtype} 

\setitemize{leftmargin=0.5cm}

\newcommand{\R}{\mathbb{R}}
\newcommand{\N}{\mathbb{N}}
\newcommand{\Z}{\mathbb{Z}}
\newcommand{\F}{\mathcal{F}^{(n)}}
\newcommand{\E}{\ensuremath\mathbb{E}}
\newcommand{\1}{\textup{$\mathds{1}$}}
\newcommand{\Pro}{\ensuremath\mathbb{P}}

\newcommand{\Var}{\ensuremath\mathrm{Var}}
\newcommand{\Cov}{\ensuremath\mathrm{Cov}}
\newcommand{\id}{\ensuremath\mathrm{id}}
\newcommand{\Disc}{\ensuremath\mathrm{Disc}}

\newcommand{\sign}{\ensuremath\mathrm{sign}}

\renewcommand{\leq}{\leqslant}
\renewcommand{\geq}{\geqslant}
\renewcommand{\bar}{\overline}
\renewcommand{\tilde}{\widetilde}

\newcommand{\tn}{\Delta t^{(n)}}

\newcommand{\vn}{\Delta v^{(n)}}

\renewcommand{\t}{t^{(n)}}

\newcommand{\htau}{\hat{\tau}}

\newcommand{\dbtilde}[1]{\accentset{\approx}{#1}}

\newtheorem{The}{Theorem}[section]
\newtheorem*{The*}{Theorem}
\newtheorem{Ass}{Assumption}

\newtheorem{Def}{Definition}
\newtheorem{Lem}[The]{Lemma}
\newtheorem{Cor}[The]{Corollary}
\newtheorem{Rmk}[The]{Remark}
\newtheorem{Prop}[The]{Proposition}

\providecommand{\keywords}[1]
{
  \small	
  \textbf{\textit{Keywords---}} #1
}

\begin{document}

\title{A cross-border market model\\with limited transmission capacities}
\date{\today}

\author[1]{D\"orte Kreher\thanks{Dörte Kreher acknowledges support from DFG CRC/TRR 388 Rough Analysis, Stochastic Dynamics and Related Fields (Project B02), from DFG IRTG 2544 Stochastic Analysis in Interaction, and from DFG Berlin Mathematics Research Center MATH+ (Project AA4-9).}}
\author[2]{Cassandra Milbradt}
\affil[1,2]{Humboldt-Universit\"at zu Berlin, Germany. E-mail: \url{doerte.kreher@hu-berlin.de}, \url{cassandra.milbradt@gmail.com}}

\maketitle

\begin{abstract}\small
We develop a cross-border market model for two countries based on a continuous trading mechanism, in which the transmission capacities that enable transactions between market participants from different countries are limited. Our market model can be described by a regime-switching process alternating between active and inactive regimes, in which cross-border trading is possible respectively prohibited. Starting from a reduced-form representation of the two national limit order books, we derive a high-frequency approximation of the microscopic model, assuming that the size of an individual order converges to zero while the order arrival rate tends to infinity. If transmission capacities are available, the limiting dynamics are as follows: the queue size processes at the top of the two limit order books follow a four-dimensional linear Brownian motion in the positive orthant with oblique reflection at the axes. Each time the two best ask queues or the two best bid queues simultaneously hit zero, the queue size process is reinitialized. The capacity process can be described as a linear combination of local times and is hence of finite variation. The analytic tractability of the limiting dynamics allows us to compute key quantities of interest. 
\end{abstract}


\raggedbottom

\section{Introduction}\label{CB:sec:motSU}

Limit order books (LOBs) are a standard tool for price formation in modern financial markets. At any given point in time, an LOB depicts the number of unexecuted buy and sell orders at different price levels. The highest price a potential buyer is willing to pay is called the best bid price, whereas the best ask price is the smallest price of all placed sell orders. Incoming limit orders can be placed at many different price levels, while incoming market orders are matched against standing limit orders according to a set of priority rules. While limit order books have provided an efficient market matching mechanism in financial markets for many years, they are nowadays also used as a matching tool in intraday electricity markets. Intraday electricity markets are short-term markets that allow the market participants to balance their position from the day-ahead market until shortly before delivery, which is especially relevant for markets with a high share of renewables. \par

The mathematical analysis of the dynamics of a limit order market is a challenging task which has recently attracted many researchers in mathematical finance. One approach to study limit order books is based on an event-by-event description of the order flow as for example done in \cite{L03, CST10, CL13}. The resulting stochastic systems typically yield realistic models as they preserve the discrete nature of the dynamics, but turn out to be computationally challenging. To overcome this issue, others deal with continuum approximations of the order book, describing its time-dependent density through  
a stochastic partial differential equation, cf.~e.g.~\cite{MTW16, ContMueller}. Combining these two approaches, one can include suitable scaling constants in the microscopic order book dynamics and study its scaling behavior when the number of orders gets large while each of them is of negligible size. Scaling limits of the full order book dynamics (prices and volumes) are e.g.~studied in Horst and Paulsen \cite{HP17}, Horst and Kreher \cite{HorstKreherFluid, HorstKreherDiffusion}, Gao and Deng \cite{GD18}, and Kreher and Milbradt \cite{KM20}. These models can explain how the shape of the volume density function changes due to the incoming order flow. On the other hand, Cont and de Larrad \cite{CL21} study a reduced-form model of a limit order book, which assumes a constant spread and only describes the order flow at the top of the book, thereby providing closer insight into the price formation process. 
\par

While the mathematical description of a single LOB is already a complex task, real-world opportunities to trade the same asset at different venues simultaneously (so called cross-listed assets) require even more complex models. A special example of this is the continuous trading mechanism of the integrated European intraday electricity market ``Single Intraday Coupling'' (SIDC), launched in 2018. Its market mechanism is based on a shared order book and allows for cross-border transactions between market participants from different countries as long as transmission capacities are available.

In this paper, we introduce a microscopic model of a cross-border market between two countries with a continuous trading mechanism and limited transmission capacities. To the best of our knowledge, this is the first attempt in the literature to build a rigorous mathematical model of an integrated limit order market based on the underlying market microstructure. Our model is a further development of the one considered in Cont and de Larrard \cite{CL21} and will allow us to analyze the influence of the market coupling on the price dynamics. In \cite{CL21} the authors prove that under heavy traffic conditions the best bid and ask queues can be approximated by a planar Brownian motion in the first quadrant with inward jumps upon hitting the axes. 
We extend their model in multiple ways. First, we analyze the reduced-form representations of two LOB models simultaneously over time. Second, we allow domestic market orders to be matched against standing limit order volume in the foreign country, which leads to cross-border trades. Finally, motivated by the limitation of transmission capacities in the SIDC, we restrict the maximum cross-border trading volume. To this end, we keep track of the origin of each incoming order and introduce a capacity process that records the net cross-border trading volume over time. As the resulting microscopic market dynamics are too complex to get a good understanding of the market mechanism and the impact of the market coupling, we impose heavy traffic conditions and prove that our model can be approximated by a tractable, continuous-time regime switching process. Using this approximation we investigate how the cross-border trading opportunities affect the price evolution analytically and via simulations.\par

In more detail, our microscopic model is described through two best bid and ask queue size processes and two price processes, one for each country, as well as a two-sided capacity process. The state of the shared order book changes due to incoming market or limit orders. While the order flow directly affects the queue size processes, it indirectly also defines the price and capacity dynamics: if a \textit{single} queue is about to become negative due to an incoming market order, the queue will be set to zero and the remaining order size will be depleted from the corresponding queue in the foreign market as long as both, enough standing volume and enough transmission capacity are available. This will lead to a cross-border trade. If the \textit{cumulative} best bid or ask queue is about to be depleted by an incoming market order, both price processes will change by one tick and all queues will be reinitialized by some random variables representing the depth of the books. Each time a cross-border trade is executed, the capacity process is updated. Now, starting in a so-called \textit{active} regime in which cross-border trading is possible, we switch to a so-called \textit{inactive} regime if the capacity process hits one of its boundaries. In the inactive regime, market participants can only execute market orders against limit orders of the same origin. Throughout, we make the simplifying assumption that there are no transmission/interconnector costs. Under these idealized conditions the best bid and ask prices of both national LOBs coincide in the active regime, while they become different when switching to an inactive regime. Assuming an efficient allocation of capacities, i.e.~only the net volume of all executed cross-border trades will actually be transacted cross-border, the submission of off-setting orders makes it possible to switch back to an active regime. To illustrate our cross-border market model, we conduct a small simulation study. \par 

Our main theorem states that in the high frequency limit, i.e.~when the number of orders goes to infinity while each of them is of negligible size, the market dynamics are the following: in an active regime, the queue size process is described by a four-dimensional linear Brownian motion in the positive orthant with oblique reflection at the axes. Each time \textit{two} queues simultaneously hit zero, the \textit{whole} process is reinitialized at a new value in the interior of $\R^4_+.$ The bid price processes of both countries agree and follow a pure jump process, whose jump times are equal to those of the volume process. The capacity process is a continuous process of finite variation, constructed from the component-wise local time process of the queue size process. During an inactive regime the queue size process is described by a four-dimensional linear Brownian motion in the positive orthant. Each time it hits \textit{one} of the axes, the depleted component as well as the component corresponding to the other side of the LOB of the same country are \textit{both} reinitialized at a new value in the interior of $\R^2_+$, while the other two components stay unchanged. The bid price process follows a two-dimensional pure jump process, whose components do almost surely not jump simultaneously and only move at hitting times of the axes by the corresponding components of the queue size process.\par 

While the high-frequency approximation of the cross-border markets dynamics during inactive regimes can readily be deduced from \cite{CL21}, the analysis of the cross-border market dynamics during active regimes is more intricate as the used methodology - rewriting the queue size process as a regulated process of the order flow and the reinitalization values -  becomes more involved. In the active regime, the queue size process follows inbetween price changes the so called generalized processor sharing (GPS) discipline, which was introduced in \cite{PG93} as a model for the efficient sharing of a single processing resource among traffic in communication networks. 
Heavy traffic limit theorems for a large class of queuing networks can be conveniently analyzed by rewriting the dynamics as the solution to an associated Skorokhod problem (SP) and applying continuous mapping techniques. For polyhedral Skorokhod problems satisfying the completely-$\mathcal{S}$ condition this has been done in \cite{CM91,R84}. The GPS SP does not satisfy this condition and solutions of it generally only exist for functions of finite variation, cf.~\cite{DR99}. However, for general c\`adl\`ag traffic processes the GPS dynamics can be represented as the unique solution of a so called extended Skorokhod problem (ESP) introduced in \cite{R06}. This representation can then be used to derive its scaling limit, cf.~\cite{RR03}. \par

To study the active dynamics under heavy traffic conditions, we first represent the bid respectively ask side components of the queue size process between price changes as a series of local solutions to the one-dimensional Skorokhod problem at alternating axes and afterwards identify this construction with the solution to the GPS ESP. This allows us to apply continuous mapping techniques to prove that in the active regime the queue size process on each side of the shared order book can be approximated by a sum-conserving semimartingale-reflecting Brownian motion (SRBM) in the positive orthant, which jumps into the interior of the orthant upon hitting the origin. SRBMs were first studied in \cite{VW85} under the completely-$\mathcal{S}$ condition. As mentioned above the GPS reflection problem does not satisfy this condition and our sum-conserivng SRBM with absorption is indeed a special case of the dual skew symmetric SRBM studied in \cite{EFH21,FR22}. Our construction via a series of one-dimensional Skorokhod problems also allows us to gain a better understanding of the reflecting behaviour. For example, we are able to characterize the survival probability of the hitting time at the opposite axis as the solution of an interface problem for the two-dimensional heat equation. This result can be used to analyze how often country $F$ imports from respectively exports to $G$. Moreover, we provide an analytical formula for the distribution of the duration until the next price change and study its dependence on the model parameters. \par  

We also conduct a simulation study to see how the market coupling affects price stability in both markets. Altogether, the main financial insights from the analysis of our cross-border market model are the following: (i) under an efficient allocation of transmission capacities, the cross-border trading volume is of much smaller magnitude than the domestic trading volume, even if transmission capacities are unbounded; (ii) in an imbalanced market setting, capacities will be rather quickly occupied in one direction and regime switches will become rare events; 
(iii) while the liquidity in the shared order book is always higher than in each national order book, the exact impact of the cross-border trading mechanism on price stability depends on the precise correlation structure of the order flow process, which may have a stabilizing, but in some cases also a drift amplifying effect. \par

The literature concerning the modelling of interconnected electricity markets is rather sparse. Christensen and Benth \cite{CB20} propose a regime-switching reduced-form model for electricity prices in a two-market setup, 
in which prices are described by multiplying a two-dimensional Ornstein-Uhlenbeck process with a regime-dependent convergence matrix whose state is determined by an independent latent stochastic process. As opposed to this, in our model regime switches are triggered endogenously by keeping track of the utilized transmission capacities. 

\par

Last but not least, we note that there is a substantial literature on optimal trading strategies for cross-listed stocks (cf.~e.g.~\cite{MPW08,TY13}) and also some recent papers on optimal cross-border trading in intraday electricity markets (cf.~e.g.~\cite{CJQ19,CFVS22}). However, these models start from an ad hoc description of the price dynamics and the incurred permanent and transient price impact cost. Our focus in this paper is different: we want to understand the impact of the market architecture on the price formation process from first principles, i.e.~from the observed order flow. Indeed, empirical studies (cf. e.g. \cite{KP17, HPW16}) suggest that the incorporation of the market microstructure in any realistic model of intraday electricity markets is crucial to gain a better understanding of the market dynamics. Studying more realistic (and more complicated) order flow dynamics as well as optimal trading problems within our model is left for future research.

\subsection{Empirical motivation}\label{CB:sec:empEv}

In this subsection, we provide some empirical motivation for our microscopic model dynamics, looking at SIDC 
trade book data\footnote{The data set is publicly available at \href{https://www.epexspot.com/}{EPEX SPOT}. We consider the German and Austrian market areas.} (from June 30, 2020) and capacity flow data\footnote{The data set is publicly available at \href{https://www.omie.es/en/market-results/monthly/continuous-intradaily-market/capacity-occupation-import-export?scope=monthly&year=2021&month=1&country=3}{OMIE} and focuses on France and Spain.} (from January 2021). Since the start of SIDC (formerly XBID) in 2018, the traded volume as well as the overall number of trades on SIDC is steadily increasing. In November 2022, the overall monthly traded volume has exceeded 10 million MWh, while the number of trades per quarter has exceeded 20 million, cf.~\cite{SIDC}. 
Even though the SIDC market is still very illiquid compared to many financial markets, the development of the past years shows a clear trend toward significantly more liquid intraday electricity markets, which suggests that scaling limit approximations of the market dynamics will gain plausibility in the future.

\begin{figure}[h]
\centering 
\begin{subfigure}[b]{0.48 \textwidth}
\centering
\includegraphics[width = \textwidth]{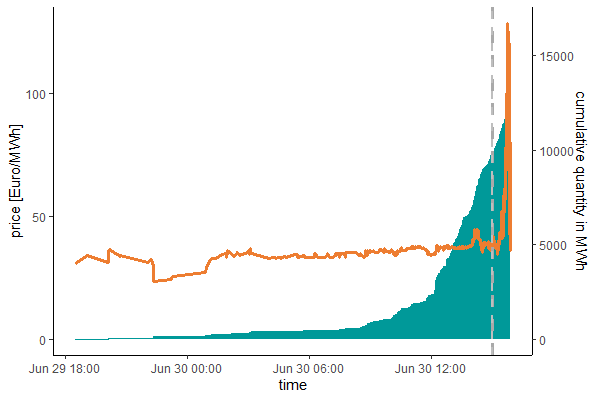}
\caption{\small Germany}
\label{CB:fig:transFlowDE}
\end{subfigure}
\hfill
\begin{subfigure}[b]{0.48 \textwidth}
\centering
\includegraphics[width = \textwidth]{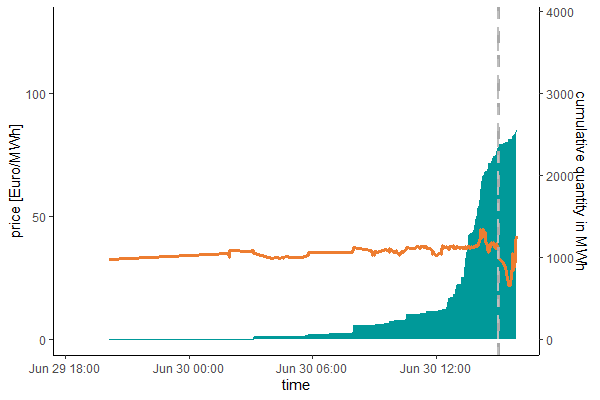}
\caption{\small Austria}
\label{CB:fig:transFlowAT}
\end{subfigure}
\caption{\small Transactions for the German and Austrian market area for delivery time $7$\,pm.\\ The $y$-axis displays all executed trades: prices (left), cumulative volume (right).}
\label{CB:fig:transFlow}
\end{figure}

Cross-border trading opportunities play an important role in the European intraday electricity market. For example, in our analyzed data set, Austria exported four times more electricity for delivery time 7pm to Germany as it transacted domestically. In contrast, Germany exported comparably few electricity to other countries for the same delivery window and most of its transactions were done domestically. Studying the transaction flow over time reveals the same picture. In Figure \ref{CB:fig:transFlow}, we depict the cumulative trading volumes and execution prices as a function of the execution time.  
The dashed gray lines in Figure \ref{CB:fig:transFlow} denote the last cross-border trade between both market areas. In both pictures, we observe that the execution prices rapidly change after the last transaction between Germany and Austria. 

\begin{figure}[h]
    \centering
    \includegraphics[scale = 0.53]{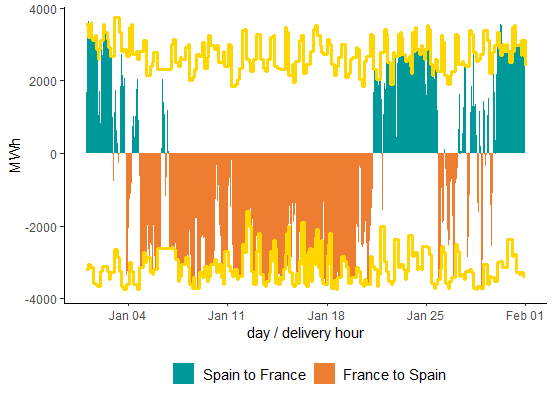}
    \caption{\small Available (in yellow) and occupied transmission capacity in January 2021.}
\label{CB:fig:SIDCmarketData5}
\end{figure}

To move toward a fully integrated European electricity market, the major obstacle to overcome is the limitation of transmission capacities. 
The rapid change of the prices in Figure \ref{CB:fig:transFlow} after the last cross-border trade between Germany and Austria indicates a full occupation of the transmission capacity in our considered data set.  Indeed, the available transmission capacities are typically occupied at the end of a trading session: Figure \ref{CB:fig:SIDCmarketData5} shows that for most continuously traded intraday electricity products, the transmission capacities between France and Spain were occupied at the end of the trading session in January 2021. Hence, incorporating the restriction of cross-border trading possibilities into the market model is necessary to develop a good understanding of the SIDC market dynamics.

\subsection{Outline}\label{CB:sec:outline}

The remainder of this paper is structured as follows: in Section \ref{CB:sec:disTimeMM} we define the discrete-time cross-border market dynamics $S^{(n)}$. Especially, we introduce the active market dynamics $\tilde{S}^{(n)}$ and the inactive market dynamics $\dbtilde{S}^{(n)}$, describing the evolution of the cross-border market dynamics when cross-border trades are possible respectively prohibited. We then state our main result (Theorem \ref{CB:res:mainTheorem}). A small simulation study of the cross-border market model is provided in Section \ref{CB:sec:simStudyCB}. In Section \ref{CB:sec:anaAD} we analyze the active dynamics $\tilde{S}^{(n)}$, identify their limit $\tilde{S}$, and study the price behaviour as well as the cross-border trading activity (cf.~Subsection \ref{CB:sec:ProbLimit}) for the limiting dynamics. In Section \ref{CB:sec:anaID}, we state the convergence result for the inactive dynamics $\dbtilde{S}^{(n)}$. Finally, Section \ref{CB:sec:proofmainResult} contains the proof of our main result, while Section \ref{CB:sec:conclusion} 
discusses a possible extension of our model to more than two countries.  Some technical details and proofs are transferred to the Appendix. \newline

\textbf{Notations.} $D(\R_+, \R^d)$ denotes the Skorokhod space of c\`adl\`ag functions $\omega: \R_+\rightarrow \R^d$ endowed with the usual Skorokhod metric $d_{J_1}$, cf.~\cite[Section 16]{B99}. 
For $k\in\N$ and $1\leq j\leq k$ the $j$-th projection map is given by $\pi^{(k)}_j : D(\R_+, \R^k) \rightarrow D(\R_+, \R)$, i.e.~$\pi^{(k)}_j \omega = \omega_j$ for $\omega = (\omega_1, \cdots, \omega_k) \in D(\R_+, \R^k)$. 
Similarly, $\pi^{(k)}_{i,j}: D(\R_+, \R^k) \rightarrow D(\R_+, \R^2)$, $k \geq 2,$ denots the $(i,j)$-th projection map, i.e.  $\pi^{(k)}_{i,j} \omega = (\omega_i, \omega_j)$ for $1\leq i, j \leq k.$ 
By a slight abuse of notation, we also write $\pi^{(k)}_j x = x_j$ and $\pi^{(k)}_{i,j} x = (x_i,x_j)$ for $x = (x_1, \cdots, x_k) \in \R^k$ and $1 \leq i,j \leq k$. 
Moreover, we introduce the short hand notations $\pi_j := \pi^{(4)}_j$ and $\pi_{i,j} := \pi^{(4)}_{i,j}$ for $1 \leq i,j \leq 4$ as well as $\pi_F := \pi_{1,2},$ $\pi_G := \pi_{3,4},$ $\pi_b := \pi_{1,3},$ $\pi_a := \pi_{2,4}$.\par
In the following, we make extensive use of the two summation functions $h$ and $h_1$ defined  
as 
\begin{align*}
  h_1&: D(\R_+, \R^2) \rightarrow D(\R_+, \R), &\quad h_1&: \omega \mapsto (\pi^{(2)}_1 \omega + \pi^{(2)}_2 \omega),\\
    h&: D(\R_+, \R^4) \rightarrow D(\R_+, \R^2), &\quad h&: \omega \mapsto 
    (h_1(\pi_b\omega),h_1(\pi_a\omega))=(\pi_1 \omega + \pi_3 \omega,\, \pi_2 \omega + \pi_4 \omega).
\end{align*} 
Throughout we set $\inf\emptyset:=\infty$ and $\infty-\infty:=\infty$. Furthermore, if two stochastic processes $X$ and $Y$ have the same finite-dimensional distributions, we write $X \simeq Y$. The symbol $\Rightarrow$ denotes weak convergence of probability measures resp.~convergence in distribution of random variables.

\section{Setup and main result}\label{CB:sec:disTimeMM}

We consider two neighboring countries $F$ (France) and $G$ (Germany), each of them having a national LOB through which it can trade its goods domestically. Moreover, as long as enough transmission capacities remain, market orders can also be matched against  standing volumes in the foreign LOB. In the following we describe our cross-border market dynamics by extending the \textit{reduced-form representation} of an LOB introduced in \cite{CL21} to two interacting LOBs: for each $n\in\N$ we consider on a stochastic basis $(\Omega^{(n)},\F,\Pro^{(n)})$\footnote{For notational convenience, we will write $\Pro$ instead of $\Pro^{(n)}$ throughout the paper.} a stochastic process $S^{(n)} = (S^{(n)}(t))_{t\geq 0}$ taking values in $E:= \R^2 \times \R^4_+ \times \R$ , i.e.
\begin{equation*}
S^{(n)}(t) =\left(B^{(n)}(t),\, Q^{(n)}(t),\, C^{(n)}(t)\right),
\end{equation*}
where the $\R^2$-valued process $B^{(n)} = (B^{F,(n)}, B^{G,(n)})$ describes the best bid prices in $F$ and $G$, the $\R^4_+$-valued process $Q^{(n)} = (Q^{b,F,(n)}, Q^{a,F,(n)},$ $Q^{b,G,(n)}, Q^{a,G,(n)})$ describes the standing volumes at the best bid/ask price in $F$ and $G,$ and the $\R$-valued process $C^{(n)}$ describes the net volume of all cross-border trades. We will refer to $C^{(n)}$ as the capacity process of the cross-border market model $S^{(n)}$. For all $n\in \N$, the bid price process $B^{(n)}$ takes values in $(\delta \Z)^2 \subset \R^2$, where $\delta>0$ denotes the tick size. As in \cite{CL21} we make the simplifying assumption that the spread is constant and equal to one tick at all times, i.e.~the best ask price processes are given by $A^{I,(n)} := B^{I,(n)} + \delta$ for $I = F,G$.

The state of the cross-border market changes due to arriving market and limit orders at the best bid/ask queues in $F,G$. For simplicity, we assume that the time intervals between two consecutive order arrivals are of equal length $\tn > 0$, i.e.~at each time $\t_k := k \tn, \ k \in\N$, an order book event takes place. Note however that an extension to randomly spaced arrival times as in \cite{HP17, CL21} is possible by invoking a time-change argument. We assume that $\tn\rightarrow0$ as $n\rightarrow \infty$, i.e.~the number of orders explodes as $n\rightarrow \infty.$ Moreover, all orders are assumed to be of constant size $\vn$, which also tends to zero as $n\rightarrow \infty$.

To study the cross-border market dynamics with limited transmission capacities we denote by $\kappa_- > 0$ the maximum transmission capacity in direction $F$ to $G$ (i.e.~\textit{exports} from $F$) and by $\kappa_+ > 0$ the maximum capacity in direction $G$ to $F$ (i.e.~\textit{imports} to $F$). Without loss of generality we suppose that $\kappa_+$ and $\kappa_-$ are multiples of $\vn$ for all $n\in\N$. We suppose that at time $t=0$ the bid prices of $F$ and $G$ coincide, i.e.~$B^{F,(n)}_0 = B^{G,(n)}_0$. Hence, in the $n$-th model the initial state of the cross-border market model is given by the $E$-valued random variable $S^{(n)}_0 =(B_0^{(n)}, Q^{(n)}_0, C^{(n)}_0)$, where $B^{(n)}_0=(B^{F,(n)}_0,B^{G,(n)}_0) \in (\delta,\delta) \Z$, $Q^{(n)}_0 = (Q^{b,F,(n)}_0,$ $Q^{a,F,(n)}_0, Q^{b,G,(n)}_0,$ $Q^{a,G,(n)}_0) \in (\vn \N)^4,$ and $C^{(n)}_0 \in[-\kappa_-,\kappa_+]\cap\vn\Z$ $\Pro$-a.s.

\begin{Ass}[Convergence of initial state]\label{CB:ass:initialCon}
	There exists a $(\delta,\delta)\Z\times(0,\infty)^4\times[-\kappa_-,\kappa_+]$-valued random variable $S_0 = \big(B_0, Q_0, C_0\big)$ 
    such that $S^{(n)}_0\Rightarrow S_0$ as $n\rightarrow\infty$.    
\end{Ass}

\subsection{Event types, order sizes, and the depth of the limit order books}

The state of the cross-border market model changes due to the submission of limit and market orders. While limit orders are added to the queues of orders awaiting execution, market orders reduce the queue lengths at the top of the book. To each order event we assign two random variables $\phi^{(n)}_k \in \{b,a\}$  and $\psi^{(n)}_k \in \{F,G\}$, which record the affected side of the LOB and the origin of the order: if the $k$-th order event affects the bid (resp.~ask) side of the LOB, then $\phi^{(n)}_k = b$ (resp.~$\phi^{(n)}_k = a$). If  the $k$-th order has origin $F$ (resp.~$G$), then $\psi^{(n)}_k = F$ (resp.~$\psi^{(n)}_k = G$). Hence, each event falls into one of the following four categories:
\begin{enumerate}[leftmargin = 1.5cm]
 \setlength\itemsep{0em}
	\item[$\textbf{(b,F).}$] A market sell or a limit buy order with origin $F$ arrives, i.e. $(\phi^{(n)}_k, \psi^{(n)}_k) = (b,F)$.
        \item[$\textbf{(a,F).}$] A market buy or a limit sell order with origin $F$ arrives, i.e.
        $(\phi^{(n)}_k, \psi^{(n)}_k) = (a,F)$.
        \item[$\textbf{(b,G).}$] A market sell or a limit buy order with origin $G$ arrives, i.e.
        $(\phi^{(n)}_k, \psi^{(n)}_k) = (b,G)$.
	\item[$\textbf{(a,G).}$] A market buy or a limit sell order with origin $G$ arrives, i.e. $(\phi^{(n)}_k, \psi^{(n)}_k) = (a, G)$.
\end{enumerate}
Market orders of type $(b,F)$ and $(a,G)$ correspond to possible transactions from $F$ to $G$, i.e. \textit{exports} from $F$, whereas $(a,F)$ and $(b,G)$ correspond to possible transactions from $G$ to $F$, i.e. \textit{imports} to $F$. In all cases, however, domestic transactions are preferred to foreign ones. \par

To distinguish between limit and market orders, we introduce the random sequence $(U^{(n)}_k)_{k\geq 1}$, where $U^{(n)}_k \in \{-\vn, \vn\}$ for all $k \geq 1$. If $U_k^{(n)}$ is positive, i.e.~$U^{(n)}_k = \vn$, the $k$-th order is a limit order placement. Conversely, if $U_k^{(n)}$ is negative, i.e.~$U^{(n)}_k = -\vn$, the $k$-th order is a market order. For notational convenience, the type-dependent order sizes will be denoted by 
\begin{align*}
V^{i,I,(n)}_k &:= U^{(n)}_k\1^{i,I}_{k,n} \quad\text{with}\quad \1^{i,I}_{k,n}:= \1_{\left\{\left(\phi^{(n)}_k, \psi^{(n)}_k\right) = (i, I)\right\}},\quad (i,I) \in \{b,a\} \times \{F,G\},\quad k \in\N,
\end{align*}
and the four-dimensional order size vector by $V^{(n)}_k := (V^{b,F,(n)}_k, V^{a,F,(n)}_k, V^{b,G,(n)}_k, V^{a,G,(n)}_k)$, $k\in\N$.\par

In order to derive a heavy traffic approximation for the cross-border market model $S^{(n)}$, we must impose further assumptions on the mean and covariance structure of the order size vectors $(V^{(n)}_k)_{k\geq 1}$. 
Combined with the right scaling relation between $\vn$ and $\tn$ (cf.~Assumption \ref{CB:ass:scaling}), this will guarantee that the partial sums of the order sizes verify some version of Donsker's theorem. 

\begin{Ass}[Order sizes]\label{CB:ass:prob} \,
 \begin{enumerate}
\item[i)]  For all $n \in \N$, $V^{(n)}=(V^{(n)}_k)_{k \geq 1}$ 
is a stationary, uniform mixing array of random variables which are independent of $S^{(n)}_0$.
 	\item[ii)] For all $n\in\N$ and $(i,I), (j,J) \in \{b,a\} \times \{F,G\}$ with $(i,I) \neq (j,J)$ 
  there exist $\mu^{i,I,(n)}\in \R$, $\sigma^{i,I,(n)} > 0$, and $\sigma^{(i,I),(j,J),(n)} \in \R$ such that
	\begin{align*}
	\E\left[V^{i,I,(n)}_1\right] &= \left(\vn\right)^2 \mu^{i,I,(n)},\\
	\Var\left[V^{i,I,(n)}_1\right] + 2 \sum_{k=2}^{\infty}\Cov\left[V^{i,I,(n)}_1, V^{i,I,(n)}_k\right] &= \left(\vn\right)^2 \left(\sigma^{i,I,(n)}\right)^2,\\
	 2 \sum_{k=1}^{\infty} \Cov\left[V^{i,I,(n)}_1, V^{j,J,(n)}_k\right]	 + 2 \sum_{k=2}^{\infty} \Cov\left[V^{i,I,(n)}_k, V^{j,J,(n)}_1\right] &= \left(\vn\right)^2\sigma^{(i,I), (j,J),(n)}.
	\end{align*}
	\item[iii)] There exist $\mu\in \R^4$ and a positive definite matrix $\Sigma\in \R^{4\times4}$ denoted by
   \begin{align*}
\mu = \left(\begin{array}{c}
\mu^{b,F}\\
\mu^{a,F}\\
\mu^{b,G}\\
\mu^{a,G}\end{array}\right),\quad
\Sigma = \left(\begin{array}{cccc}
\left(\sigma^{b,F}\right)^2 & \sigma^{(b,F),(a,F)} & \sigma^{(b,F),(b,G)} & \sigma^{(b,F),(a,G)}\\
\sigma^{(b,F),(a,F)} & \left(\sigma^{a,F}\right)^2 & \sigma^{(a,F),(b,G)} & \sigma^{(a,F),(a,G)}\\
\sigma^{(b,F),(b,G)} & \sigma^{(a,F),(b,G)} & \left(\sigma^{b,G}\right)^2 & \sigma^{(b,G),(a,G)}\\
\sigma^{(b,F),(a,G)} & \sigma^{(a,F),(a,G)} & \sigma^{(b,G),(a,G)} &\left(\sigma^{a,G}\right)^2
\end{array}\right).
\end{align*}
such that as $n \rightarrow \infty$ it holds for all $(i,I), (j,J) \in \{b,a\} \times \{F,G\}$ with $(i,I) \neq (j,J)$,
	 $$\left(\mu^{i,I,(n)}, \sigma^{i,I,(n)}\right) \rightarrow \left(\mu^{i,I}, \sigma^{i,I}\right) \quad \text{and} \quad \sigma^{(i,I),(j,J),(n)} \rightarrow\sigma^{(i,I),(j,J)}=:\rho^{(i,I),(j,J)}\sigma^{i,I} \sigma^{j,J}.$$
\end{enumerate}
\end{Ass}

The dependence structure introduced in Assumption \ref{CB:ass:prob} is for example satisfied if the sequence of scaled order sizes $(\vn)^{-1} (V^{(n)}_k)_{k\geq 1}$ becomes independent as $n\rightarrow \infty$ or if it is $m$-dependent for large $n\in \N$, i.e.~$(\vn)^{-1}(V^{(n)}_1, \cdots, V^{(n)}_k)$ and $(\vn)^{-1}(V^{(n)}_{k+j}, \cdots, V^{(n)}_{k+j+l})$ are independent for all $k, l \in \N$ whenever $j > m$. Note that the uniform mixing condition can be replaced by a $\rho$-mixing condition as discussed in \cite[Theorem 19.2]{B99} and is therefore stronger than the $\alpha$-mixing condition.\par

\begin{Ass}[Scaling]\label{CB:ass:scaling}
	As $n\rightarrow\infty$,	$\frac{\tn}{\left(\vn\right)^2} \rightarrow1$.
\end{Ass}

Under Assumptions \ref{CB:ass:prob} and \ref{CB:ass:scaling} we obtain convergence of the underlying net order flow process
\begin{align*}\label{CB:def:netorderflow}
X^{(n)} := \left(X^{b,F,(n)}, X^{a,F,(n)}, X^{b,G,(n)}, X^{a,G,(n)}\right),
\end{align*}
where for each $(i,I) \in \{b,a\} \times \{F,G\}$ and $t \geq 0$,
$$X^{i,I,(n)}(t): = \sum_{k=1}^{\infty} X^{i,I,(n)}_k \1_{\left[\t_k, \t_{k+1}\right)}(t) \quad \text{and} \quad X^{i,I,(n)}_k := \sum^k_{j=1} V_j^{i,I,(n)}.$$
The components of $X^{(n)}(t)$ represent the sum of all market and limit orders at each side of the order book in countries $F$ and $G$, which have been submitted up to time $t$. In Section \ref{CB:sec:anaAD}  we will derive a pathwise representation of $S^{(n)}=(B^{(n)},Q^{(n)},C^{(n)})$ in terms of the the underlying net order flow process $X^{(n)}$ and the reinitialization values after price changes defined below, cf.~Theorem \ref{CB:res:repactive}. 

\begin{Prop}[Convergence of the net order flow process]\label{CB:res:fCLT}
Under Assumptions \ref{CB:ass:prob}, \ref{CB:ass:scaling} the process $X^{(n)}$ converges weakly in the Skorokhod topology to a four-dimensional linear Brownian motion, i.e.
	\begin{align*}
	X^{(n)} \Rightarrow X := \left(\Sigma^{1/2} W(t) + \mu t\right)_{t\geq 0},
	\end{align*}
where $W$ is a standard four-dimensional Brownian motion.
\end{Prop}

The proof of Proposition \ref{CB:res:fCLT} follows from a classical functional convergence result for sums of weakly dependent random variables (cf.~\cite[Theorem 19.1]{B99}) and is hence omitted. \par

To describe the dynamics of the queue size process $Q^{(n)}$ we must define its value after price changes. Motivated by pegged limit orders, we follow the approach from \cite{CL21} and allow the queue sizes after a price change to depend on the state of the queue size process just before the price change. In the following let $(\tau^{(n)}_k)_{k\geq 1}$ denote the sequence of stopping times, at which we observe price changes in the cross-border market $S^{(n)}$, and let $R^{+,(n)}_k$ (resp.~$R^{-,(n)}_k$) denote the queue size after the $k$-th price change in case of a price increase (resp.~decrease).

\begin{Ass}[Reinitialization after price changes]\label{CB:ass:R} $ $
	\begin{enumerate}
	\item[i)] For all $n \in \N$ there exist two independent sequences $(\epsilon^{+,(n)}_k)_{k\geq 1}$ and $(\epsilon^{-,(n)}_k)_{k\geq 1}$ of $(0,\infty)^4$-valued, i.i.d. random variables on $(\Omega^{(n)},\F,\Pro^{(n)})$, which are independent of $V^{(n)}$ and $S^{(n)}_0$, and a function $\Phi^{(n)} : \R^4_+ \times \R^4_+ \rightarrow (\vn (\N+1))^4$ such that the reinitialization value of the queue size process after a price increase/decrease is given by
    \begin{align*}
    \begin{aligned}
    R^{+,(n)}_k := \Phi^{(n)}\left(Q^{(n)}(\tau^{(n)}_k-), \epsilon^{+,(n)}_k\right),\qquad
     R^{-,(n)}_k := \Phi^{(n)}\left(Q^{(n)}(\tau^{(n)}_k-), \epsilon^{-,(n)}_k\right).
    \end{aligned}
    \end{align*}	
    Moreover, there exists $\alpha>0$ such that for $j = 1,2,3,4$ and for all $(x,y) \in \R^4_+ \times \R^4_+$,
    	\[\pi_j\Phi^{(n)}(x,y) \geq \alpha \pi_j y.\]
	In the following, let $f^+_n$ (resp.~$f^-_n$) denote the distribution of $\epsilon^{+,(n)}_k$ (resp.~$\epsilon^{-,(n)}_k$), $k\in\N$. 
	\item[ii)]  There exist probability distributions $f^+,f^-$ on $(0,\infty)^4$ and $\Phi \in C\left(\R^4_+ \times \R^4_+, (0,\infty)^4 \right)$ such that
		\begin{align*}
		f^+_n \Rightarrow f^+,\quad f^-_n \Rightarrow f^-, \quad \text{ and } \quad \|\Phi^{(n)} - \Phi\|_{\infty} \rightarrow 0 \quad \text{ as } n \rightarrow \infty.
		\end{align*}
	\end{enumerate}
\end{Ass}

\subsection{Microscopic market dynamics}

To describe the dynamics of the cross-border market model $S^{(n)}$, we first make the following observation: if enough transmission capacities remain, the national LOBs are \textit{coupled}, i.e.~incoming market orders can be matched against standing volumes of the national \textit{and} the foreign LOB. As  the transmission capacities in both directions are bounded by $\kappa_-, \kappa_+ >0$, it may happen that the national LOBs \textit{decouple}, in which case market orders can only be matched against limit orders with the same origin. Hence, our cross-border market model switches between the following two regimes:

\begin{itemize}
 \setlength\itemsep{0em}
    \item the \textbf{active} regime in which the LOBs of $F$ and $G$ are coupled, and
    \item the \textbf{inactive} regime in which the LOBs of $F$ and $G$ are decoupled.
\end{itemize}
To describe the behaviour of $S^{(n)}$ during the two regimes, we define the \textbf{active dynamics} $\boldsymbol{\tilde{S}^{(n)}}$ in Section \ref{CB:sec:defAD} and the \textbf{inactive dynamics} $\boldsymbol{\dbtilde{S}^{(n)}}$ in Section \ref{CB:sec:defID}, which describe the evolution of the two national LOBs as if we were in the active respectively inactive regime for all times.

\subsubsection{Active dynamics}\label{CB:sec:defAD}

We define for fixed $n\in\N$ the active dynamics given by the piecewise constant interpolation
\[\tilde{S}^{(n)}(t) = \sum_{k=0}^{\infty}\tilde{S}^{(n)}_k \1_{\left[\t_k, \t_{k+1}\right)}(t), \quad t \geq 0,\]
of the $E$-valued random variables 
\[\tilde{S}^{(n)}_k := \left(\tilde{B}^{(n)}_k, \tilde{Q}^{(n)}_k, \tilde{C}^{(n)}_k\right), \quad k \in \N_0,\]
where $\tilde{B}^{(n)}_k$ denotes the bid prices of $F$ and $G,$ $\tilde{Q}^{(n)}_k$ denotes the sizes of the best bid and ask queues in $F$ and $G,$ and $\tilde{C}^{(n)}_k$ denotes the net volume of all cross-border trades after $k$ order events. We assume that there are \textbf{no capacity constraints} ($\kappa_+=\kappa_-=\infty$) and therefore the two limit order markets of $F$ and $G$ are always coupled. Especially, prices  agree throughout, i.e.~$\tilde{B}^{F,(n)}_k = \tilde{B}^{G,(n)}_k$ for all $k\geq 0$. 
In the shared order book, the cumulative queue size processes are given by
$$\tilde{Q}^{i,(n)}:= \tilde{Q}^{i,F,(n)}+ \tilde{Q}^{i,G,(n)},\quad i=b,a.$$

\begin{Def}[Evolution of the active dynamics $\tilde{S}^{(n)}_{k-1}\leadsto\tilde{S}^{(n)}_k$]\label{CB:def:active}
Let $I \in \{F,G\}$ denote the domestic country, $J \in \{F, G\} \setminus I$ the foreign country, and $i\in \{b,a\}$ the affected side of the order book. If the $k$-th order type is $(\phi_k^{(n)},\psi_k^{(n)})=(i,I)$ and if $(l-1)\in \N_0$ price changes have occurred before, then:
\begin{itemize}
\item If $V^{i,I,(n)}_k=+\vn$, a limit order placement occurs and $\tilde{Q}_k^{(n)}=\tilde{Q}_{k-1}^{(n)}+V_k^{(n)}$, $\tilde{B}^{(n)}_k=\tilde{B}^{(n)}_{k-1}$,  $\tilde{C}^{(n)}_k=\tilde{C}^{(n)}_{k-1}$.
    \item If $V_k^{i,I,(n)}=-\vn$ and 
    $\tilde{Q}^{i,I,(n)}_{k-1} \geq \vn$, a domestic trade occurs and $\tilde{C}^{(n)}_k=\tilde{C}^{(n)}_{k-1}$.\\
    If moreover $\tilde{Q}^{i,(n)}_{k-1}>\vn$, no price change occurs and $\tilde{B}^{(n)}_k=\tilde{B}^{(n)}_{k-1}$, $\tilde{Q}_k^{(n)}=\tilde{Q}_{k-1}^{(n)}+V_k^{(n)}$. 
\item If $V_k^{i,I,(n)}=-\vn$ and the domestic queue is empty, i.e. $\tilde{Q}^{i,I,(n)}_{k-1}=0$, a cross-border trade occurs. In this case, $\tilde{C}^{(n)}_k = \tilde{C}^{(n)}_{k-1} + \vn (\1_{(i,I)\in\{(a,F),(b,G)\}}- \1_{(i,I)\in\{(a,G),(b,F)\}})$.\\
    If moreover $\tilde{Q}^{i,J,(n)}_{k-1} > \vn$, no price change occurs, i.e.~$\tilde{B}^{(n)}_k=\tilde{B}^{(n)}_{k-1}$, and we have $\tilde{Q}^{i,J,(n)}_k = \tilde{Q}^{i,J,(n)}_{k-1} + V^{i,I,(n)}_k$ and $\tilde{Q}^{i,I,(n)}_k = \tilde{Q}^{i,I,(n)}_{k-1} =0$, $\tilde{Q}^{j,H,(n)}_k=\tilde{Q}^{j,H,(n)}_{k-1}$ for $j\neq i$ and $H=F,G$.

\item     If $V_k^{i,I,(n)}=-\vn$ and $\tilde{Q}^{i,(n)}_{k-1} = \vn$, the $l$-th price change occurs. In this case, 
    $\tilde{B}^{I,(n)}_k = \tilde{B}^{J,(n)}_k = \tilde{B}^{I,(n)}_{k-1} + \delta (\1_{i=a} - \1_{i=b})$ and all queues are reinitialized, i.e.~$\tilde{Q}^{(n)}_{k}=\tilde{R}^{+,(n)}_{l} \1_{i=b}+\tilde{R}^{-,(n)}_{l} \1_{i=a}$, where
$\tilde{R}^{\pm,(n)}_{l} := \Phi^{(n)}\left(\tilde{Q}^{(n)}_{k-1}, \epsilon^{\pm,(n)}_l\right)$.

\end{itemize}
\end{Def}

The definition of the active dynamics corresponds to an efficient cross-border trading mechanism in the sense that market orders are only matched against foreign standing volumes if no domestic volumes are available for the same price. 
Moreover, we keep track of the net volume of cross-border trades over time via the (unbounded) capacity process $\tilde{C}^{(n)}$. Its evolution can be described by
\begin{equation}\label{CB:def:capacity1}
\tilde{C}^{(n)}_k = \tilde{C}^{(n)}_0 + \underbrace{\left(M^{b,G, (n)}_k + M^{a,F, (n)}_k\right)}_{\text{trades from $G$ to $F$}} - \underbrace{\left(M^{b,F, (n)}_k + M^{a,G, (n)}_k\right)}_{\text{trades from $F$ to $G$}},
\end{equation}
where for $(i,I) \in \{b,a\} \times \{F,G\}$ the type-dependent cumulative cross-border volume is defined by
\begin{align*}
M_k^{i,I,(n)} &:= \vn \sum_{j=1}^k \1_{\left\{\tilde{Q}^{i,I,(n)}_{j-1} = 0, \, V^{i,I,(n)}_j = -\vn\right\}},\quad k\in\N.
\end{align*}
In the following, we denote 
\begin{equation}\label{CB:def:sNoCB}
\begin{split}
M^{(n)}_k &:= \left(M^{b,F,(n)}_k, M^{a,F,(n)}_k, M^{b,G,(n)}_k, M^{a,G,(n)}_k\right),\quad k\geq 1,\\
M^{(n)}(t)&:=\sum_{k=1}^{\infty}M_k^{(n)}\1_{\big[t^{(n)}_k,t^{(n)}_{k+1}\big)}(t),\quad t\geq 0.
\end{split}
\end{equation}

\begin{Rmk}
Between price changes the queue size dynamics in the active regime follows the so called generalized processor sharing (GPS) discipline, which was introduced in \cite{PG93} as a model for efficient, work-conserving resource sharing. 
\end{Rmk}

\subsubsection{Inactive dynamics}\label{CB:sec:defID}

In this subsection, we describe the inactive dynamics given by the piecewise constant interpolation
\[\dbtilde{S}^{(n)}(t) = \sum_{k=0}^{\infty} \dbtilde{S}^{(n)}_k \1_{\left[\t_k, \t_{k+1}\right)}(t), \quad t\geq 0,\]
of the $E$-valued random variables 
\[\dbtilde{S}^{(n)}_k := \left(\dbtilde{B}^{(n)}_k, \dbtilde{Q}^{(n)}_k, \dbtilde{C}^{(n)}_k\right), \quad k\in \N_0,\]
where $\dbtilde{B}^{(n)}_k$ denotes the bid prices in $F$ and $G,$ $\dbtilde{Q}^{(n)}_k$ denotes the sizes of the best bid and ask queues in $F$ and $G$, and $\dbtilde{C}^{(n)}_k$ denotes the net volume of all cross-border trades after $k$ order events. Since the national order books are decoupled in the inactive dynamics, the depletion of a single national order queue already causes a price change in the corresponding national LOB. Hence, the best bid prices in $F$ and $G$ will generally not coincide. Finally, since no cross-border trades are possible, the capacity process $\dbtilde{C}^{(n)}$ is constant for the whole trading period, i.e.~$\dbtilde{C}^{(n)}_k = \dbtilde{C}^{(n)}_0$ for all $k \geq 0$.

\begin{Def}[Evolution of the inactive dynamics $\dbtilde{S}^{(n)}_{k-1}\leadsto\dbtilde{S}^{(n)}_k$]\label{CB:def:inactive} 
Let $I \in \{F,G\}$ denote the domestic country, $J \neq I$ the foreign country, and $i\in \{b,a\}$ the affected side of the order book. If the $k$-th order type is $(\phi_k^{(n)},\psi_k^{(n)})=(i,I)$ and if $(l-1)\in \N_0$ price changes have occurred before, then:
\begin{itemize}
\item No cross-border trades are possible and no changes in the foreign market occur, i.e. $\dbtilde{C}^{(n)}_k = \dbtilde{C}^{(n)}_{k-1}$, $\dbtilde{B}^{J,(n)}_k=\dbtilde{B}_{k-1}^{J,(n)}$ and $\tilde{Q}^{i,J,(n)}_k=\tilde{Q}^{i,J,(n)}_{k-1}$ for $i=b,a$.
\item If $V^{i,I,(n)}_k=+\vn$, a limit order placement occurs and $\dbtilde{B}^{F,(n)}_k=\dbtilde{B}^{F,(n)}_{k-1},$ $\dbtilde{Q}_k^{(n)}=\dbtilde{Q}_{k-1}^{(n)}+V_k^{(n)}$.
    \item If $V_k^{i,I,(n)}=-\vn$ and 
    $\dbtilde{Q}^{i,I,(n)}_{k-1} \geq \vn$, a domestic trade occurs.
    If moreover $\dbtilde{Q}^{i,I,(n)}_{k-1}>\vn$, no price change occurs and $\dbtilde{B}^{F,(n)}_k=\dbtilde{B}^{F,(n)}_{k-1}$, $\dbtilde{Q}_k^{(n)}=\dbtilde{Q}_{k-1}^{(n)}+V_k^{(n)}$. 
    \item If $V_k^{i,I,(n)}=-\vn$ and 
    $\dbtilde{Q}^{i,I,(n)}_{k-1}= \vn$, the $l$-th price change occurs. In this case, $\dbtilde{B}^{I,(n)}_k = \dbtilde{B}^{I,(n)}_{k-1} + \delta ( \1_{i=a} - \1_{i=b})$ and both national queues are reinitialized, i.e.~$\pi_I\tilde{Q}^{(n)}_{k}=\pi_I\dbtilde{R}^{+,(n)}_{l} \1_{i=b}+\pi_I\dbtilde{R}^{-,(n)}_{l} \1_{i=a}$, where
$\dbtilde{R}^{\pm,(n)}_{l} := \Phi^{(n)}\left(\dbtilde{Q}^{(n)}_{k-1}, \epsilon^{\pm,(n)}_l\right)$.
\end{itemize}
\end{Def}
\noindent
Note that the inactive dynamics evolve analogously to those in \cite{CL21} for a single LOB.

\subsubsection{Cross-border market dynamics} \label{CB:sec:CBMD}

In this subsection we finally introduce the regime switching dynamics of our cross-border market model $S^{(n)}$. To this end, we have to specify when to switch from the active dynamics $\tilde{S}^{(n)}$ to the inactive dynamics $\dbtilde{S}^{(n)}$ and vice versa. While it seems natural to switch from an active to an inactive regime as soon as the unconstrained capacity process $\tilde{C}^{(n)}$ leaves the interval $[-\kappa_-,\kappa_+]$, the capacity process $\dbtilde{C}^{(n)}$ stays constant during an inactive regime and can therefore not be used to decide when to switch back. However, recall that in an active regime prices $\tilde{B}^{F,(n)}$ and $\tilde{B}^{(n),G}$ are always equal, while they diverge immediately upon entering an inactive regime. Hence, before reverting back to an active regime, $\dbtilde{B}^{F,(n)}$ and $\dbtilde{B}^{G,(n)}$ have to coincide again. 
Therefore, the price difference $\dbtilde{B}^{F,(n)}-\dbtilde{B}^{G,(n)}$ can be used to indicate when to switch back to an active regime. Note however that in the active regime the price difference $\tilde{B}^{F,(n)}-\tilde{B}^{G,(n)}$ always equals zero and hence it cannot be used to indicate a regime switch to an inactive regime.

To formalize this idea we define a sequence of random variables $(\tilde{Z}^{(n)}_k)_{k \geq 1} := (\tilde{Z}^{F,(n)}_k, \tilde{Z}^{G,(n)}_k)_{k\geq 1}$, which determine the state of the cross-border market dynamics at the start of an inactive regime, by
\begin{equation}\label{CB:def:Zn}
\begin{split}
    \tilde{Z}^{I,(n)}_k := \begin{cases}
    \textcolor{white}{+}1&:\ \pi_IQ^{(n)}_{k-1} + \pi_I V_k^{(n)} \in \vn\N_0\times \{-\vn\}\\    
    -1 &:\ \pi_IQ^{(n)}_{k-1} + \pi_I V_k^{(n)} \in \{-\vn\} \times \vn\N_0\\    
    \textcolor{white}{+}0 &:\ \text{otherwise } 
    \end{cases}
\end{split}
\end{equation}
and denote their piecewise constant interpolation by $\tilde{Z}^{(n)} (t)=\sum_{k=1}^{\infty}\tilde{Z}_k^{(n)}\1_{\big[t^{(n)}_k,t^{(n)}_{k+1}\big)}(t)$ for $t\geq 0.$ 
Moreover, we denote for all $t\geq0$ by 
\begin{align*}
    &l^{(n)}(t) \text{ -- the random number of price changes in $S^{(n)}$ up to time $t$},\\
    & \tilde{S}^{(n),t} \text{ -- the active dynamics with starting value $S^{(n)}(t)$},\\
    & \dbtilde{S}^{(n),t} \text{ -- the inactive dynamics with starting value $S^{(n)}(t)$.}
\end{align*}

\begin{Def}[Cross-border market dynamics]\label{CB:def:Sn}
For every $n\in\N$, we set 
\[S^{(n)}(t):=\sum_{k=1}^\infty \tilde{S}^{(n), \rho^{(n)}_{k-1}}\left(t- \rho^{(n)}_{k-1}\right)\1_{\left[\rho^{(n)}_{k-1}, \sigma^{(n)}_k\right)}(t)+\dbtilde{S}^{(n), \sigma^{(n)}_k}\left(t-\sigma^{(n)}_k\right)\1_{\left[\sigma^{(n)}_k,\rho^{(n)}_k\right)}(t),\]
where $\rho^{(n)}_0 := 0,$ and for $k\in\N$, the $k$-th inactive resp.~active regime starts at 
        \begin{align*}
                \sigma^{(n)}_k &:= \rho^{(n)}_{k-1}+\inf\Big\{t_l^{(n)} > 0: \tilde{C}^{(n), \rho^{(n)}_{k-1}}_l \notin [-\kappa_-,\kappa_+] \Big\},\\
                \rho^{(n)}_k &:= \sigma_k^{(n)}+\inf\left\{t_l^{(n)} > 0: \dbtilde{B}^{F,(n),\sigma_k^{(n)}}_{l-1} = \dbtilde{B}^{G,(n),\sigma_k^{(n)}}_{l-1}\right\}.
        \end{align*}
Moreover, $S^{(n)}(0):=S^{(n)}_0$ and for $k\in\N$, the starting value of the $k$-th inactive resp.~active regime is given for $I=F,G$ by
        \begin{align*}
    B^{I,(n)}(\sigma^{(n)}_k) &:= B^{I,(n)}(\sigma^{(n)}_k-) + \delta\left(\1_{\{\tilde{Z}^{I,(n)}(\sigma^{(n)}_k) = 1\}} -\1_{\{\tilde{Z}^{I,(n)}(\sigma^{(n)}_k) = -1\}}\right),\\
\pi_IQ^{(n)}(\sigma^{(n)}_k) &:=
    \pi_I Q^{(n), \rho^{(n)}_{k-1}}\big(\sigma^{(n)}_k-\big)\1_{\{\tilde{Z}^{I,(n)}(\sigma^{(n)}_k) = 0\}}+
    \pi_I R^{+,(n)}_{l^{(n)}(\sigma^{(n)}_k-)+1}\1_{\{\tilde{Z}^{I,(n)}(\sigma^{(n)}_k) = 1\}}\\
    &\qquad\qquad+
    \pi_I R^{-,(n)}_{l^{(n)}(\sigma^{(n)}_k-)+1} \1_{\{\tilde{Z}^{I,(n)}(\sigma^{(n)}_k) = -1\}}
+ \pi_IV^{(n)}_{\sigma_k^{(n)}/\tn},\\
C^{(n)}(\sigma^{(n)}_k) &:= C^{(n)}(\sigma^{(n)}_k-),\\
\end{align*}
respectively
        \begin{align*}
        S^{(n)}(\rho^{(n)}_k) &:= \left(B^{(n)}(\rho^{(n)}_k-),\; Q^{(n)}(\rho^{(n)}_k-) + V^{(n)}_{\rho^{(n)}_k/\tn},\, C^{(n)}(\rho^{(n)}_k -)\right).        
        \end{align*}
\end{Def}

\begin{Rmk}\label{CB:rem:regime-dependent-parameters}
	In our model, a regime switch only changes the order matching mechanism. However, it might be interesting to also consider a change in the order flow. This could be incorporated by a slight generalization of Assumptions \ref{CB:ass:prob} and \ref{CB:ass:R}, assuming that the distribution of the order sizes and of the reinitialization values after price changes depend on the type of the current regime.
\end{Rmk}

\subsection{Main result}

Our main result states the weak convergence of the sequence $(S^{(n)})_{n\in\N}$ to some $S = (B, Q, C)$. To describe the dynamics of $S$, let us introduce the continuous-time analogues $(\rho_k)_{k\geq 0}$ and $(\sigma_k)_{k\geq 1}$ of the $\tn\N_0$-valued stopping times $(\rho^{(n)}_k)_{k\geq 0}$ and $(\sigma^{(n)}_k)_{k\geq 1}.$ 
We set $\rho_0 := 0$ and for $k\geq 1$, 
\begin{equation}\label{CB:eq:sigmarho}
\sigma_k := \inf\left\{t\geq\rho_{k-1}:\ C(t)\notin[-\kappa_-,\kappa_+]\right\},\qquad
\rho_k := \inf\left\{t \geq \sigma_k: B^F(t) = B^G(t)\right\}.
\end{equation}

\begin{The}[Main result]\label{CB:res:mainTheorem}
    Let Assumptions \ref{CB:ass:initialCon}--\ref{CB:ass:R} be satisfied. Then the microscopic cross-border market models $S^{(n)} = (S^{(n)}(t))_{t\geq 0},\ n\in\N,$ converge weakly in the Skorokhod topology on $D(\R_+, E)$ to a continuous-time regime switching process $S$, whose dynamics is as follows: for all $k \geq 0,$
    \begin{itemize}
    		\item on the interval $[\rho_k, \sigma_{k+1})$ the queue size process $Q$ is a four-dimensional linear Brownian motion in the positive orthant with oblique reflection at the axes.  
      Each time the two bid queues or the two ask queues simultaneously hit zero, the process is reinitialized at a new value in the interior of $\R^4_+.$ The bid prices $B$ follow a two-dimensional pure jump process with jump times equal to those of the queue size process. In particular, we have $B^F \equiv B^G.$ The capacity process $C$ is a continuous process of finite variation with values in $[-\kappa_-, \kappa_+].$ 
    		\item on the interval $[\sigma_{k+1}, \rho_{k+1})$ the queue size process $Q$ behaves like a four-dimensional linear Brownian motion in the interior of $\R^4_+.$ Each time it hits one of the axes, the two components corresponding to the origin of the depleted component are reinitialized at a new value in $(0,\infty)^2$ while the others stay unchanged. The bid prices $B$ follow a two-dimensional pure jump process whose components jump at hitting times of the axes by the corresponding components of the queue size process. In particular, $B^F$ and $B^G$ are pure jump processes which do almost surely not jump simultaneously. The capacity process $C$ stays constant at either $-\kappa_-$ or $\kappa_+.$
    \end{itemize}
\end{The}

To prove our main result, we will first analyze the active and inactive dynamics separately in Sections \ref{CB:sec:anaAD} and \ref{CB:sec:anaID} below. The final proof of Theorem \ref{CB:res:mainTheorem} can be found in Section \ref{CB:sec:proofmainResult}. Before we turn to the theoretical analysis of our model, we present a small simulation study of the model dynamics in the next section.

\section{Simulation study}\label{CB:sec:simStudyCB}

In this section, we consider a finite time horizon $[0,1]$ and simulate our market dynamics for different parameter constellations of the order flow process. We choose $\tn = n^{-1}$, $\vn = n^{-1/2}$, and simulate $n = 10,000$ time steps. For simplicity, we assume that the order size vectors $(V_k^{(n)})_{k\geq 1}$ are independent over time.

\subsection{Simulation of the full cross-border market dynamics}

In this subsection, we assume that $\kappa_+ = \kappa_- = 0.5$, $\delta = 0.1,$ and that the reinitialization values, i.e.~the queue sizes after price changes, are drawn from the uniform distribution on $\{j\vn: j = 10,\dots, 20\}$. \newline

\begin{figure}[h]
    \centering
    \includegraphics[scale = 0.3]{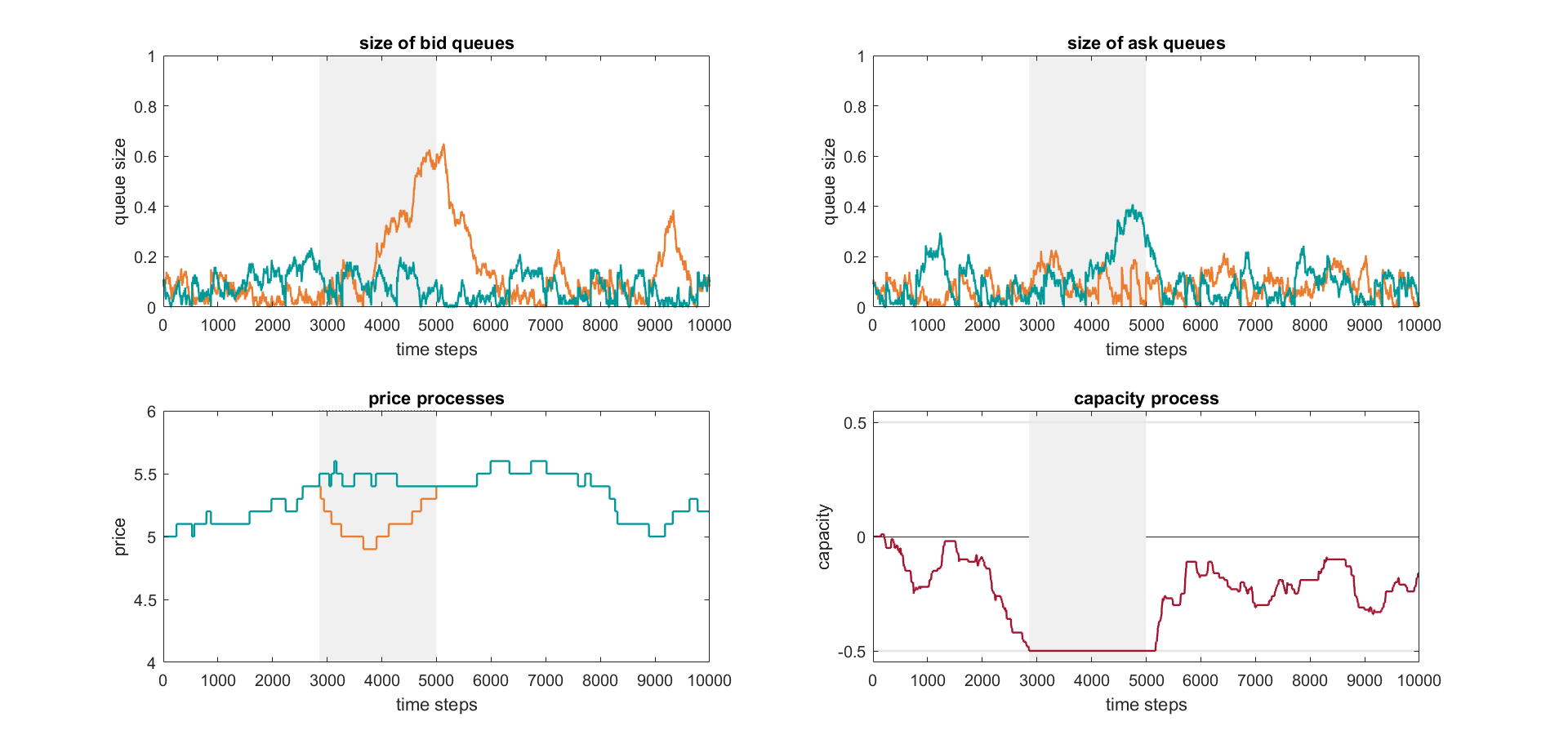}
    \caption{\small Simulation of bid/ask queues and prices of $F$ (orange) and $G$ (turquoise) and of the capacity process in a balanced setting. The white (resp. gray) areas represent the active (resp. inactive) regimes.}
    \label{CB:fig:balOB}
\end{figure}

In a first simulation, we assume that different order types occur equally likely, while the frequency of market orders is slightly higher than the frequency of limit order placements: for $(i,I) \in \{b,a\} \times \{F,G\}$, we choose
\[\Pro\left[(\phi^{(n)}_1, \psi^{(n)}_1) = (i,I)\right] = \frac{1}{4} \quad \text{and} \quad \Pro\left[V^{i,I,(n)}_1 = -\vn \,\right|\left.\, (\phi^{(n)}_1, \psi^{(n)}_1) = (i,I)\right]
= \frac{1}{2} + 5\vn.\]
This choice yields $\mu^{i,I, (n)} = -2.5$ and $(\sigma^{i,I,(n)})^2 = 0.25 - n^{-1}(\mu^{i,I,(n)})^2$ for $(i,I) \in \{b,a\} \times \{F,G\}$ and hence corresponds to a \textit{balanced} market setting. We present one realization of this balanced cross-border market model in Figure \ref{CB:fig:balOB}. For $1,000$ replications, we obtain an empirical probability of around $0.53$ to observe at least one regime switch for the chosen model parameters. Since we have chosen the frequency of market orders higher than the frequency of limit orders, we tend to observe several price changes and a fluctuating capacity process. If we changed this relation in favor of order placements, we would observe less price changes, less variation of the capacity process, and also less regime switches on average. 
\newline

\begin{figure}[h]
	\centering
	 \includegraphics[scale = 0.3]{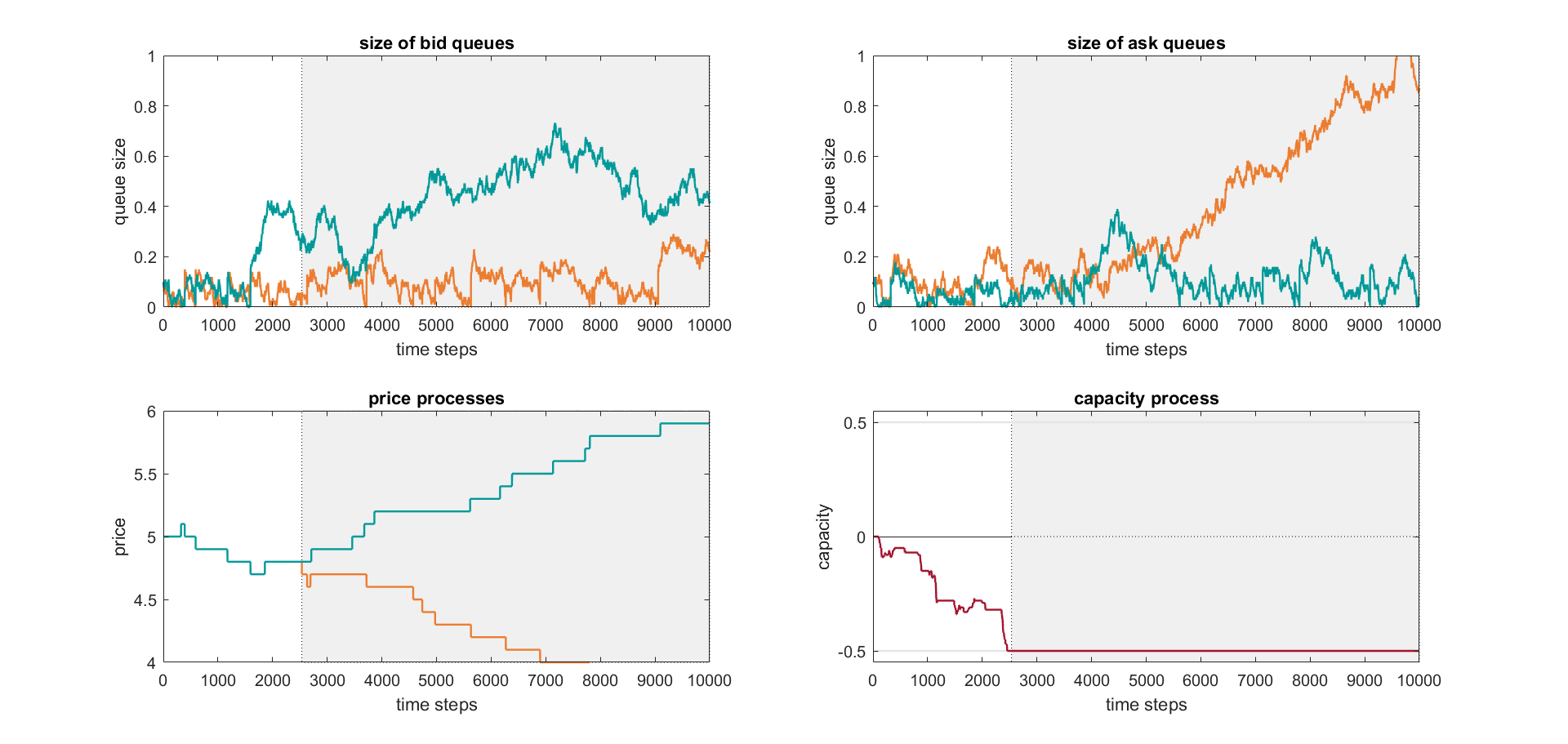}
    \caption{\small Simulation of bid/ask queues and bid prices of $F$ (orange) and $G$ (turquoise) as well as the capacity process in the imbalanced setting. The white (resp. gray) areas represent the active (resp. inactive) regimes.}
    \label{CB:fig:imbOB-1}
\end{figure}

In a second simulation, we choose the frequency of cross-border trades in direction $G$ to $F$ smaller than the frequency of cross-border trades in direction $F$ to $G$: for $(i,I) \in \{b,a\} \times \{F,G\},$ 
\begin{align*}
\Pro[(\phi^{(n)}_1, \psi^{(n)}_1) = (i,I)] &= 0.25,\\
\Pro\left[V^{i,I,(n)}_1 = -\vn \,\right|\,\left. (\phi^{(n)}_1, \psi^{(n)}_1) = (i,I)\right] &= 
\begin{cases}
  0.5 + 5\vn&: (i,I) \in \{(b,F),(a,G)\}\\  0.5 &:(i,I) \in \{(a,F), (b,G)\}.  
\end{cases}
\end{align*}
This choice yields the parameters $\mu^{b,F,(n)} = \mu^{a,G,(n)} = -2.5,$ $\mu^{a,F,(n)} = \mu^{b,G,(n)}=0,$ and $(\sigma^{i,I,(n)})^2 = 0.25 - n^{-1}(\mu^{i,I,(n)})^2$ for all $(i,I)\in \{b,a\} \times \{F,G\}$ and hence corresponds to an \textit{imbalanced} market setting. We present one realization of this imbalanced cross-border market model in Figure \ref{CB:fig:imbOB-1}, where we observe a single regime switch from the active to the inactive regime. Since we have chosen the frequency of market orders corresponding to possible exports from $F$ higher than the frequency of those corresponding to possible imports to $F$, the dynamics of the capacity process has a higher probability to move downward. Since this imbalance is maintained after the regime switch, the best bid price of $F$ decreases whereas the price of $G$ increases.

\subsection{Numerical analysis of price stability}

The economic idea behind a market coupling is to increase liquidity, price stability, and overall trading volume. To study the influence of the market coupling on price stability, we simulate the active and inactive dynamics from the same underlying order flow process and compare the empirical mean number of price changes and the bid price range, i.e.~the difference between the maximum and minimum price attained during the interval $[0,1]$ averaged over $m = 1,000$ replications. At $t=0$ we initialize all four queues at size $10\vn$. If a queue is depleted, we draw its reinitialization value from the uniform distribution on $\{j\vn:\ j=10,\dots,20\}$, while the queue on the opposite side is reinitialized by adding an independent random variable, being uniformly distributed on $\{j\vn:\ j=2,2.2,2.4,\dots,4\}$, to the previous queue size. Order types are assumed to be equally likely, i.e.~$\Pro[(\phi^{(n)}_k, \psi^{(n)}_k) = (i,I)] = 0.25$ for all $k\in\N$ and $(i,I) \in \{b,a\} \times \{F,G\}.$ 
We study four different scenarios for the means of the different order types:
\begin{align*}
    \text{a)} & \quad \mu^{b,F,(n)} = \mu^{a,F,(n)} = \mu^{b,G,(n)} = \mu^{a,G,(n)} = 0,\\
    \text{b)} & \quad \mu^{b,F,(n)} = \mu^{b,G,(n)} = -2.5, \quad \mu^{a,F,(n)} = \mu^{a,G,(n)} = 0,\\
    \text{c)} & \quad \mu^{b,F,(n)} = \mu^{a,F,(n)} = -2.5, \quad \mu^{b,G,(n)} = \mu^{a,G,(n)} = 0,\\
    \text{d)} & \quad \mu^{b,F,(n)} = \mu^{a,G,(n)} = -2.5, \quad \mu^{a,F,(n)} = \mu^{b,G,(n)} = 0.
\end{align*}

In the following, we denote by $\mathcal{N}^{\circ},$ $\mathcal{N}^{F},$ and $\mathcal{N}^{G}$ the mean number of price changes and by $\mathcal{R}^{\circ}$, $\mathcal{R}^F,$ and $\mathcal{R}^G$ the average bid price ranges in the shared order book respectively the two national order books, obtained from simulations of the active respectively inactive order book dynamics from the same order flow. We note that $\mathcal{R}^{\circ},$ $\mathcal{R}^F,$ and $\mathcal{R}^G$ are denoted as multiples of the tick size. 

\begin{table}[H]
\centering
\begin{tabular}{c||c|c|c||c|c|c}
	 & $\mathcal{N}^{\circ}$ & $\mathcal{N}^F$ & $\mathcal{N}^G$ & $\mathcal{R}^{\circ}$ & $\mathcal{R}^F$ & $\mathcal{R}^G$\\
	\hline \hline
	a) & $6.88$ & $11.91$ & $11.86$ & $3.10$ & $4.48$ & $4.44$\\
	b) &  $27.91$ & $34.99$ & $35.36$ &  $18.36$ & $15.03$ & $15.13$\\
	c) & $23.39$ & $50.34$ & $11.72$ & $6.85$ & $10.25$ & $ 4.43$\\
	d) & $23.6$ & $35.19$ & $35.31$ & $6.71$ & $14.93$ & $15.48$
\end{tabular}
\caption{\small Mean number of price changes and average price range (in ticks) for different parameter constellations.}
\label{CB:fig:tableNPC}
\end{table}

As the shared market has higher liquidity than each of the national markets, one would intuitively expect that the market coupling reduces the total number of price changes observed in both countries. This is confirmed in our simulations, cf.~Table \ref{CB:fig:tableNPC}. Except from scenario b) the market coupling also reduces the average price range. While the coupling always leads to an increase of the order volume at the best bid and ask queues, the mean and the variance of the cumulative order flow process differs from those of the individual order flow processes. This difference can amplify or cancel out the effect caused by the increased trading volume. In scenario b) the Sharpe ratio of the cumulative order flow process is higher than the Sharpe ratio of each individual queue due to their independence, which leads to a larger price range.

\section{Analysis of the active market dynamics}\label{CB:sec:anaAD}

In this section we derive a heavy traffic diffusion limit for the active market dynamics $\tilde{S}^{(n)}$ specified in Definition \ref{CB:def:active} and analyze its dynamics. To prove the convergence of $\tilde{S}^{(n)}$, we combine ideas from \cite[Theorem 2]{CL21} and \cite[Theorem 4.14]{RR03}. We first define a function $\tilde{\Psi}$ which allows us to represent $\tilde{S}^{(n)}$ as a regulated process of the net order flow process $X^{(n)}$ and the reinitialization sequences $\tilde{R}^{+,(n)}$ and $\tilde{R}^{-,(n)}$, i.e.
\[\tilde{S}^{(n)}=\tilde{\Psi}\left(\tilde{Q}^{(n)}_0+X^{(n)},\tilde{R}^{+,(n)},\tilde{R}^{-,(n)}\right)+\varepsilon_n,\quad n\in\N,\]
where $\varepsilon_n$ denotes some small error term. Having constructed $\tilde{\Psi}$, we will study its continuity set and derive a convergence result for $\tilde{S}^{(n)}$ with the help of the continuous mapping theorem. 

\subsection{The active market dynamics as a regulated process of the net order flow}\label{CB:sec:regProQ}

The goal of this subsection is to construct the function $\tilde{\Psi}$. 
To ease notation, we denote for all $\omega \in D(\R_+, \R^2)$ the component-wise reflection at zero by
\[\ell^{(2)}_t(\omega) := \left(\sup_{s\leq t}\left(-\pi^{(2)}_1\omega(s)\right)^+,\, \sup_{s\leq t}\left(-\pi^{(2)}_2 \omega(s)\right)^+\right), \quad t\geq0,\]
where $x^+ := \max\{x, 0\}.$ We start by defining a function $g: D(\R_+, \R^2) \rightarrow D(\R_+, \R^2_+)$, which describes the queue size process on each side of the shared order book before the first price change, when applied to the corresponding net order flow process.

\begin{Def}\label{CB:def:g}
Let $\omega \in D(\R_+, \R^2)$. We inductively define $g(\omega)\in D(\R_+, \R^2_+)$ as follows:
\begin{itemize}
	\item Set $g_1(\omega):= \omega$ and $\htau_1 := \htau_1(\omega) := \inf\{t \geq 0: \exists \, i \in \{1,2\} \text{ with } \pi^{(2)}_i \omega(t) \leq 0\}.$
	\item For $k \geq 2$ and $t\geq \htau_{k-1} $ set
	\begin{align*}
	g_k(\omega)(t) := g_{k-1}(\omega)(\htau_{k-1}-) + \omega(t)- \omega(\htau_{k-1}-) + \ell^{(2)}_t\big(g_{k-1}(\omega)(\htau_{k-1}-) + \omega- \omega(\htau_{k-1}-)\big)\mathcal{R}
	\end{align*}
    with reflection matrix $\mathcal{R} := \begin{pmatrix*}[r] 1&-1\\-1&1\end{pmatrix*}$ and
	\begin{align*}
	\htau_k := \htau_k(\omega) := \inf\Big\{t \geq \htau_{k-1} :& \, \pi^{(2)}_i g_k(\omega)(t)\1_{\left\{\pi^{(2)}_i g_{k-1}(\omega)(\htau_{k-1}) > 0\right\}}\leq 0\text{ for }i=1,2\Big\}.
	\end{align*}
	\item Set $\htau_0 := 0,\ \htau_{\infty}:=\htau_{\infty}(\omega) := \lim_{k\rightarrow \infty}\htau_k(\omega),$ and for $t\geq0$,
	\[
 g(\omega)(t):=\sum_{k=1}^\infty g_k(\omega)(t)\1_{[\htau_{k-1},\htau_k)}(t).\]
\end{itemize}
Moreover, we define for $t\geq0$,
	\begin{align*}
 \bar{g}(\omega)(t):=\sum_{k=1}^\infty \ell^{(2)}_{t\wedge\htau_k-}\big(g_{k-1}(\omega)(\htau_{k-1}-) + \omega- \omega(\htau_{k-1}-)\big)\ \1_{[\htau_{k-1},\infty)}(t).
 \end{align*}
\end{Def}

Note that on the event that no price change occurs during $[0,t)$, we have for all $s\in[0,t)$ and $i=b,a$ the representation
\begin{align*}
\pi_i \tilde{Q}(s)=g\left(\pi_i\tilde{Q}_0^{(n)}+\pi_iX^{(n)}\right)(s)\quad\text{ and }\quad\pi_{i} M^{(n)}(s) = \bar{g}\left(\pi_i\tilde{Q}^{(n)}_0 + \pi_iX^{(n)}\right)(s).
\end{align*}

By construction on each interval $[\hat{\tau}_{k-1}, \hat{\tau}_k)$, $k\geq 2$, only one component of the two-dimensional reflection process $\ell^{(2)}(g_{k-1}(\omega)(\hat{\tau}_{k-1}-) + \omega - \omega(\hat{\tau}_{k-1}-))$ increases, while the other stays at zero. Moreover, the increasing component alternates over successive intervals $[\hat{\tau}_{k-1}, \hat{\tau}_k),\ k\geq2$. Hence, if we denote by $\Gamma_1$ the one-dimensional Skorokhod map, i.e.~for any $\psi\in D(\R_+,\R)$,
\[\Gamma_1(\psi)(t) := \psi(t) + \sup_{s\leq t} (-\psi(s))^+,\quad t\geq0,\]
then we have on $[\htau_{k-1},\htau_k)$,
\begin{equation}\label{CB:eq:SP1d}
\begin{split}
     \pi^{(2)}_i g(\omega) &=
        \Gamma_1\left(\pi^{(2)}_i(g_{k-1}(\omega)(\hat{\tau}_{k-1}-) + \omega - \omega(\hat{\tau}_{k-1}-))\right),\\
  \pi^{(2)}_jg(\omega)& = 
         \pi^{(2)}_j(g_{k-1}(\omega)(\hat{\tau}_{k-1}-) + \omega - \omega(\hat{\tau}_{k-1}-))\\
         &\qquad\qquad\qquad- (\Gamma_1-\id)\left(\pi^{(2)}_i(g_{k-1}(\omega)(\hat{\tau}_{k-1}-) + \omega - \omega(\hat{\tau}_{k-1}-)) \right),
         \end{split}
\end{equation}   
where $i$ denotes the increasing component of the reflection process on $[\htau_{k-1},\htau_k)$ and $j\neq i$ the other component. Therefore, on $[0,\htau_\infty)$ the function $g$ is iteratively constructed as a series of local solutions to the one-dimensional Skorokhod problem on alternating axes.

\begin{Lem}\label{CB:res:sumCompg}
Let $\omega\in D(\R_+,\R^2)$. Then:  
\begin{enumerate}
    \item[i)] $g(\omega)=\omega+\bar{g}(\omega)\mathcal{R}$ on $[0,\htau_\infty)$.
    \item[ii)] $g$ is sum-conserving before time $\htau_\infty$, i.e.~$h_1(g(\omega))=h_1(\omega)$ on $[0,\htau_\infty)$.
    \item[iii)] If $\htau_k<\htau_\infty<\infty$ for all $k\in\N$, the map $t\mapsto g(\omega)(t)$ is continuous at $\htau_\infty$. 
    \item[iv)] $\pi_i^{(2)}\bar{g}(\omega)(t)=\int_0^t \1_{\{\pi_i^{(2)}g(\omega)=0\}}d\pi_i^{(2)}\bar{g}(\omega)(s)$ for all $t\in [0,\htau_\infty)$ and $i=1,2$.    
\end{enumerate}
\end{Lem}

The proof of Lemma \ref{CB:res:sumCompg} can be found in Appendix \ref{CB:app:proofs}. Lemma \ref{CB:res:sumCompg} ii) implies that 
$$\htau_\infty=\inf\{t\geq0:\, g(\omega)(t)=(0,0)\}=\inf\{t\geq0:\, h_1(g(\omega))(t)=0\}=\inf\{t\geq0:\, h_1(\omega)(t)\leq0\}.$$
Moreover, items i) and iv) of Lemma \ref{CB:res:sumCompg} imply that for any $\omega\in D(\R_+,\R^2)$ with $\omega(0)\in\R_+^2$ the pair $(g(\omega),g(\omega)-\omega)=(g(\omega),\bar{g}(\omega)\mathcal{R})$ is on $[0,\htau_\infty)$ the unique solution to the so called two-dimensional general processor sharing (extended) Skorokhod problem (GPS (E)SP). The definition of the GPS Skorokhod problem and its extended version can be found in \cite[Section 3.2]{R06}. Especially, it is shown in \cite[Theorem 3.6]{R06} that there always exists a unique solution to the GPS extended Skorokhod problem on $[0,\infty)$ and that the solution map is Lipschitz continuous. Moreover, according to \cite[Theorem 1.3]{R06} the GPS SP and the GPS ESP agree on $[0,\htau_\infty)$.

\begin{Rmk}
    If $W$ is a two-dimensional linear Brownian motion with positive definite covariance matrix and $\Pro[W(0)=(0,0)]=0$, then $\Pro[\htau_k=\htau_\infty<\infty]=0$ for all $k\in\N$.
    Clearly, this is true for $k=1$ as planar Brownian motion does not hit any given point except the starting point almost surely. Now suppose the claim is true for some $k\in\N$. Then on the event $\{\htau_k(W)<\infty\}$, $g(W)$ behaves by construction on $[\htau_k,\htau_{k+1}]$ like an obliquely reflected linear Brownian motion on a half-space started outside the origin and hence does not hit the origin almost surely, cf.~\cite[Theorem 2]{R86}.   
    Therefore also $\Pro[\htau_{k+1}=\htau_\infty<\infty\, |\, \htau_k<\infty]=0$ and the claim follows by induction.
\end{Rmk} 

Next we define $G : D(\R_+, \R^4) \rightarrow D(\R_+, \R^4_+)$ and $\bar{G}: D(\R_+, \R^4) \rightarrow D(\R_+, \R)$ by
\begin{align}\label{CB:def:G}
\begin{aligned}
	G(\omega) &:= \Big(\pi^{(2)}_1g(\pi_{b}\omega),\, \pi^{(2)}_1 g(\pi_{a}\omega),\, \pi^{(2)}_2 g(\pi_{b}\omega), \,\pi^{(2)}_2g(\pi_{a}\omega) \Big),\\
\bar{G}\left(\omega\right) &:= \left(\pi_2^{(2)}\overline{g}(\pi_b\omega) +  \pi_1^{(2)}\overline{g}(\pi_a\omega) \right) - \left(\pi_1^{(2)}\overline{g}(\pi_b\omega)+ \pi_2^{(2)}\overline{g}(\pi_a\omega)\right). 
\end{aligned}
\end{align}
On the event that no price change occurs during $[0,t)$, the function $G$ (resp.~$\bar{G}$) applied to the the net order flow process $X^{(n)}$ describes the dynamics of the queue size process $\tilde{Q}^{(n)}$ (resp.~the capacity process $\tilde{C}^{(n)}$) of the shared order book, i.e.~for all $s\in[0,t)$,
\begin{align*}
    \tilde{Q}^{(n)}(s) = G\left(\tilde{Q}^{(n)}_0 + X^{(n)}\right)(s) \quad\text{ and }\quad \tilde{C}^{(n)}(s) = \tilde{C}^{(n)}_0 + \bar{G}\left(\tilde{Q}^{(n)}_0 + X^{(n)}\right)(s).
\end{align*}
Moreover, let us introduce the first hitting time maps $\tilde{\tau}_b,\ \tilde{\tau}_a,\ \tilde{\tau}: D(\R_+, \R^4) \rightarrow \R_+$ given by
\begin{align}
    \begin{aligned}\label{CB:def:tildetauab}
    \tilde{\tau}_b(\omega)&:= \inf\left\{t\geq 0: \pi^{(2)}_1 h(\omega)(t) \leq 0\right\},\\
    \tilde{\tau}_a(\omega)&:= \inf\left\{t\geq 0: \pi^{(2)}_2 h(\omega)(t) \leq 0\right\},\\
    \tilde{\tau}(\omega)&:= \tilde{\tau}_b(\omega) \wedge \tilde{\tau}_a(\omega).
    \end{aligned}
\end{align} 
With the help of $\tau,\ G$, and $\bar{G}$ we are now able to define functions $\tilde{\Psi}^Q$, $\tilde{\Psi}^C$, and $\tilde{\Psi}^B$, which allow to construct the queue size process $\tilde{Q}^{(n)}$, the capacity process $\tilde{C}^{(n)}$, and the bid price process $\tilde{B}^{(n)}$ from the net order flow process $X^{(n)}$ and the sequences of reinitialization values $\tilde{R}^{+,(n)}$ and $\tilde{R}^{-,(n)}$. 

\begin{Def}\label{CB:def:PsiQ1}
Let $\omega \in D(\R_+, \R^4)$ and let $r^+ = (r^+_n)_{n\geq 1}, r^- = (r^-_n)_{n\geq 1} \in ((0,\infty)^4)^{\N}$.  For $k\in \N_0,$ we define functions $\tilde{\Psi}^Q_k(\omega, r^+, r^-),\ k\in\N_0,$ and $\tilde{\Psi}^Q(\omega, r^+, r^-)$ as follows:
\begin{itemize}
    \item Set $\tilde{\Psi}^Q_0(\omega, r^+, r^-) := G(\omega).$
    \item Let $k\geq 1$ and $\tilde{\Psi}^Q_{k-1} := \tilde{\Psi}^Q_{k-1}(\omega, r^+, r^-).$ Define iteratively\footnote{$\tilde{\Psi}_k(\omega,r^+,r^-)$ depends on $r^+,r^-$ only through the first $k$ reinitialization values. Therefore, with a slight abuse of notation, we may also write $\tilde{\Psi}^Q_k(\omega,r^+,r^-)=\tilde{\Psi}^Q_k(\omega,\tilde{\Psi}^Q_k(\tilde{\tau}_1),\dots,\tilde{\Psi}^Q_k(\tilde{\tau}_k))$.}
    \begin{align*}
        \tilde{\Psi}^Q_k(\omega,r^+, r^-) &:= \tilde{\Psi}^Q_{k-1}\1_{\left[0, \tilde{\tau}(\tilde{\Psi}^Q_{k-1})\right)}\\
        &\qquad + \1_{\left[ \tilde{\tau}(\tilde{\Psi}^Q_{k-1}),\infty\right)}\Bigg\{\1_{\left\{\tilde{\tau}(\tilde{\Psi}^Q_{k-1}) = \tilde{\tau}_a(\tilde{\Psi}^Q_{k-1})\right\}} G\left(r^+_k + \omega - \omega\left(\tilde{\tau}(\tilde{\Psi}^Q_{k-1})\right)\right)\\
        &\hspace{3cm} +\1_{\left\{\tilde{\tau}(\tilde{\Psi}^Q_{k-1}) \neq \tilde{\tau}_a(\tilde{\Psi}^Q_{k-1})\right\}} G\left(r^-_k + \omega - \omega\left(\tilde{\tau}(\tilde{\Psi}^Q_{k-1})\right)\right)\Bigg\}.
    \end{align*}
    \item Finally, set $\tilde{\tau}_0 := 0,$ $\tilde{\tau}_k := \tilde{\tau}\big(\tilde{\Psi}^Q_{k-1}(\omega, r^+, r^-)\big)$ for $k \geq 1,$ and
    \[\tilde{\Psi}^Q(\omega, r^+, r^-)(t) = \sum_{k=1}^{\infty} \tilde{\Psi}^Q_{k-1}(\omega, r^+, r^-)(t) \1_{\left[\tilde{\tau}_{k-1}, \tilde{\tau}_k\right)}(t), \quad  t \geq0.\]
\end{itemize}
Moreover, we define
\footnote{Since $r^+, r^- \in ((0,\infty)^4)^{\N},$ we have $\bar{G}(\tilde{\Psi}^Q(\omega, r^+,r^-)(\tilde{\tau}_k) + \omega - \omega(\tilde{\tau}_k))(\tilde{\tau}_k) = 0$ for all $k\in \N$ and hence $\tilde{\Psi}^C(\omega, r^+, r^-)(\tilde{\tau}_k) = \tilde{\Psi}^C(\omega, r^+, r^-)(\tilde{\tau}_k+),$ i.e.~$\tilde{\Psi}^C(\omega, r^+, r^-)$ is right-continuous.}
	\[\tilde{\Psi}^{C}(\omega, r^+, r^-)(t) = \sum_{k=0}^\infty \bar{G}\Big(\tilde{\Psi}^Q(\omega, r^+, r^-)(\tilde{\tau}_k) + \omega - \omega(\tilde{\tau}_k)\Big)(t\wedge\tilde{\tau}_{k+1})\ \1_{\left(\tilde{\tau_k},\infty\right)}(t),\quad t\geq0,\]
as well as 
	\[N_a(\omega, r^+, r^-)(t) := \sum_{k:\ \tilde{\tau}_k \leq t} \1_{\left\{\pi^{(2)}_2 \left(h \circ \tilde{\Psi}^Q_{k-1}\right)(\tilde{\tau}_k) \leq 0\right\}},\hspace{0.2cm} N_b(\omega, r^+, r^-)(t) := \sum_{k:\ \tilde{\tau}_k \leq t} \1_{\left\{\pi^{(2)}_1 \left(h \circ \tilde{\Psi}^Q_{k-1}\right)(\tilde{\tau}_k) \leq 0\right\}},\]
	and
	\[\tilde{\Psi}^B(\omega, r^+, r^-)(t) := \left(N_a(\omega, r^+, r^-)(t) - N_b(\omega, r^+, r^-)(t)\right) (\delta,\delta),\quad t\geq0.\]
\end{Def}

According to the above definition, $\tilde{\Psi}^{Q}(\omega, r^+, r^-)$ is obtained by ``regulating'' $\omega \in D(\R_+, \R^4)$ according to the function $G$ and the sequences $r^+,r^-$:  between two hitting times $\tilde{\tau}_k$ and $\tilde{\tau}_{k+1}$, the function $\tilde{\Psi}^{Q}(\omega, r^+, r^-)$ behaves as $G(r^+_{k} + \omega - \omega(\tilde{\tau}_k))$ or $G(r^-_k + \omega - \omega(\tilde{\tau}_k))$ depending on which component of $(h \circ \tilde{\Psi}^Q_{k-1})(\omega, r^+, r^-)$ first hits zero. Moreover, at times $\tilde{\tau}_1, \tilde{\tau}_2, \cdots$, the process jumps to a new position inside $\R^4_+$ taken from the sequence $(r^+_n)_{n\geq 1}$ or $(r^-_n)_{n\geq 1}$. Note that $\tilde{\Psi}^{Q}(\omega, r^+, r^-)$ is right-continuous, but may not have a left-limit at $\tilde{\tau}_\infty:=\lim_{k\rightarrow\infty}\tilde{\tau}_k$.

\begin{The}\label{CB:res:repactive}
For each $n \in \N$ and $t\geq0$,
\begin{align*}\label{CB:eq:constM}
\tilde{Q}^{(n)}(t)  &= \tilde{\Psi}^Q\left(\tilde{Q}^{(n)}_0 + X^{(n)}, \tilde{R}^{+,(n)}, \tilde{R}^{-,(n)}\right)(t),\\
\tilde{B}^{(n)}(t)  &= \tilde{B}^{(n)}_0 + \tilde{\Psi}^B\left(\tilde{Q}^{(n)}_0 + X^{(n)}, \tilde{R}^{+,(n)}, \tilde{R}^{-,(n)}\right)(t),\\
\tilde{C}^{(n)}(t) &= \tilde{C}^{(n)}_0 + \tilde{\Psi}^{C}\left(\tilde{Q}^{(n)}_0 + X^{(n)}, \tilde{R}^{+,(n)}, \tilde{R}^{-,(n)}\right)(t)+\vn c_n(t),
\end{align*}
where
\begin{itemize}
\setlength\itemsep{0em}
	\item $(\tilde{B}^{(n)}_0,\tilde{Q}^{(n)}_0, \tilde{C}^{(n)}_0) \in  (\delta,\delta) \Z \times (\vn \N)^4 \times \vn \Z$ are the initial values,
	\item $X^{(n)}$ is the piecewise constant interpolation of the net order flow process, 
	\item $\tilde{R}^{+,(n)} \in (\vn(\N+1)^4)^{\N}$ is the sequence of reinitialization values after price increases in $\tilde{S}^{(n)}$,
	\item $\tilde{R}^{-,(n)} \in (\vn(\N+1)^4)^{\N}$ is the sequence of reinitialization values after price decreases in $\tilde{S}^{(n)}$.
\end{itemize}
The $k$-th price change in $\tilde{S}^{(n)}$ occurs at the stopping time
\begin{align*}
   \tilde{\tau}^{(n)}_k := \tilde{\tau}\left(\tilde{\Psi}^{Q}_{k-1}\left(\tilde{Q}^{(n)}_0 + X^{(n)}, \tilde{R}^{+,(n)}, \tilde{R}^{-,(n)}\right)\right) \quad \text{for } k\geq 1,
\end{align*} 
and the error term satisfies
\[|c_n(t)|\leq\max\left\{k: \tilde{\tau}^{(n)}_k\leq t\right\},\qquad t\geq0.\]
\end{The}

\begin{proof}
The result follows from a careful, but straight-forward inspection of the previous definitions. The only thing that needs explanation is the occurrence of the term $\vn c_n$: this is due to the fact that in Definition \ref{CB:def:g} the function $\overline{g}(\omega)$ does not count the possibly last reflection at $\htau_\infty$ (corresponding to the last cross-border trade at the time of a price change) if $\htau_\infty=\htau_k$ for some $k\in\N$.
\end{proof}

Finally, we set $\tilde{\tau}^{(n)}_0 := 0$ and define for $(\omega, r^+, r^-) \in D(\R_+, \R^4) \times ((0,\infty)^4)^{\N} \times ((0,\infty)^4)^{\N}$ and $s_0 = (b_0, q_0, c_0) \in \R^2 \times (0,\infty)^4 \times \R$,
\begin{equation*}
\tilde{\Psi}(s_0, \omega, r^+, r^-):= \left(b_0 + \tilde{\Psi}^B(q_0 + \omega, r^+, r^-), \,\tilde{\Psi}^Q(q_0 + \omega, r^+, r^-),\, c_0 + \tilde{\Psi}^C(q_0+\omega, r^+, r^-)\right).
\end{equation*}
By Theorem \ref{CB:res:repactive} we have
$$\tilde{S}^{(n)}= \left(\tilde{B}^{(n)}, \tilde{Q}^{(n)}, \tilde{C}^{(n)}\right) = \tilde{\Psi}\left(\tilde{S}^{(n)}_0, X^{(n)}, R^{+,(n)}, R^{-,(n)}\right)+\left(0,0,\vn c_n\right).$$

\subsection{Heavy traffic approximation of the active dynamics}\label{CB:sec:HTAppr}

In this subsection, we derive a limit theorem for the active dynamics $\tilde{S}^{(n)} = (\tilde{B}^{(n)}, \tilde{Q}^{(n)}, \tilde{C}^{(n)})$ under the heavy traffic condition that $\vn,\tn\rightarrow0$. Proposition \ref{CB:res:fCLT} and Theorem \ref{CB:res:repactive} strongly suggest to use the continuous mapping theorem to prove the convergence of $\tilde{S}^{(n)}$, for which we have to determine the continuity set of $\tilde{\Psi}$. For $\ell\in\N$, we endow the space $D(\R_+, \R^\ell)$ with the Skorokhod $J_1$-topology,  the set $(\R^\ell_+)^{\N}$ with the topology induced by the cylindrical semi-norm given through 
\[(R^n_k)_{k\in\N}\rightarrow (R_k)_{k\in\N}\quad\Leftrightarrow\quad \forall k\in\N:\ R_k^n\rightarrow R_k,\]
and the space $D(\R_+, \R^\ell) \times (\R^\ell_+)^{\N} \times (\R^\ell_+)^{\N}$ with the corresponding product topology.

\begin{The}[Continuity of $\tilde{\Psi}$]\label{CB:res:contPsi}
 Let $(\omega_0, r^+, r^-) \in D(\R_+, \R^4) \times (\R^4_+)^{\N} \times (\R^4_+)^{\N}$ satisfy the following conditions:
	\begin{enumerate}
		\item[i)] $\tilde{\tau}_k=\tilde{\tau}\big(\tilde{\Psi}^Q_{k-1}(\omega_0,r^+,r^-)\big)\rightarrow\infty$ as $k\rightarrow\infty$. 
		\item[ii)] For all $k\in\N_0$, $h\big(\tilde{\Psi}^Q(\omega_0,r^+,r^-)(\tilde{\tau}_k)+\omega_0(\cdot+\tilde{\tau}_k)-\omega_0(\tilde{\tau}_k)\big)\in C'_0(\R_+,\R^2)$. Here $C'_0(\R_+,\R^2)$ contains all continuous functions whose components cross $0$ as soon as they hit it, cf.~\eqref{CB:def:C'}.  
		\item[iii)] For all $k\in\N$, $\tilde{\tau}_a\big(\tilde{\Psi}^Q_{k-1}(\omega_0,r^+,r^-)\big)\neq\tilde{\tau}_b\big(\tilde{\Psi}^Q_{k-1}(\omega_0,r^+,r^-)\big)$.  
		\end{enumerate}
	Then $\tilde{\Psi}^Q(\omega_0,r^+,r^-)\in D(\R_+,\R_+^4),\ \tilde{\Psi}^B(\omega_0,r^+,r^-)\in D(\R_+,\R^2),\ \tilde{\Psi}^C(\omega_0,r^+,r^-)\in D(\R_+,\R)$ and the functions $\tilde{\Psi}^{Q},\ \tilde{\Psi}^{B},\ \tilde{\Psi}^{C}$ defined on $D(\R_+,\R^4)\times(\R_+^4)^\N\times(\R_+^4)^\N$ are continuous at $(\omega_0, r^+, r^-)$. Hence, $\tilde{\Psi}$ is continuous at $(s_0,\omega_0,r^+,r^-)$ for any $s_0\in\R^2\times\R_+^4\times\R$.
\end{The}

The proof of Theorem \ref{CB:res:contPsi} can be found in Section \ref{CB:app:contProp} of the Appendix. 
Condition i) ensures that the function $\tilde{\Psi}^{Q}(\omega_0, r^+, r^-)$ is reinitialized only a finite number of times on every compact interval. Condition ii) ensures that the times of price changes are continuous at $(\omega_0,r^+,r^-)$. Finally, condition iii) guarantees that the reinitialization of  $\tilde{\Psi}^{Q}(\omega_0, r^+, r^-)$ at $\tilde{\tau}_k=\tilde{\tau}\big(\tilde{\Psi}_{k-1}^Q(\omega_0,r^+,r^-)\big)$  by either $r^+_k$ or $r^-_k$ is continuous at $(\omega_0, r^+, r^-).$ \par

We will apply the map $\tilde{\Psi}$ to an $\R^2\times(0,\infty)^4\times\R$-valued random variable $\tilde{S}^{(n)}_0$, an independent four-dimensional linear Brownian motion $Y$, and two sequences $\tilde{R}^+$ and $\tilde{R}^-$ of $(0,\infty)^4$-valued random variables to be specified below. In Step 2 of the proof of Theorem \ref{CB:res:convQM-1} we will show that $(\tilde{S}_0,X, \tilde{R}^+, \tilde{R}^-)$ fulfills conditions i)-iii) of Theorem \ref{CB:res:contPsi} with probability one and hence lies in the continuity set of $\tilde{\Psi}$ $\Pro$-a.s.

\begin{The}[Limit theorem for the active dynamics $\tilde{S}^{(n)}$]\label{CB:res:convQM-1}
Let Assumptions \ref{CB:ass:prob}, \ref{CB:ass:scaling}, and \ref{CB:ass:R} be satisfied and suppose that for every $n\in\N$, $\tilde{S}^{(n)}_0=(\tilde{B}^{(n)}_0,\tilde{Q}^{(n)}_0,\tilde{C}^{(n)}_0)$ is a $(\delta,\delta) \Z \times \vn\N^4 \times \vn\Z$-valued random variable on $(\Omega^{(n)},\F,\Pro^{(n)})$ such that $(\tilde{S}^{(n)}_0,X^{(n)})\Rightarrow (\tilde{S}_0,X)$, where $X$ is as in Theorem \ref{CB:res:fCLT} and $\tilde{S}_0=(\tilde{B}_0,\tilde{Q}_0,\tilde{C}_0)$ with $\tilde{Q}_0\in(0,\infty)^4$. Then
\begin{align*}
	\tilde{S}^{(n)} \Rightarrow \tilde{S}:=(\tilde{B}, \tilde{Q}, \tilde{C}):=\tilde{\Psi}\left(\tilde{S}_0, X, \tilde{R}^{+}, \tilde{R}^{-}\right)
\end{align*}
in the Skorokhod topology on $D(\R_+, E),$ 
\begin{align*}
\tilde{R}^+_k := \Phi\left(\tilde{Q}(\tilde{\tau}^*_k-), \epsilon^+_k\right), \quad \tilde{R}^-_k := \Phi\left(\tilde{Q}(\tilde{\tau}^*_k-), \epsilon^-_k\right),\quad k\in\N,
\end{align*}
for independent sequences of i.i.d. random variables $(\epsilon^+_k)_{k\geq 1}$ and $(\epsilon^-_k)_{k\geq 1}$ with $\epsilon^+_k \sim f^+$ and $\epsilon^-_k \sim f^-$, and
\begin{align*}
    \tilde{\tau}^*_0 := 0,\quad \tilde{\tau}^*_k: =
    \inf\left\{t>\tilde{\tau}_{k-1}^*:\ \pi_b\tilde{Q}(t-)=(0,0)\ \text{ or }\  \pi_a\tilde{Q}(t-)=(0,0)\right\} \quad \text{for } k \geq 1.
\end{align*}
Especially, there are only finitely many price changes on every compact interval, i.e.~$\lim_{k\rightarrow\infty}\tilde{\tau}_k^*=\infty$ $\Pro$-a.s.
\end{The}

\begin{proof}
By Skorokhod's representation theorem we may assume that all $S^{(n)}_0,X^{(n)},\epsilon^{+,(n)},\epsilon^{-,(n)}$ as well as $\tilde{S}_0$ and $X$ are defined on a common probability space $(\Omega, \mathcal{F}, \Pro)$ supporting two independent sequences of i.i.d. random variables $(\epsilon^+_k)_{k\geq 1}$, $(\epsilon^-_k)_{k\geq 1}$ with $\epsilon^+_k \sim f^+,$ $\epsilon^-_k \sim f^-$ such that
\[\Pro\left[ \tilde{S}^{(n)}_0\rightarrow\tilde{S}_0,\ X^{(n)}\rightarrow X \text{ and for all } k \geq 1,\, \epsilon^{+,(n)}_k \rightarrow \epsilon^+_k,\, \epsilon^{-,(n)}_k \rightarrow \epsilon^-_k\right] = 1.\]
In the following we denote $Y^{(n)}:=\tilde{Q}_0^{(n)}+X^{(n)}$ and $Y:=\tilde{Q}_0+X$. 
\vskip6pt

\noindent
\underline{Step 1:} Iterative construction of $\tilde{Q},\tilde{R}^+$, and $\tilde{R}^-$.\\
Since for $i=b,a$ the process $\pi_i Y$ defines a planar Brownian motion, it satisfies the assumptions of Lemma \ref{CB:res:contg} $\Pro$-a.s.~and we may conclude that $G$ is $\Pro$-a.s.~surely continuous at $Y.$ Hence, the continuous mapping theorem implies that $G(Y^{(n)}) \rightarrow G(Y)$ $\Pro$-a.s.~in the Skorokhod topology. Next we construct the process $\tilde{Q}$ by induction: 
let $\tilde{\tau}^*_1 := \tilde{\tau}(G(Y))$ be the first time that $(h \circ G)(Y)$ hits the axes and set $\tilde{Q}(t) = G(Y)(t)$ for $t< \tilde{\tau}^*_1.$ Since $(h \circ G)(Y) = h(Y)$ on $[0, \tilde{\tau}^*_1)$ by Lemma \ref{CB:res:sumCompg} ii), $\tilde{\tau}^*_1=\tilde{\tau}(Y)$, i.e.~$\tilde{\tau}^*_1$ equals the first hitting time of $h(Y)$ of the axes. As $h(Y)$ is a planar Brownian motion, $h(Y)\in C_0'(\R_+,\R^2)$ $\Pro$-a.s. Hence, Lemma \ref{CB:res:cont-firstExit-lastValue} and the continuous mapping theorem imply that $\tilde{\tau}^{(n)}_1\rightarrow\tilde{\tau}^*_1$ $\Pro$-a.s. As $G(Y)$ is continuous, the convergence of $G(Y^{(n)})$ towards $G(Y)$ is actually uniformly on compact intervals. 
Therefore, we have on $\{\tilde{\tau}^*_1\geq\tilde{\tau}^{(n)}_1\}$ that $\Pro$-almost surely
\begin{align*}
\left|\tilde{Q}^{(n)}(\tilde{\tau}^{(n)}_1-)-\tilde{Q}(\tilde{\tau}^*_1-)\right|\leq
\left|G(Y)(\tilde{\tau}^{(n)}_1-)-G(Y)(\tilde{\tau}^*_1-)\right|+
\left|G\big(Y^{(n)}\big)(\tilde{\tau}^{(n)}_1-)-G\big(Y\big)(\tilde{\tau}^{(n)}_1-)\right|
\rightarrow0.
\end{align*}
Similarly, one can show that the convergence also holds on the event $\{\tilde{\tau}^*_1<\tilde{\tau}_1^{(n)}\}$. Hence, $\tilde{Q}^{(n)}(\tilde{\tau}^{(n)}_1-)\rightarrow\tilde{Q}(\tilde{\tau}_1^*-)$ $\Pro$-a.s. 
Let us set
\[\tilde{Q}(\tilde{\tau}^*_1) := \Phi\big(\tilde{Q}(\tilde{\tau}^*_1-), \epsilon^+_1\big) \1_{\left\{\tilde{\tau}_{a,1} = \tilde{\tau}^*_1\right\}} + \Phi\big(\tilde{Q}(\tilde{\tau}^*_1-), \epsilon^-_1\big) \1_{\left\{\tilde{\tau}_{b,1} =\tilde{\tau}^*_1\right\}},\]
where
\begin{equation}\label{CB:eq:hitab}
\tilde{\tau}_{a,1} := \tilde{\tau}_a(G(Y))=\tilde{\tau}_a(Y) \quad \text{and} \quad \tilde{\tau}_{b,1} := \tilde{\tau}_b(G(Y))=\tilde{\tau}_b(Y).
\end{equation}
As $h(Y)$ is a planar Brownian motion started outside the origin, it does not hit the origin $\Pro$-a.s. Thus $\tilde{\tau}_{b,1}\neq\tilde{\tau}_{a,1}$ $\Pro$-a.s.~and it follows from Lemma \ref{CB:res:contH}, Assumption \ref{CB:ass:R}, and the continuous mapping theorem that
\[\tilde{Q}^{(n)}(\tilde{\tau}^{(n)}_1)\rightarrow \tilde{Q}(\tilde{\tau}^*_1) \quad \Pro\text{-a.s.}\]
Now assume that we have already defined $\tilde{Q}$ on $[0, \tilde{\tau}^*_k]$ and that
\begin{align*}
\Big(\tilde{\tau}^{(n)}_1, \cdots, \tilde{\tau}^{(n)}_{k}, \tilde{Q}^{(n)}(\tilde{\tau}^{(n)}_1), \cdots, \tilde{Q}^{(n)}(\tilde{\tau}^{(n)}_{k})\Big)\rightarrow \left(\tilde{\tau}^*_1, \cdots, \tilde{\tau}^*_{k}, \tilde{Q}(\tilde{\tau}^*_1), \cdots, \tilde{Q}(\tilde{\tau}^*_{k})\right) \quad \Pro\text{-a.s.}
\end{align*}
Define the stopping time $\tilde{\tau}^*_{k+1} :=\tilde{\tau}\big(\tilde{\Psi}^{Q}_{k}(Y, \tilde{Q}(\tilde{\tau}^*_1), \cdots, \tilde{Q}(\tilde{\tau}^*_{k}))\big)$ and set 
\begin{align*}
\tilde{Q}(t) &:= \tilde{\Psi}^{Q}_{k}(Y, \tilde{Q}(\tilde{\tau}^*_1), \cdots, \tilde{Q}(\tilde{\tau}^*_{k}))(t) \quad \text{for } t < \tilde{\tau}^*_{k+1},\\
\tilde{Q}(\tilde{\tau}^*_{k+1}) &:= \Phi\left(\tilde{Q}(\tilde{\tau}^*_{k+1}-), \epsilon^+_{k+1}\right) \1_{\left\{\tilde{\tau}_{a,k+1} = \tilde{\tau}^*_{k+1}\right\}}+ \Phi\left(\tilde{Q}(\tilde{\tau}^*_{k+1}-), \epsilon^-_{k+1}\right) \1_{\left\{\tilde{\tau}_{b,k+1} = \tilde{\tau}^*_{k+1}\right\}},
\end{align*}
where
\[\tilde{\tau}_{a,k+1} := \tilde{\tau}_a\left(\tilde{\Psi}^{Q}_{k}(Y, \tilde{Q}(\tilde{\tau}^*_1), \cdots, \tilde{Q}(\tilde{\tau}^*_{k}))\right), \quad \tilde{\tau}_{b,k+1} := \tilde{\tau}_b\left(\tilde{\Psi}^{Q}_{k}(Y, \tilde{Q}(\tilde{\tau}^*_1), \cdots, \tilde{Q}(\tilde{\tau}^*_{k}))\right).\]
Note that $\tilde{W}_k := \tilde{Q}(\tilde{\tau}^*_k) + Y(\cdot + \tilde{\tau}^*_k) - Y(\tilde{\tau}^*_k)$ is a four-dimensional linear Brownian motion with positive definite covariance matrix. Hence, $h(\tilde{W}_k)$ is a two-dimensional linear Brownian motion with positive definite covariance matrix and the conditions of Lemma \ref{CB:res:contPsik} are fulfilled $\Pro$-a.s. Therefore, we obtain 
 by similar arguments as above that
\begin{align*}
\tilde{\Psi}^{Q}_k\left(Y^{(n)}, \tilde{Q}^{(n)}(\tilde{\tau}^{(n)}_1), \cdots, \tilde{Q}^{(n)}(\tilde{\tau}^{(n)}_k)\right) &\rightarrow \tilde{\Psi}^{Q}_k\left(Y, \tilde{Q}(\tilde{\tau}^*_1), \cdots, \tilde{Q}(\tilde{\tau}^*_k)\right) \quad \Pro\text{-a.s.}\\
\left(\tilde{\tau}^{(n)}_{k+1}, \tilde{Q}^{(n)}(\tilde{\tau}^{(n)}_{k+1})\right) &\rightarrow \left(\tilde{\tau}^*_{k+1},\tilde{Q}^{(n)}(\tilde{\tau}_{k+1})\right) \quad \Pro\text{-a.s.}
\end{align*}
In the following, let us denote $\tilde{R}^+_k := \Phi(\tilde{Q}(\tilde{\tau}^*_k-), \epsilon^+_k)$ and $\tilde{R}^-_k := \Phi(\tilde{Q}(\tilde{\tau}^*_k-), \epsilon^-_k)$ for all $k \geq 1.$ Then we showed that $(\tilde{R}^{+,(n)}, \tilde{R}^{-,(n)}) \rightarrow (\tilde{R}^+, \tilde{R}^-)$ $\Pro$-a.s.
\vskip6pt

\noindent
\underline{Step 2:} $(Y, \tilde{R}^+, \tilde{R}^-)$ lies in the continuity set of $\tilde{\Psi}$ $\Pro$-a.s.\\
We have to verify conditions i)-iii) of Theorem \ref{CB:res:contPsi}. As noted above, conditions ii) and iii) are satisfied $\Pro$-a.s.~because every $h(\tilde{W}_k)$ is a two-dimensional linear Brownian motions with positive definite covariance matrix. 
To verify condition i), let us show that there are only finitely many price changes on every compact interval. For this note that the distribution of the first hitting time of a planar Brownian motion is well-known, cf.~\eqref{CB:def:condDist-Duration}. Hence, for any $T  > 0$, $y \in (0,\infty)^4$, and  $k \geq 0$, we know that
\begin{equation*}\label{CB:eq:hittingPBM}
\Pro\left[\left.\tilde{\tau}(\tilde{W}_k) \leq T \,\right.|\, \tilde{W}_k(0) = y\right] =: \beta(T,y) \in (0,1)
\end{equation*}
is decreasing in $y$ component-wise for fixed $T$. For $\varepsilon > 0$ define $U(\varepsilon):= \{x\in \R^4_+: \pi_j x \leq \varepsilon\  \forall\,  j \}$ and $\epsilon_k := \epsilon_k^+\1_{\{\tilde{\tau}_{a,k} = \tilde{\tau}^*_{k}\}} + \epsilon^-_k\1_{\{\tilde{\tau}_{b,k}=\tilde{\tau}^*_k\}},\ k\in\N$. For any $\varepsilon>0$ and $m\geq0$ we have
by Assumption \ref{CB:ass:R}, 
\begin{align*}
   	\Pro\left[\left.\tilde{\tau}(\tilde{W}_m) \leq T\ \right|(\tilde{S}(u))_{u\leq\tilde{\tau}^*_m}\right]
	&\leq\Pro\left[\left.\tilde{\tau}(\tilde{W}_m) \leq T \,\right|\, \tilde{W}_m(0)\right] \1_{\{\tilde{W}_m(0) \notin U(\varepsilon)\}} + \1_{\{\tilde{W}_m(0) \in U(\varepsilon)\}}\\
	&\leq\beta\left(T,(\varepsilon,\varepsilon,\varepsilon,\varepsilon)\right)  \1_{\{\tilde{W}_m(0) \notin U(\varepsilon)\}} + \1_{\{\tilde{W}_m(0) \in U(\varepsilon)\}}\\
	&\leq \left(1-\beta\left(T,(\varepsilon,\varepsilon,\varepsilon,\varepsilon)\right) \right)\1_{\{\alpha \epsilon_m \in U(\varepsilon)\}}+\beta\left(T,(\varepsilon,\varepsilon,\varepsilon,\varepsilon)\right)=:\beta_m(T,\varepsilon),
\end{align*}
where $\beta_m(T,\varepsilon)$ is independent of $\{\tilde{\tau}(\tilde{W}_k)\leq T\ \forall\ k\leq m-1\}$. As $\Pro[\alpha\epsilon_m\in U(\varepsilon)]$ does not depend on $m$ and is strictly smaller than $1$ for $\varepsilon$ small enough by Assumption \ref{CB:ass:R} ii), we may conclude that
\begin{equation*}
\begin{split}
	\Pro\left[\tilde{\tau}(\tilde{W}_k) \leq T \,\, \forall \, k\leq m\right] &= \E\left[\prod_{k\leq m-1}  \1_{\left\{\tilde{\tau}(\tilde{W}_k) \leq T\right\}} \cdot\Pro\left[\left.\tilde{\tau}(\tilde{W}_m)\leq T \ \right|\, (\tilde{S}(u))_{u\leq\tilde{\tau}^*_m} \right]\right]\\
	&\leq \Pro\left[\tilde{\tau}(\tilde{W}_k) \leq T \,\, \forall \, k\leq m-1\right]\cdot \E\big[\beta_m(T,\varepsilon)\big]\leq \left(\E\big[\beta_m(T,\varepsilon)\big]\right)^m\overset{m\rightarrow \infty}{\longrightarrow} 0.
\end{split}
\end{equation*}
Hence, for every $T>0$ there exists a finite, $\N$-valued random variable $N_T$ such that 
\begin{equation*}\label{CB:def:NT}
\tilde{\Psi}^{Q}_{N_T}\big(Y, \tilde{R}^+, \tilde{R}^-\big) = \tilde{\Psi}^{Q}\big(Y, \tilde{R}^+, \tilde{R}^-\big) \quad \Pro\text{-a.s. on } [0,T].
\end{equation*}
\vskip6pt

\noindent
\underline{Step 3:} Finally, the continuous mapping theorem allows us to conclude that
\[\tilde{\Psi}(\tilde{S}^{(n)}_0, X^{(n)}, \tilde{R}^{+,(n)} ,\tilde{R}^{-,(n)}) \rightarrow \tilde{\Psi}(\tilde{S}_0, X,\tilde{R}^+, \tilde{R}^-) \quad \Pro\text{-a.s.}\]
As there are only finitely many price changes on each compact interval (cf.~Step 2), we have for any $\varepsilon>0$ and $t\geq0$,
\[\Pro\left[\vn |c_n(t)|>\varepsilon\right]=\Pro\left[|c_n(t)|>\frac{\varepsilon}{\vn}\right]
\leq \Pro\left[\tilde{\tau}^*_{\lfloor \varepsilon/\vn\rfloor}\leq t\right]\rightarrow 0\quad\text{as}\quad n\rightarrow\infty.\]
Thus, Theorem \ref{CB:res:repactive} yields that also $\tilde{S}^{(n)} \rightarrow \tilde{\Psi}(\tilde{S}_0, X,\tilde{R}^+, \tilde{R}^-)$  $\Pro$-a.s. 
\end{proof}

\subsection{Identification of the limit}

Below we identify the process $g(W)$ for a planar Brownian motion $W$ as a special two-dimensional semimartingale reflecting Brownian motion (SRBM) absorbed at the origin, which we will call a sum-conserving SRBM absorbed at the origin.

\begin{Def}\label{CB:def:SRBM}
Let $W$ be a two-dimensional linear Brownian motion on a filtered probability space $(\Omega,\F,\mathbb{F},\Pro)$ with $W(0)\neq(0,0)$ $\Pro$-a.s. The pair $(Z,l)$ of continuous, adapted, $\R^2$-valued processes is a sum-conserving SRBM with absorption associated with $W$ if 
\begin{equation*}
        Z(t) = \begin{cases}
        W(t) + l(t)\mathcal{R} &\quad \text{fo } t \leq \tau:=\inf\{t\geq 0:Z(t) = (0,0)\}\\
        0 &\quad \text{for } t > \tau
        \end{cases},
    \end{equation*}
    where
    \begin{equation*}\label{CB:def:reflMatrix}
        \mathcal{R} = \begin{pmatrix*}[r] 1&-1\\-1&1\end{pmatrix*}
    \end{equation*} 
    and for $i \in \{1,2\}$ the $i$-th component $l_i$ of $l$ satisfies
        \begin{itemize}
            \item[a)] $l_i$ is non-decreasing with $l_i(0) = 0$,
            \item[b)] $\int_0^t Z_i(t) dl_i(t) = 0$ for all $t \geq 0,$ and
            \item[c)] $l_i(t) = l_i(\tau)$ for all $t\geq \tau.$
        \end{itemize}
\end{Def}

Multidimensional reflecting Brownian motion has been extensively studied in the literature for various choices of reflection matrices, starting with the pioneering works \cite{HR81,VW85,W85,T90,TW93}. We note, however, that the reflection matrix $\mathcal{R}$ considered above is singular and does in particular not satisfy the completely-$\mathscr{S}$ property often encountered in the literature on reflecting Brownian motion as it allows a pathwise construction of SRBM by means of the Skorokhod map on $[0,\infty)$ and thus ensures the semimartingale property for all times. For $\mathcal{R}$ as defined above, this does not work anymore and it is indeed important to set $Z(t) = (0,0)$ for all $t > \tau$ for $Z$ to remain a semimartingale after time $\tau$. In fact, Definition \ref{CB:def:SRBM} is a special case of the dual skew symmetric SRBM introduced in \cite{EFH21,FR22}, which requires the reflection vectors to point in opposite directions. The existence and uniqueness of this process is shown in \cite[§4.2 and §4.3]{TW93}, cf.~also \cite[Proposition 2]{FR22}. 
In particular, it is known that the process $(Z,l)$ defined above is indeed a continuous semimartingale (cf.~\cite[Corollary 2]{W85}) and a strong Markov process (cf.~\cite[Theorem 4.3]{R06}). In the following, with a slight abuse of terminology, we will also refer to $Z$ itself as a sum-conserving SRBM with absorption associated with $W$.

\par

The next result shows that for a planar Brownian motion $W$, the process $g(W)$ is a sum-conserving SRBM with absorption associated with $W$. In the following, we denote for any $\R^2$-valued continuous semimartingale $X$ the component-wise local time of $X$ at zero by
\begin{equation*}\label{CB:def:localTime2}
    L^{(2)}_t(X) := \left(L_t(\pi^{(2)}_1 X), L_t(\pi^{(2)}_2 X)\right), \quad t \geq 0,
\end{equation*}
i.e.~$L(\pi^{(2)}_iX)$ denotes the local time of $\pi^{(2)}_iX$ at zero for $i=b,a$. 

\begin{Prop}\label{CB:res:idgW}
    Let $W=(W_1,W_2)$ be a two-dimensional linear Brownian motion with $W(0)\neq(0,0)$ $\Pro$-a.s.  
    Then $g(W)$ is a sum-conserving SRBM absorbed at the origin associated with $W$. Moreover, $\htau_\infty(W)=\inf\{t\geq0:\ g(W)=(0,0)\}$ $\Pro$-a.s.~and 
    \begin{equation}\label{CB:def:newRepgW}
    g(W)(t) = W(t) + \frac{1}{2}L^{(2)}_t(g(W))\mathcal{R} \quad\text{on}\quad[0,\htau_\infty(W)].
    \end{equation}
\end{Prop}

The proof of Proposition \ref{CB:res:idgW} can be found in Appendix \ref{CB:app:SRBM}. 
With the help of Proposition \ref{CB:res:idgW} we can now characterize the limit $\tilde{S}:=\tilde{\Psi}(\tilde{S}_0,X,\tilde{R}^+,\tilde{R}^-)$ from Theorem \ref{CB:res:convQM-1}.

\begin{The}
\label{CB:res:idQ-1}
    Let the assumptions of Theorem \ref{CB:res:convQM-1} be satisfied and define a sequence $(\tilde{W}_k)_{k\in \N_0}$ of four-dimensional linear Brownian motions, each with mean $\mu$ and covariance matrix $\Sigma$, given by
    \begin{equation*}\label{CB:def:tildeB}
    \tilde{W}_k := \tilde{Q}(\tilde{\tau}^*_k) + X(\cdot + \tilde{\tau}^*_k) - X(\tilde{\tau}^*_k),\quad k\in\N_0.
    \end{equation*}
        Then $\lim_k\tilde{\tau}_k^*=\infty$ $\Pro$-a.s.~and the weak limit $\tilde{S} = (\tilde{B},\tilde{Q},\tilde{C})$ of $\tilde{S}^{(n)}$ satisfies:
    \begin{enumerate}
        \item[i)] For all $k\geq 0$ and $t\in [0, \tilde{\tau}^*_{k-1}-\tilde{\tau}^*_k),$
        \[\left(\tilde{Q}^{i,F}, \tilde{Q}^{i,G}\right)(t + \tilde{\tau}^*_k) = \pi_i \tilde{W}_k(t) + \frac{1}{2}L^{(2)}_{t}\left(g(\pi_i \tilde{W}_k)\right)\mathcal{R}\quad\text{for}\quad i=b,a,\]
        i.e.~on $[\tilde{\tau}^*_k, \tilde{\tau}^*_{k+1})$ the process $(\tilde{Q}^{i,F}, \tilde{Q}^{i,G})$ is a sum-conserving SRBM for $i=b,a$. 
        \item[ii)] For all $k\geq 0$, 
        $\tilde{Q}$ is reinitialized at time $\tilde{\tau}_k^*$ at value 
        $\tilde{Q}(\tilde{\tau}^*_k) = \tilde{R}^{+}_k \1_{\{\tilde{\tau}_{a,k} = \tilde{\tau}^*_k\}} +\tilde{R}^-_k\1_{\{ \tilde{\tau}_{b,k} = \tilde{\tau}^*_k\}}$.
	\item[iii)] For all $k\geq0$ and $t \in [0, \tilde{\tau}^*_{k+1} - \tilde{\tau}^*_k],$
	\begin{equation}\label{CB:eq:Ci}
    \begin{split}
        &\tilde{C}(t + \tilde{\tau}^*_k) - \tilde{C}(\tilde{\tau}^*_k)=\\
        & \qquad \frac{1}{2} \left(L_{t}(\pi^{(2)}_2 g(\pi_b \tilde{W}_k)) + L_{t}(\pi^{(2)}_1 g(\pi_a \tilde{W}_k))\right) - \frac{1}{2}\left( L_{t}(\pi^{(2)}_1 g(\pi_b \tilde{W}_k)) + L_{t}(\pi^{(2)}_2 g(\pi_a \tilde{W}_k))\right).
\end{split}
\end{equation}
	In particular, $\tilde{C}$ is a continuous process of finite variation.
\item[iv)] $\tilde{B}=(\tilde{B}^F,\tilde{B}^F)$ is a piecewise constant, càdlàg process which
\begin{itemize}
	\item increases by $\delta$ at time $\tilde{\tau}_k^*$ if $\tilde{\tau}_{a,k} = \tilde{\tau}^*_k$, 
	i.e.~if $h(\tilde{W}_k)$ hits first the $x$-axis,
	\item decreases by $\delta$ at time $\tilde{\tau}_k^*$ if $\tilde{\tau}_{b,k} = \tilde{\tau}^*_k$, 
		i.e.~if $h(\tilde{W}_k)$ hits first the $y$-axis,
\end{itemize}
i.e.
\[\tilde{B}^F(t)  = \tilde{B}^F_0 + \delta \sum_{0\leq s \leq t}\left(\1_{\left\{\pi^{(2)}_2 \left(h \circ \tilde{Q}\right)(s-) = 0\right\}} - \1_{\left\{ \pi^{(2)}_1 \left(h\circ \tilde{Q}\right)(s-) = 0\right\}}\right),\quad t\geq 0.\]
    \end{enumerate}
\end{The}

\begin{proof}
According to Theorem \ref{CB:res:convQM-1} we have $\lim_k\tilde{\tau}_k^*=\infty$ $\Pro$-a.s. By definition we have for all $k\geq 0$,
$$\tilde{Q}\left(\cdot + \tilde{\tau}^*_k\right)  = \tilde{\Psi}^{Q}\left(\tilde{Q}_0 + X, \tilde{R}^+, \tilde{R}^-\right)\left(\cdot + \tilde{\tau}^*_k\right) = G(\tilde{W}_k) \quad \text{on } \left[0, \tilde{\tau}^*_{k+1} - \tilde{\tau}^*_k\right).$$
Hence, part i) follows from the definition of $G$ and Proposition \ref{CB:res:idgW}. Part ii) follows directly from Theorem \ref{CB:res:convQM-1}. By the definition of $\tilde{\Psi}^{C}$, we have $\tilde{C}(\cdot + \tilde{\tau}^*_k) - \tilde{C}(\tilde{\tau}^*_k) = \bar{G}(\tilde{W}_k)$ on $ [0, \tilde{\tau}^*_{k+1} - \tilde{\tau}^*_k]$ for $k \geq 0$. Hence, part iii) follows from the definition of $\bar{G}$, Lemma \ref{CB:res:sumCompg}, and Proposition \ref{CB:res:idgW}. This implies that $\tilde{C}(\cdot + \tilde{\tau}^*_{k}) - \tilde{C}(\tilde{\tau}^*_k)$ is a process of finite variation on each interval $[0, \tilde{\tau}^*_{k+1} - \tilde{\tau}^*_k]$. Since there are only finitely many price changes on each compact interval, $\tilde{C}$ is the sum of finitely many continuous processes of finite variation on each compact interval and thus itself of finite variation. Part iv) follows directly from the definition of $\tilde{\Psi}^B.$
\end{proof}

\begin{Rmk}
\label{CB:rem:cross-border}
    By Theorem \ref{CB:res:idQ-1} the queue size process $\tilde{Q}^{(n)}$ is approximated by a semimartingale with a non-trivial martingale part, whereas the capacity process $\tilde{C}^{(n)}$ is approximated by a continuous process of finite variation. Therefore, the domestic trading activity is of much greater magnitude than the cross-border trading activity. 
\end{Rmk}

\subsection{Analysis of the limiting dynamics}\label{CB:sec:ProbLimit}

In the microscopic order book model, the price and queue size processes result from a complex interplay of the order flow processes in both countries and therefore quantitative properties cannot be easily derived. However, Theorem \ref{CB:res:idQ-1} shows that the limiting dynamics are analytically quite tractable and therefore invite explicit computations. 

\subsubsection{Price changes}\label{CB:sec:simStudy-priceChanges}

In this subsection, 
we compute the distribution of the duration until the next price change analytically and study its dependence on the parameters. 
According to Theorem \ref{CB:res:idQ-1} the limiting price process $\tilde{B}^F$ can be described through the hitting times of a planar Brownian motion with reinitializations inside the positive orthant. As this process is a Markov process, 
we only study the time of the first price change $\tilde{\tau}_1^*$ in the following. 

\begin{Prop}
\label{CB:cor:pricechanges}
For any $t\geq0$, the survival probability of $\tilde{\tau}_1^*$ is given by 
\begin{equation}\label{CB:eq:survival}
\Pro[\tilde{\tau}^*_1>t]=k\left((h\circ\tilde{Q})(0),\mu_h,\Sigma_h;t\right),
\end{equation}
where the function $k$ is defined in \eqref{CB:def:condDist-Duration} and
\begin{align*}
\mu_h:= 
\begin{pmatrix}
    \mu_b\\
    \mu_a
\end{pmatrix}
:=
\begin{pmatrix}
    \mu^{b,F}+\mu^{b,G}\\
    \mu^{a,F}+\mu^{a,G}  
\end{pmatrix}, \quad
\Sigma_h:=
\begin{pmatrix}
    \sigma^2_{b} & \rho_h \sigma_{b} \sigma_{a}\\
    \rho_h\sigma_{b}\sigma_{a} & \sigma^2_{a}\\
\end{pmatrix}
\end{align*}
with
\begin{align*}
    \sigma^2_{b} &:= (\sigma^{b,F})^2 + 2\sigma^{(b,F),(b,G)} + (\sigma^{b,G})^2,\\
    \sigma^2_{a} &:= (\sigma^{a,F})^2+2\sigma^{(a,F),(a,G)}+(\sigma^{a,G})^2,\\
    \rho_h \sigma_{b} \sigma_{a} &:= \sigma^{(b,F),(a,F)}+\sigma^{(b,F),(a,G)}+\sigma^{(b,G),(a,F)}+\sigma^{(b,G),(a,G)}.
\end{align*}
Moreover, if $a_t=a_t(\mu_h,\Sigma_h)\leq 0$, where the constant $a_t$ is defined in Theorem \ref{CB:thm:2dBM},
then $\tilde{\tau}_1^*$ is finite $\Pro$-almost surely and the probability of the first price movement being upwards is gven by
$$\Pro\left[\pi_1^{(2)}(h \circ \tilde{Q})(\tilde{\tau}_1^*-)>0\right]=\int_0^\infty\int_0^\infty \exp\left(-z\left(\frac{\mu_{b}}{\sigma_{b}^2}-\frac{\rho_h d_2}{\sigma_{b}}\right)\right)p(z,t)dtdz,$$
where $p(z,t)$ and $d_2$ are defined in Theorem \ref{CB:thm:2dBM}. Especially, if $\mu_h = (0,0)$, we have 
\begin{equation*}
\Pro[\tilde{\tau}_1^* > t]=\frac{2}{t}\exp\left(- \frac{U}{2t}\right)\sum_{j = 1}^{\infty} \sin\left(\frac{j\pi \theta_0}{\alpha}\right)\frac{\big(1-(-1)^j\big)}{j\pi}\int_0^{\infty} r \exp\left(-\frac{r^2}{2t}\right) I_{j\pi/\alpha}\left(\frac{r \sqrt{U}}{t}\right) dr
\end{equation*}
and
$$\Pro\left[\pi_1^{(2)}(h \circ \tilde{Q})(\tilde{\tau}_1^*-)>0\right]=\frac{\alpha-\theta_0}{\alpha},$$
where $\alpha$, $\theta_0$, $U$ are defined in Theorem \ref{CB:thm:2dBM} and $I$ is the modified Bessel function of the first kind. 
\end{Prop}

\begin{proof}
Theorems \ref{CB:res:convQM-1} and \ref{CB:res:idQ-1} imply that
$\tilde{\tau}^*_1=\inf\left\{t\geq0: (\pi_1Y(t)+\pi_3Y(t))(\pi_2Y(t)+\pi_4Y(t))=0\right\}$,
where $Y$ is a linear four-dimensional Brownian motion started inside the positive orthant. Hence, the time of the first price change equals the first exit time of a planar Brownian motion from the positive orthant and can be computed analogously to \cite{CL21}, using Theorem \ref{CB:thm:2dBM}.

Moreover, the first price change is upwards if and only if the cumulative bid queue is strictly positive at that time. Hence, the probability of observing an upwards price change equals
\[\Pro\left[(\pi_1^{(2)}\circ h)(\tilde{Q})(\tilde{\tau}_1^*-)>0\right]=\Pro\left[\pi_1Y(\tilde{\tau}_1^*)+\pi_3Y(\tilde{\tau}_1^*)>0\right]\]
and can be explicitly computed via formula \eqref{CB:def:jointdensitytau} in the Appendix. The simplified expressions for $\mu_h = (0,0)$ also follow directly from Theorem \ref{CB:thm:2dBM}.
\end{proof}

\begin{figure}[H]
	\centering
	 \includegraphics[scale = 0.32]{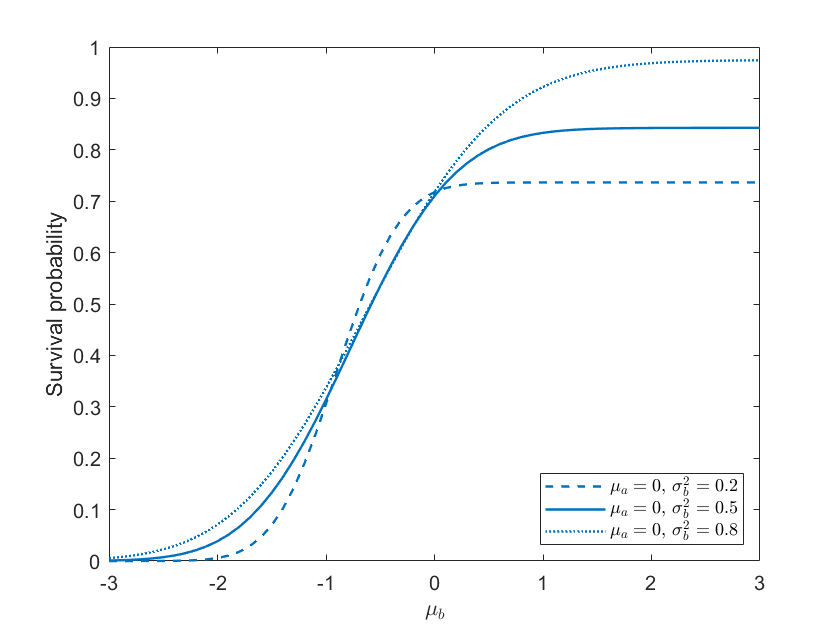} \quad
	 \includegraphics[scale = 0.32]{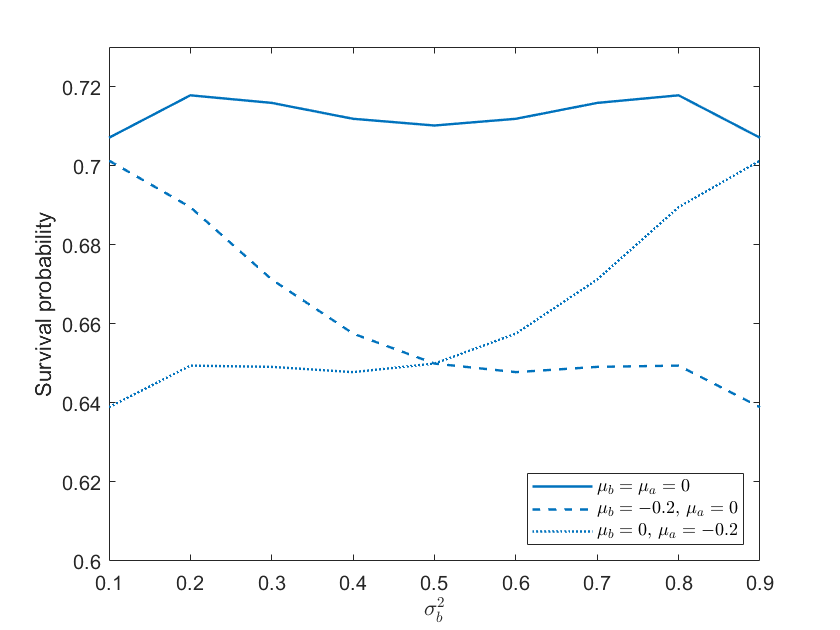}
    \caption{\small Survival probability \eqref{CB:eq:survival} for an independent order flow and different values of $\mu_b,\mu_a$ and $\sigma_b,\sigma_a$.} 
    \label{CB:fig:survProp1}
\end{figure}

Next we study numerically the dependence of formula \eqref{CB:eq:survival} on the model parameters for $t = 1$ and $\tilde{Q}^b(0) = \tilde{Q}^a(0) = 1$. First, we assume that the order size vectors are independent over time. In this case, we have $\sigma_{b}^2 = (\sigma^{b,F})^2 + (\sigma^{b,G})^2$,  $\sigma_{a}^2=(\sigma^{a,F})^2 + (\sigma^{a,G})^2,$ $\rho_h = 0$, and
\[\sum_{(i,I) \in \{b,a\} \times \{F,G\}}(\sigma^{i,I,(n)})^2 + o\left(\vn\right)=\sum_{(i,I) \in \{b,a\} \times \{F,G\}}\Pro\left[(\phi^{(n)}_1, \psi^{(n)}_1) =  (i,I)\right] = 1,\quad n\in\N.\]
Hence, $\sigma_{b}^2 = 1 - \sigma_{a}^2$. As expected the survival probability increases with the mean, while the variance affects the skewness and kurtosis of the curve, cf. Figure \ref{CB:fig:survProp1} (left). 
For centered mean $\mu_{b} = \mu_{a} = 0$, the variance has only little influence on the survival probability, while the symmetry of the curve follows from the relation $\sigma^2_{b} = 1-\sigma^2_{a}$, cf. Figure \ref{CB:fig:survProp1} (right). In contrast, if $\mu_{b}<0$ (resp.~$\mu_{a}<0$) the probability is decreasing (resp.~increases) in the bid variance $\sigma^2_{b},$ since a higher bid variance results in a higher likelihood of the cumulative bid queue to hit zero and simultaneously in a lower likelihood of the cumulative ask queue to hit zero.

\begin{figure}[H]
	\centering
	 \includegraphics[scale = 0.32]{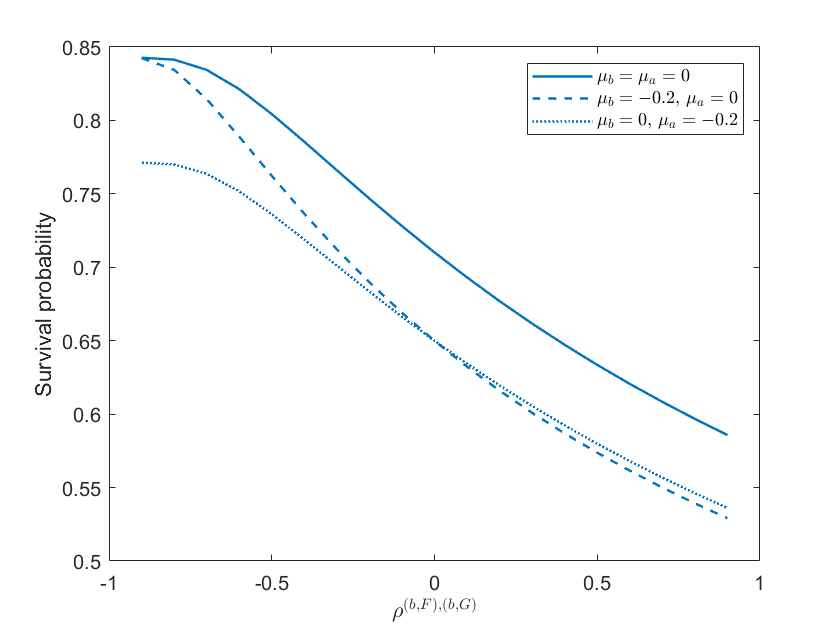} \quad
	 \includegraphics[scale = 0.32]{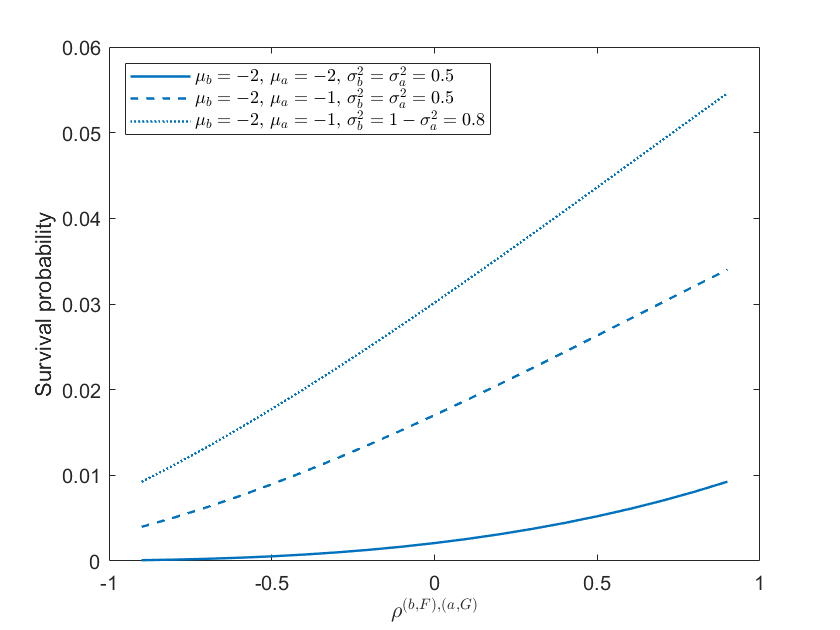}
    \caption{\small Influence of correlations in the order flow on the survival probability in \eqref{CB:eq:survival}.\\ Left: Influence of $\rho^{(b,F),(b,G)}$ for different values of $\mu_b,\mu_a$ and  $(\sigma^{i,I})^2 = 0.25$.\\ Right: Influence of $\rho^{(b,F),(a,G)}$ for different values of $\mu_a$ and $\sigma_a^2=1-\sigma_b^2$.}
    \label{CB:fig:survProp2}
\end{figure}

Second, we study the influence of correlations between order events of different origins ($F$ and $G$) on the survival probability. 
In Figure \ref{CB:fig:survProp2} (left) we set all correlation parameters except $\rho^{(b,F), (b,G)}$ to zero. In this case a negative (resp.~positive) correlation decreases (resp.~increases) the bid variance and hence decreases (resp.~increases) the probability that the cumulative bid queue will hit zero, while the ask side is unaffected. For this reason the survival probability decreases with $\rho^{(b,F),(b,G)}.$

In Figure \ref{CB:fig:survProp2} (right), we study the influence of the correlation parameter $\rho^{(b,F),(a,G)}$, while the other correlation parameters are set to zero. Since $\rho^{(b,F),(a,G)}$ only affects the correlation parameter $\rho_h$ in the equation for the survival probability, the relation $\sigma^2_{b} = 1-\sigma^2_{a}$ is again satisfied. We observe that the survival probability increases with $\rho^{(b,F),(a,G)}$ for different choices of $\mu_a$ and $\sigma_a^2$ . 

\subsubsection{Cross-border trading activity between price changes}

A key insight from our model is that the cross-border trading activity is of smaller magnitude than the domestic trading activity, cf.~Remark \ref{CB:rem:cross-border}. In this subsection, we want to study how the model parameters affect the cross-border trading activity, which mathematically amounts to an analysis of the reflecting behaviour. For this we recall from Theorem \ref{CB:res:idQ-1} that the queue size process $(\tilde{Q}^{i,F},\tilde{Q}^{i,G})$ on the bid ($i=b$) resp.~ask ($i=a$) side behaves before time $\hat{\tau}^{i,G}:=\inf\{t\geq0:\ \tilde{Q}^{i,G}(t)=0\}$ like
\begin{align*}
    \left(\tilde{Q}^{i,F},\tilde{Q}^{i,G}\right)=\left(Y^{i,F}+\sup_{s\leq t}\left(-Y^{i,F}(s)\right)^+,Y^{i,G}-\sup_{s\leq t}\left(-Y^{i,F}(s)\right)^+\right),
\end{align*}
where the process  $(Y^{i,F},Y^{i,G})$ is a two-dimensional linear Brownian motion. Our goal is to study the duration of the direction of cross-border trading, i.e.~we would like to understand when the initially importing country (w.l.o.g.~$F$) becomes the exporting country. The next proposition characterizes both, the survival probability of $\hat{\tau}^{i,G}$ as well as the distribution function of the size of the opposite queue at time $\htau^{i,G}$, i.e.~$\tilde{Q}^{i,F}(\hat{\tau}^{i,G})$, as the unique solution of an interface problem for the inhomogeneous respectively homogeneous two-dimensional heat equation.

\begin{Prop}\label{CB:thm:interface}
For any $t>0$,
\begin{equation}\label{CB:res:hitting}
\Pro[\htau^{i,G}>t]=F\left(\tilde{Q}^{i,G}(0),\, 2\tilde{Q}^{i,F}(0)+\tilde{Q}^{i,G}(0);\, t\right),
\end{equation}
where $F$ is the unique bounded $C^{1,1}(\R_+^2,\R_+)\cap C^{2,1}(\R_+^2\setminus\{(1,1)\R_+\},\R_+)$ solution of the interface problem
\begin{equation}\label{kolmogorovPDE}
F_t=\begin{cases}\frac{\sigma_1^2}{2}F_{z_1z_1}+\rho\sigma_1\sigma_2F_{z_1z_2}+\frac{\sigma_2^2}{2}F_{z_2z_2}+\mu_1F_{z_1}+\mu_2F_{z_2}&:z_2>z_1\\
\frac{\sigma_2^2}{2}F_{z_1z_1}+\rho\sigma_1\sigma_2F_{z_1z_2}+\frac{\sigma_1^2}{2}F_{z_2z_2}+\mu_2F_{z_1}+\mu_1F_{z_2}&:z_2<z_1
\end{cases}
\end{equation}
with boundary conditions 
$$F(z_1,z_2;0)=1,\qquad F(z_1,0;t)=F(0,z_2;t)=0\ \text{ for all }z_1,z_2>0,$$ and parameters
\begin{align*}
\mu_1&:=\mu^{i,G}, & \mu_2&:=\mu^{i,G}+2\mu^{i,F}, & &\\
\sigma_1^2&:=(\sigma^{i,G})^2,  &
    \sigma_2^2&:=(\sigma^{i,G})^2+4\sigma^{(i,F),(i,G)}+4(\sigma^{i,F})^2,  & 
    \rho\sigma_1\sigma_2&:=(\sigma^{i,G})^2+2\sigma^{(i,F),(i,G)}.
\end{align*}
Moreover, if $\mu^{i,G}\leq0$ or $\mu^{i,G}+\mu^{i,F}\leq0$, then $\htau^{i,G}<\infty$ $\Pro$-almost surely and for any $u>0$ we have
\[\Pro\left[\tilde{Q}^{i,F}\left(\htau^{i,G}\right)>u\right]=G^u\left(\tilde{Q}^{i,G}(0),\,2\tilde{Q}^{i,F}(0)+\tilde{Q}^{i,G}(0)\right),\]
where $G^u$ is the unique bounded $C^{1}(\R_+^2)\cap C^{2}(\R_+^2\setminus\{(1,1)\R_+\})$ solution of the interface problem
\begin{equation*}\label{PDEGu}
0=\begin{cases}\frac{\sigma_1^2}{2}G^u_{z_1z_1}+\rho\sigma_1\sigma_2G^u_{z_1z_2}+\frac{\sigma_2^2}{2}G^u_{z_2z_2}+\mu_1G^u_{z_1}+\mu_2G^u_{z_2}&:z_2>z_1>0\\
\frac{\sigma_2^2}{2}G^u_{z_1z_1}+\rho\sigma_1\sigma_2G^u_{z_1z_2}+\frac{\sigma_1^2}{2}G^u_{z_2z_2}+\mu_2G^u_{z_1}+\mu_1G^u_{z_2}&:0<z_2<z_1
\end{cases}
\end{equation*}
with boundary condition $G^u(0,z)=G^u(z,0)=\1_{(u,\infty)}(z/2),\ z>0$.
\end{Prop}

\begin{proof}
Since
\begin{align*}
   \htau^{i,G}&=\inf\left\{t\geq0:\tilde{Q}^{i,G}(t)=0\right\}=\inf\left\{t\geq0:Y^{i,G}(t)=\sup_{s\leq t}\left(-Y^{i,F}(s)\right)^+\right\},
\end{align*}
we have for any $t\geq0$,
\begin{equation}\label{CB:eq:hit1}
\begin{split}
    \Pro\left[\htau^{i,G}>t\right]&=\Pro\left[\forall s\in[0,t]: Y^{i,G}(s)>\sup_{s\leq t}\left(-Y^{i,F}(s)\right)^+\right]\\
    &=\Pro\left[\forall s\in[0,t]: Y^{i,G}(s)+Y^{i,F}(s)>Y^{i,F}(s)+\sup_{s\leq t}\left(-Y^{i,F}(s)\right)^+\right].
    \end{split}
\end{equation}
With an abuse of notation, let us denote by $(\tilde{Y}^{i,F},Y^{i,F},Y^{i,G})$ the unique weak solution of the Tanaka SDE
\[\tilde{Y}^{i,F}(t)=\tilde{Q}^{i,F}(0)+\int_0^t\sign\big(\tilde{Y}^{i,F}(s)\big)dY^{i,F}(s),\quad t\geq0,\]
where $(Y^{i,F},Y^{i,G})$ has the same law as before. From Tanaka's formula we have
$$\big|\tilde{Y}^{i,F}(t)\big|={Y}^{i,F}(t)+L_t(\tilde{Y}^{i,F})=Y^{i,F}(t)+\sup_{s\leq t}\left(-Y^{i,F}(s)\right)^+,\quad t\geq0,$$
where the second equality follows from \cite[Lemma 2.12]{Z17}. Therefore,
\begin{equation}\label{CB:eq:hit2}
\begin{split}
    \Pro\left[\htau^{i,G}>t\right]&=\Pro\left[\forall s\in[0,t]: Y^{i,G}(s)+Y^{i,F}(s)>\big|\tilde{Y}^{i,F}(s)\big|\right]\\
    &=\Pro\left[\forall s\in[0,t]: Z^{i,1}(s)Z^{i,2}(s)>0\right]
    \end{split}
\end{equation}
with
\begin{align*}
    Z^1(t):=Y^{i,F}(t)+Y^{i,G}(t)-\tilde{Y}^{i,F}(t),\qquad Z^2(t):=Y^{i,F}(t)+Y^{i,G}(t)+\tilde{Y}^{i,F}(t)
\end{align*}
satisfying the SDE
\begin{align*}
    dZ^1(t)&=dY^{i,G}(t)+\big(1-\sign(Z^2(t)-Z^1(t))\big)dY^{i,F}(t),\\
    dZ^2(t)&=dY^{i,G}(t)+\big(1+\sign(Z^2(t)-Z^1(t))\big)dY^{i,F}(t).
\end{align*}
Note that $Z=(Z^1,Z^2)$ satisfies the strong Markov property by the weak uniqueness property of the Tanaka SDE, cf.~\cite[Theorems 32.7 and 32.11]{K21}. Hence, we are looking for the probability that a two-dimensional, homogeneous diffusion with starting value $z:=(y,2x+y)\in\R^2_{+}$ does not hit the axes until time $t$.  By classical diffusion theory, $F(z_1,z_2;t):=\Pro[\htau^{i,G}>t|Z^1(0)=z_1,Z^2(0)=z_2]$ solves the Kolmogorov backward equation
\begin{equation*}
F_t=\begin{cases}\frac{\sigma_1^2}{2}F_{z_1z_1}+\rho\sigma_1\sigma_2F_{z_1z_2}+\frac{\sigma_2^2}{2}F_{z_2z_2}+\mu_1F_{z_1}+\mu_2F_{z_2}&:z_2>z_1\\
\frac{\sigma_2^2}{2}F_{z_1z_1}+\rho\sigma_1\sigma_2F_{z_1z_2}+\frac{\sigma_1^2}{2}F_{z_2z_2}+\mu_2F_{z_1}+\mu_1F_{z_2}&:z_2<z_1
\end{cases}
\end{equation*}
with boundary conditions $F(z_1,z_2;0)=1,\ F(z_1,0;t)=F(0,z_2;t)=0$ for all $z_1,z_2>0$.  
Finally, we note that the uniqueness of a bounded $C^{1,1}(\R_+^2,\R_+)\cap C^{2,1}(\R_+^2\setminus\{(1,1)\R_+\},\R_+)$ solution to the above interface problem can be established as usual by the energy method, imposing smooth fit conditions at the interface.\\
Moreover, if follows from \eqref{CB:eq:hit1} and \eqref{CB:eq:hit2} that 
\begin{align*}
    \Pro\left[\htau^{i,G}=\infty\right]\leq \Pro\left[\forall t\geq0: Y^{i,G}(t)>0\right]\wedge\Pro\left[\forall t\geq0: Y^{i,F}(t)+Y^{i,G}(t)>0\right].
\end{align*}
Hence, if $\mu^{i,G}\leq 0$ or $\mu^{i,G}+\mu^{i,F}\leq 0$, then $\htau^{i,G}<\infty$ $\Pro$-a.s. In this case, the distribution function of $\tilde{Q}^{i,F}(\htau^{i,G})$ can be characterized analogously to the survival probability of $\htau^{i,G}$.
\end{proof}

Below we plot the numerical solution of the interface problem \eqref{kolmogorovPDE} for a symmetric, uncorrelated order flow process.
\begin{figure}[H]
	\centering
	 \includegraphics[scale = 0.55]{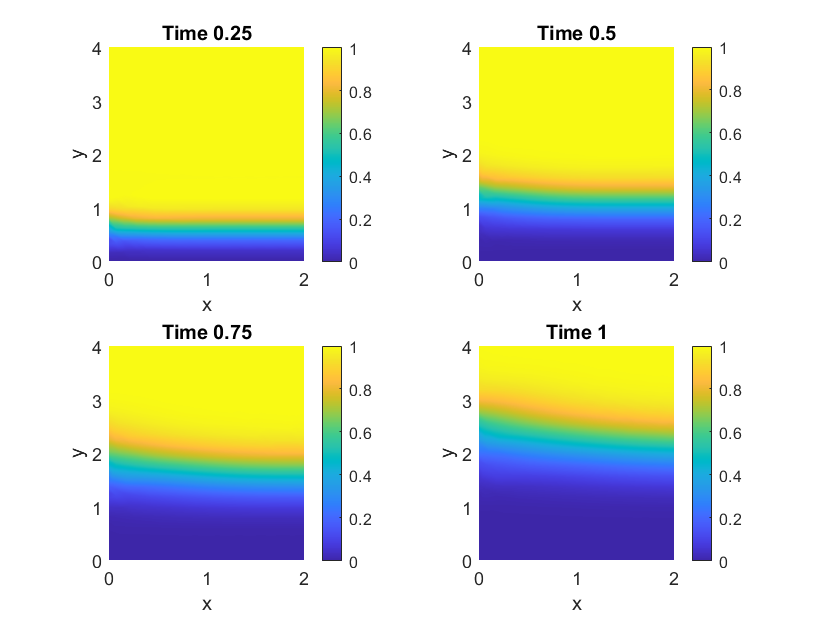}
    \caption{\small $\Pro[\hat{\tau}^{i,G} > t]$ as a function of $\tilde{Q}^{i,F}(0)$ ($x$-axis) and $\tilde{Q}^{i,G}(0)$ ($y$-axis) for $t = 0.25,{} 0.5,{}0.75,{}1$; $(\sigma^{i,F})^2 = (\sigma^{i,G})^2 = 0.25$, $\rho^{(i,F),(i,G)} = 0$, $\mu^{i,F} = \mu^{i,G} = -2.$ }
    \label{CB:fig:solPDEs}
\end{figure}

\begin{Rmk}
Note that \eqref{CB:res:hitting} also holds true for the starting values $\tilde{Q}^{i,F}=0$, $\tilde{Q}^{i,G}>0$. Therefore, exchanging the role of $F$ and $G$, one can obtain the distribution of the reflection times at alternating axes, given the value of the queue size process at the previous reflection time. 
\end{Rmk}

\section{The inactive dynamics}\label{CB:sec:anaID}

In this section we briefly describe the heavy traffic approximation for the inactive dynamics $\dbtilde{S}^{(n)}$ specified in Definition \ref{CB:def:inactive}. Similarly to the analysis of the active dynamics in Section \ref{CB:sec:anaAD}, one can define a function 
$\dbtilde{\Psi}: E \times D(\R_+, \R^4) \times (\R^4_+)^{\N} \times (\R^4_+)^{\N} \rightarrow D(\R_+, E)$
such that for all $n\in\N$,
\begin{align*}
    \dbtilde{S}^{(n)} =  \dbtilde{\Psi}\left(\dbtilde{S}^{(n)}_0, X^{(n)}, \dbtilde{R}^{+,(n)}, \dbtilde{R}^{-,(n)}\right),
\end{align*}
 where $\dbtilde{S}^{(n)}_0$ is the initial state, $X^{(n)}$ is the net order flow process, and $\dbtilde{R}^{+,(n)}$ resp.~$\dbtilde{R}^{-,(n)}$ are the reinitialization values after price increases resp.~decreases. For details on the construction of the function $\dbtilde{\Psi}$ we refer the interested reader to \cite[Definition 1]{CL21} and \cite[Chapter 2.4]{PhDCM}. The following theorem is a straight-forward extension of \cite[Theorem 2, Proposition 1]{CL21} to two correlated, non-interacting LOBS and follows by similar arguments as Theorem \ref{CB:res:convQM-1}.\footnote{
In \cite[Proposition 1]{CL21} the weak convergence of prices is established in the weaker Skorokhod $M_1$-topology. However, by similar arguments as in the proof of Theorem \ref{CB:res:convQM-1}, one even obtains weak convergence of prices in the stronger Skorokhod $J_1$-topology.}

\begin{The}[Limit theorem for the inactive dynamics $\dbtilde{S}^{(n)}$]\label{CB:res:convS2} Let Assumptions \ref{CB:ass:prob}, \ref{CB:ass:scaling}, and \ref{CB:ass:R} be satisfied and suppose that for every $n\in\N$, $\dbtilde{S}^{(n)}_0=(\dbtilde{B}^{(n)}_0,\dbtilde{Q}^{(n)}_0,\dbtilde{C}^{(n)}_0)$ is a $(\delta \Z)^2 \times \vn \N^4 \times \vn\Z$-valued random variable on $(\Omega^{(n)},\F,\Pro^{(n)})$ such that $(\dbtilde{S}^{(n)}_0,X^{(n)})\Rightarrow (\dbtilde{S}_0,X)$, where $X$ is as in Proposition \ref{CB:res:fCLT} and  $\dbtilde{S}_0=(\dbtilde{B}_0,\dbtilde{Q}_0,\dbtilde{C}_0)$ with $\dbtilde{Q}_0\in(0,\infty)^4$. Then
    \begin{align*}
        \dbtilde{S}^{(n)} \Rightarrow \dbtilde{S}=(\dbtilde{B}, \dbtilde{Q}, \dbtilde{C}) :=\dbtilde{\Psi}\left(\dbtilde{S}_0, X, \dbtilde{R}^+, \dbtilde{R}^-\right)
    \end{align*}
    in the Skorokhod topology on the space $D(\R_+, E)$, 
    \begin{align*}
    \dbtilde{R}^+_k := \Phi\left(\dbtilde{Q}(\dbtilde{\tau}^*_k-), \epsilon^+_k\right)  \quad \text{ and } \quad \dbtilde{R}^-_k := \Phi\left(\dbtilde{Q}(\dbtilde{\tau}^*_k-), \epsilon^-_k\right),\quad k\in\N,
    \end{align*}
    for independent sequences $(\epsilon^+_k)_{k\geq 1}$, $(\epsilon^-_k)_{k\geq 1}$ of i.i.d.~random variables with $\epsilon^+_k \sim f^+$, $\epsilon^-_k \sim f^-$ and 
    $$\dbtilde{\tau}_0^*:=0,\qquad\dbtilde{\tau}_k^*:=\inf\left\{t>\dbtilde{\tau}_{k-1}^*: \min_{i=1,\dots,4}\pi_i\dbtilde{Q}(t-)=0\right\},\quad k\in\N .$$ 
    Moreover, $\lim_{k\rightarrow\infty}\dbtilde{\tau}_k^*=\infty$ $\Pro$-a.s.~and the dynamics of $\dbtilde{S}$ can be characterized as follows:
    \begin{itemize}
    \item $\dbtilde{C}(t)=\dbtilde{C}_0$ for all $t\geq 0$.
        \item On each interval $[\dbtilde{\tau}_k^*,\dbtilde{\tau}_{k+1}^*)$, the process $\dbtilde{Q}$ is a four-dimensional linear Brownian motion starting in the interior of $\R_+^4$ with drift $\mu$ and covariance matrix $\Sigma$.
        \item $\dbtilde{B} = (\dbtilde{B}^{F}, \dbtilde{B}^G)$ is a piecewise constant c\`adl\`ag process whose components do almost surely not jump simultaneously; for all $t\geq 0$,
    \begin{align*}
        \dbtilde{B}(t) = \dbtilde{B}_0 & + \delta \sum_{0 \leq s \leq t}\Bigg(\1_{\left\{\pi_2\dbtilde{Q}(s-) = 0 \right\}} - \1_{\left\{\pi_1\dbtilde{Q}(s-) = 0\right\}} , \1_{\left\{\pi_4\dbtilde{Q}(s-) = 0\right\}} - \1_{\left\{\pi_3\dbtilde{Q}(s-) = 0\right\}} \Bigg).
    \end{align*}
    \end{itemize}
\end{The}

\section{Proof of Theorem \ref{CB:res:mainTheorem}}\label{CB:sec:proofmainResult}

In this section we prove our main result. The proof will follow by induction, combining the two convergence theorems for the active dynamics (Theorem \ref{CB:res:convQM-1}) and for the inactive dynamics (Theorem \ref{CB:res:convS2}). Similarly to the discrete-time setting, we denote by $(\tau_k)_{k\geq 1}$ the sequence of stopping times at which we observe a price change in the limit process $S$, by $l(t)$ the random number of price changes in $S$ up to time $t$, and by $R^+_k:= \Phi(Q(\tau_k-), \epsilon^+_k),$ $R^-_k:=\Phi(Q(\tau_k-), \epsilon^-_k),$ $k \geq 1,$ the reinitialization values, i.e.~the queue sizes after price changes. Finally, for all $t\geq 0$, the processes
\begin{align*}
\tilde{S}^t &:= \tilde{\Psi}\left(S(t), X(\cdot + t) - X(t), (R^+_{l(t) + j})_{j\geq 1}, (R^-_{l(t) + j})_{j\geq 1}\right),\\
\dbtilde{S}^t &:= \dbtilde{\Psi}\left(S(t), X(\cdot +t) - X(t), (R^+_{l(t)+j})_{j\geq 1}, (R^-_{l(t)+j})_{j\geq 1}\right),
\end{align*}
denote the active respectively inactive dynamics starting from $S(t)$. With this notation Theorem \ref{CB:res:mainTheorem} can be reformulated as follows:

\begin{The}\label{CB:res:mainTheorem2}
    Let Assumption \ref{CB:ass:initialCon}--\ref{CB:ass:R} be satisfied. Then the microscopic cross-border market models $S^{(n)} = (S^{(n)}(t))_{t\geq 0},\ n\in\N,$ converge weakly in the Skorokhod topology on $D(\R_+, E)$ to a continuous-time regime switching process $S=(B,Q,C)$ such that for all $k\in\N_0$,
    \[S\simeq \tilde{S}^{\rho_k}\,\,\text{on}\,\,\, [\rho_k, \sigma_{k+1})\qquad\text{and}\qquad S \simeq \dbtilde{S}^{\sigma_{k+1}}\,\,\text{on}\,\,\, [\sigma_{k+1}, \rho_{k+1}),\]
    where $\rho_k,\sigma_k$ were defined in \eqref{CB:eq:sigmarho}.
\end{The}

\begin{proof}
    As $S_0^{(n)},X^{(n)},\epsilon^{+,(n)},\epsilon^{-,(n)}$ are independent we may assume by Skorokhod representation that all processes are defined on a common probability space and 
    $$\Pro\left[S_0^{(n)}\rightarrow S_0,\ X^{(n)} \rightarrow X\text{ and for all }k\geq 1,\, \epsilon^{+,(n)}_{k} \rightarrow \epsilon^+_{k},\, \epsilon^{-,(n)}_{k} \rightarrow \epsilon^-_{k}\right] = 1,$$
    where $\epsilon^+_k \sim f^+,$ $\epsilon^-_k \sim f^-$ for $k\in\N$. For all $n\in\N$ and $t\geq 0$, we set
    \begin{align*}
    \tilde{S}^{(n), t} &:= \tilde{\Psi}\left(S^{(n)}(t), \,X^{(n)}(\cdot + t) - X^{(n)}(t), \left(R^{+,(n)}_{l^{(n)}(t)+j}\right)_{j\geq 1}, \left(R^{-,(n)}_{l^{(n)}(t)+j}\right)_{j\geq 1}\right),\\
	\dbtilde{S}^{(n),t} &:= \dbtilde{\Psi}\left(S^{(n)}(t), \,X^{(n)}(\cdot + t) - X^{(n)}(t), \left(R^{+,(n)}_{l^{(n)}(t)+j}\right)_{j\geq 1}, \left(R^{-,(n)}_{l^{(n)}(t)+j}\right)_{j\geq 1}\right).
	\end{align*}
 By Assumption \ref{CB:ass:initialCon} and Theorem \ref{CB:res:convQM-1}, 
    \begin{equation}\label{CB:eq:convS0}
    \tilde{S}^{(n),0} =(\tilde{B}^{(n),0},\tilde{Q}^{(n),0},\tilde{C}^{(n),0})\rightarrow \tilde{S}^{0}=(\tilde{B}^{0},\tilde{Q}^{0},\tilde{C}^{0})\quad \Pro\text{-a.s.}
    \end{equation}
    Note that we have $\sigma^{(n)}_1 =  \sigma(\tilde{C}^{(n),0})$ and $\sigma_1=\sigma(C)$, where  $\sigma(\omega) := \inf\left\{t \geq 0: \omega(t) \notin [-\kappa_-,\kappa_+] \right\}$. Let us set $\overline{\sigma}_1 := \sigma(\tilde{C}^{0})$,   
    $S= \tilde{S}^{0}$ on $[0 ,\overline{\sigma}_1)$, and $C(\overline{\sigma}_1):=C(\overline{\sigma}_1-)$. By definition of $\sigma_1$ and the left-continuity of $C$ at $\overline{\sigma}_1$, we conclude that $\sigma_1=\overline{\sigma}_1$. 

   From Lemma \ref{CB:res:cumQatSwitching} and Theorem \ref{CB:res:idQ-1} we know that in a neighborhood of $\sigma_1$ all but one of the local time processes in \eqref{CB:eq:Ci} stay constant, i.e.~locally $\tilde{C}$ behaves like
   $|d\tilde{C}(t)|=\frac{1}{2}|dL_t(\pi_j^{(2)}g(\pi_i\tilde{W}_k))|$ for some $k\in\N$, $j\in\{1,2\}$, and $i\in\{b,a\}$.
   Hence, on the event $\{\tilde{C}^0(0)\notin\{\kappa_+-\kappa_-\}\}$ the paths of $\tilde{C}^{0}$ take their values in $C_{\kappa_+}'(\R_+,\R)\cap C_{-\kappa_-}'(\R_+,\R)$ $\Pro$-a.s.~and Lemma \ref{CB:res:cont-firstExit-lastValue} yields  $\sigma^{(n)}_1\rightarrow\sigma_1$ $\Pro$-a.s. Define $\tilde{\tau}^+:=\inf\{t>0: \tilde{Q}^{a,F}(t)\tilde{Q}^{b,G}(t)=0\}$ and $\tilde{\tau}^-:=\inf\{t>0: \tilde{Q}^{a,G}(t)\tilde{Q}^{b,F}(t)=0\}$. By Lemma \ref{CB:res:cumQatSwitching} it holds $\Pro\left[\tilde{\tau}^+=\tilde{\tau}^-<\infty\right]=0$.
   On the event $\{\tilde{C}^0(0)=\kappa_+\}\cap\{\tau^+>\tau^-\}$ we have $\sigma_1>\tau^+$ and $\tilde{C}^0(\tilde{\tau}^+)<\kappa_+$ $\Pro$-a.s.~and hence by the strong Markov property we can argue as above and obtain $\sigma_1^{(n)}\rightarrow\sigma_1$ $\Pro$-a.s.
   On the other hand, on the event $\{\tilde{C}^0(0)=\kappa_+\}\cap\{\tau^+<\tau^-\}$ we have $\sigma=\tau^+=\inf\{t>0: Y^{a,F}(t) Y^{b,G}(t)=0\}=\lim_n\inf\{t>0: Y^{a,F,(n)}(t) Y^{b,G,(n)}(t)=0\}$ $\Pro$-a.s.~by Lemma \ref{CB:res:cont-firstExit-lastValue}. Similarly, one can treat the event $\{\tilde{C}^0(0)=-\kappa_-\}$. Altogether we thus obtain that $$\sigma^{(n)}_1\rightarrow\sigma_1\ \Pro\text{-a.s.}$$
   
    By Lemma \ref{CB:res:cumQatSwitching} $\Pro[\sigma_1=\tilde{\tau}_k^*<\infty]=0$ for all $k\in\N$. Hence, $\tilde{Q}$ is almost surely continuous at time $\sigma_1$ and $\tilde{Q}^{(n),0}(\sigma_1^{(n)}-)\rightarrow\tilde{Q}^0(\sigma_1-)$ $\Pro$-a.s. Similarly, $\tilde{C}^{(n),0}(\sigma^{(n)}_1-)\rightarrow\tilde{C}^0(\sigma_1-)$ $\Pro$-a.s. From the proof of Theorem \ref{CB:res:convQM-1} we already know that $\tilde{\tau}_k^{(n)}\rightarrow\tilde{\tau}_k^*$ $\Pro$-a.s.~for all $k\in\N$. Hence, not only $\sigma_1$ avoids all $\tilde{\tau}_k^*$, but for large enough $n$ also $\sigma^{(n)}_1$ avoids all $\tilde{\tau}_k^{(n)}$.
   As $\tilde{B}^0$ is constant on each interval $[\tilde{\tau}^*_{k-1},\tilde{\tau}^*_k)$ and $\tilde{B}^{(n)}$ is constant on each interval $[\tilde{\tau}^{(n)}_{k-1},\tilde{\tau}_k^{(n)})$, we conclude that $\tilde{B}^{(n),0}(\sigma^{(n)}_1-)\rightarrow\tilde{B}^0(\sigma_1-)$ $\Pro$-a.s.
   
    Next define
    \begin{equation*}\label{CB:def:Z}
    Z(x) := \begin{cases}
    \textcolor{white}{+}1 &: x \in (0,\infty) \times \R_-\\
    -1 &: x \in \R_- \times (0,\infty),\\
        \textcolor{white}{+}0 &: \text{otherwise}
    \end{cases}
    \end{equation*}
    and $\tilde{Z}^{I}_1 := Z(\pi_I \tilde{Q}^0(\sigma_1))$ for $I = F,G.$ Now set
    \begin{align*}
        \pi_IQ(\sigma_1) &= \pi_IQ(\sigma_1-)\1_{\{\tilde{Z}^I_1 = 0\}} + \pi_I R^+_{l(\sigma_1-) +1} \1_{\{\tilde{Z}^I_1 = 1\}} + \pi_IR^{-}_{l(\sigma_1-)+1}\1_{\{\tilde{Z}^I_1 = -1\}},\\
        B^{I}(\sigma_1) &= B^{I}(\sigma_1-) + \delta\left(\1_{\{\tilde{Z}^I_1 = 1\}} - \1_{\{\tilde{Z}^I_1 = -1\}}\right).
    \end{align*}
       Let us show that $\tilde{Z}^{I,(n)}_1(\sigma^{(n)}_1) \rightarrow \tilde{Z}^I_1$ $\Pro$-almost surely on the event $\{\sigma_1<\infty\}$ for $I = F,G$, where the random variables $\tilde{Z}^{I,(n)}_1$, $I = F,G$, were defined in \eqref{CB:def:Zn}. 
    From Lemma \ref{CB:res:cumQatSwitching} we know that $\Pro$-a.s.~exactly one component of $\tilde{Q}^0(\sigma_1)$ is in  $\R_-$ while the remaining components are in $(0,\infty).$ Since $\tilde{Q}^{(n),0}(\sigma_1^{(n)}-)\rightarrow\tilde{Q}^0(\sigma_1-)$ $\Pro$-a.s., we conclude that for almost all paths for large $n$ the same component of $\tilde{Q}^{(n),0}(\sigma^{(n)}_1)$ is in $\R_-$, while the others are in $(0,\infty).$ Hence, for $I =F,G$, $\pi_I\tilde{Q}^{(n),0}(\sigma^{(n)}_1)$ and $\pi_I\tilde{Q}^0(\sigma_1)$ are both either in $(0,\infty)^2$, $\R_- \times (0,\infty),$ or $(0,\infty)  \times \R_-$, which implies that $\tilde{Z}^{I,(n)}_1(\sigma^{(n)}_1) \rightarrow \tilde{Z}^I_1$ $\Pro$-almost surely. With the help of Lemma \ref{CB:res:contNT} and Lemma \ref{CB:res:cumQatSwitching} we further conclude that $l^{(n)}(\sigma^{(n)}_1-) \rightarrow l(\sigma_1-)$ $\Pro$-a.s. Hence, Assumption \ref{CB:ass:R} implies that
    \begin{align}\label{CB:eq:convSwitchState1}
        S^{(n)}(\sigma^{(n)}_1)\rightarrow S(\sigma_1) \quad \mathbb{P}\text{-a.s.}
    \end{align}
    In order to extend the definition of $S$ beyond $[0, \sigma_1]$, we define
    $${\rho}(\omega) := \inf\left\{t \geq 0: \left|\pi^{(2)}_1\omega(t) - \pi^{(2)}_2\omega(t)\right| < \frac{\delta}{2}\right\}$$
    and set $\bar{\rho}_1:=(\sigma_1+\rho(\dbtilde{B}^{\sigma_1}))$. Further, we set $S:= \dbtilde{S}^{\sigma_1}(\cdot-\sigma_1)$ on $[\sigma_1, \bar{\rho}_1]$. Then by definition $\bar{\rho}_1=\rho_1$. Moreover, from Theorems \ref{CB:res:convS2} and \eqref{CB:eq:convSwitchState1} we obtain that
    $\dbtilde{S}^{(n), \sigma_1} \rightarrow \dbtilde{S}^{\sigma_1}$ $\Pro$-a.s. Since  $\rho_1^{(n)}=\sigma_1^{(n)}+\rho(\dbtilde{B}^{(n),\sigma^{(n)}_1})+\tn$, we conclude with Lemma \ref{CB:res:cont-firstExit-lastValue} that
    $\rho^{(n)}_{1} \rightarrow {\rho}_{1} \ \Pro\text{-a.s.}$ 
    Note that at time $\rho^{(n)}_1-\tn$ respectively at time $\rho_1$ a price change occurs. Hence, similarly to Step 1 in the proof of Theorem \ref{CB:res:convQM-1} one can show that $S^{(n)}(\rho^{(n)}_1-)\rightarrow S(\rho_1-)$ $\Pro$-a.s. and hence by definition also $S^{(n)}(\rho^{(n)}_1)\rightarrow S(\rho_1)$ $\Pro$-a.s.
    
    Now we can proceed iteratively and construct the process $S$ on $[0,\lim_k\rho_k)$ such that for all $k\in\N$ we have $S=\tilde{S}^{\rho_{k-1}}(\cdot-\rho_{k-1})$ on $[\rho_{k-1},\sigma_k)$ and $S=\dbtilde{S}^{\sigma_{k}}(\cdot-\sigma_{k})$ on $[\sigma_k,\rho_k)$ and
    \begin{align*}
        &\left(\sigma^{(n)}_1, \rho^{(n)}_1, \cdots, \sigma^{(n)}_k, \rho^{(n)}_k, S^{(n)}(\sigma^{(n)}_1), S^{(n)}(\rho^{(n)}_1), \cdots, S^{(n)}(\sigma^{(n)}_k), S^{(n)}(\rho^{(n)}_k)\right)\\
        &\hspace{4cm} \rightarrow \left(\sigma_1, \rho_1,\cdots, \sigma_k, \rho_k, S(\sigma_1), S(\rho_1), \cdots, S(\sigma_k), S(\rho_k)\right)\quad \Pro\text{-a.s.}
    \end{align*}
    It remains to show that $\Pro[\rho_k \rightarrow \infty]=1$: at the time of a regime switch at least the two queues of the country which triggered the regime change get reinitialized. Hence, at time $\rho_{2k}$ one of the countries, say $I$, has triggered at least $k$ regime switches from an inactive to an active regime and its queues have been reinitialized at least $k$ times
    at either $\pi_I R^+_k$ or $\pi_I R^-_k$. By Assumption \ref{CB:ass:R} we know that each $R_k^\pm$ is component-wise bounded from below by  $\alpha \epsilon_k^\pm$, where the $(\epsilon_k^\pm)$ are i.i.d.~and $\alpha>0$. As the starting time $\rho_k$ of an active regime equals the hitting time at zero of one of the components of $Q$, we can thus bound $\rho_{2k}$ from below by the sum of the first hitting times at the axes of $k$ independent planar Brownian motions, each started from $\alpha\cdot\min\{\pi_F\epsilon_k^+,\pi_F\epsilon_k^-,\pi_G\epsilon_k^+,\pi_G\epsilon_k^-\}$. But this sum goes to infinity almost surely as $k\rightarrow\infty$, so that also $\Pro[\rho_k\rightarrow \infty] = 1$.
    
    Altogether, we finally obtain from Theorems \ref{CB:res:convQM-1} and \ref{CB:res:convS2} that $S^{(n)} \rightarrow S$ $\mathbb{P}\text{-a.s.}$
\end{proof}

\section{Conclusion}\label{CB:sec:conclusion}

In this paper we presented a microscopic reduced form model of a cross-border market for two countries with restricted transmission capacities based on a shared LOB. We derived its high frequency dynamics which turn out to be analytically quite tractable: in the limit the dynamics are characterized by a continuous-time regime switching process, whose queue size dynamics between price changes is given by a linear four dimensional Brownian motion in the positive orthant during inactive regimes and a four-dimensional SRBM in the positive orthant during active regimes.

A natural question is whether our analysis can be extended to more than two countries. The main obstacle for such an extension is the definition of the reflection matrix in a multiple countries setting, i.e.~the specification of the cross-border trading mechanism. While our mathematical analysis can easily be extended to a multidimensional setting for a given reflection matrix, the situation becomes more complicated when the reflection matrix is state dependent on the queue size processes and the available capacities. Mathematically, this would lead to a new class of reflected SDEs, in which the coefficients for the reflection term depend on the reflection process itself. Moreover, any such choice of a reflection matrix would still be ad-hoc, while the optimal design of the cross-border trading mechanism in the presence of transmission capacity constraints remains an open and challenging question.

\appendix

\section{Technical details}\label{CB:app:techTools}

\subsection{Proof of Lemma \ref{CB:res:sumCompg}} \label{CB:app:proofs}

\begin{proof}[Proof of Lemma \ref{CB:res:sumCompg}]
$ $
\begin{enumerate}
\item[i)] By definition, we have for any $t\in[\htau_{k-1},\htau_k)$,
\begin{align*}
	g(\omega)(t) = g(\omega)(\htau_{k-1}-) + \omega(t)- \omega(\htau_{k-1}-) + \big(\bar{g}(\omega)(t)-\bar{g}(\omega)(\htau_{k-1}-)\big)\mathcal{R}.
	\end{align*}
Summing over $l=2,\dots,k$ and noting that $g(\omega)(\htau_1-)=\omega(\htau_1-)$, $\bar{g}(\omega)(\htau_1-)=0$, we obtain
\[g(\omega)(t) = g(\omega)(\htau_{1}-) + \omega(t)- \omega(\htau_{1}-) + \big(\bar{g}(\omega)(t)-\bar{g}(\omega)(\htau_{1}-)\big)\mathcal{R}=\omega(t)+\bar{g}(\omega)(t)\mathcal{R}.\]

\item[ii)] Apply $h_1$ to i) and note that $h_1(x\mathcal{R})=0$ for all $x\in\R^2$: \[ h_1(g(\omega))=\pi^{(2)}_1g(\omega)+\pi^{(2)}_2g(\omega)=\pi^{(2)}_1\omega+\pi^{(2)}_2\omega=h_1(\omega).\]

\item[iii)]  Let $\htau_k<\htau_\infty<\infty$ for all $k\in\N$. Without loss of generality, suppose that $\pi_1^{(2)}\omega(\htau_1)=0$. Then $\pi_1^{(2)}g(\omega)(\htau_{2j-1})=0$ and $\pi_2^{(2)}g(\omega)(\htau_{2j})=0$ for all $j\in\N$. Hence, for any $k\in\N$,
    \begin{align*}
        \sup_{t \in [\htau_k,\htau_\infty)} \pi^{(2)}_1 g(\omega)(t) &\leq 
        \sup_{j\geq \lfloor k/2 \rfloor}\sup_{t\in [\htau_{2j}, \htau_{2j+2})} \left|\pi^{(2)}_1 g(\omega)(\htau_{2j+1}) - \pi^{(2)}_1 g(\omega)(t)\right|.
    \end{align*}
    By i) and the definition of $\bar{g}$,
    \begin{align*}
        \sup_{t\in [\htau_{2j}, \htau_{2j+1})}&\left|\pi^{(2)}_1 g(\omega)(\htau_{2j+1}) - \pi^{(2)}_1 g(\omega)(t)\right|\\
        &= \sup_{t\in [\htau_{2j}, \htau_{2j+1})}\left|\pi_1^{(2)}\omega(\htau_{2j+1}) - \pi_1^{(2)}\omega(t) - \pi_2^{(2)}\bar{g}(\omega)(\htau_{2j+1}) + \pi_2^{(2)}\bar{g}(\omega)(t)\right|\\
        &\leq \sup_{t\in [\htau_{2j}, \htau_{2j+1})} \left|\pi_1^{(2)}\omega(\htau_{2j+1})-\pi_1^{(2)}\omega(t)\right| + \sup_{s,t \in [\htau_{2j}, \htau_{2j+1})} \left|\pi_2^{(2)}\omega(t) - \pi_2^{(2)}\omega(s)\right|
    \end{align*}
and
    \begin{align*}
        \sup_{t\in [\htau_{2j+1}, \htau_{2j+2})}&\left|\pi^{(2)}_1 g(\omega)(\htau_{2j+1}) - \pi^{(2)}_1 g(\omega)(t)\right|\\
        &= \sup_{t\in [\htau_{2j+1}, \htau_{2j+2})}\left|\pi_1^{(2)}\omega(\htau_{2j+1}) - \pi_1^{(2)}\omega(t)+\pi_1^{(2)}\bar{g}(\omega)(\htau_{2j+1}) - \pi_1^{(2)}\bar{g}(\omega)(t)\right|\\
        &\leq \sup_{t\in [\htau_{2j+1}, \htau_{2j+2})} \left|\pi_1^{(2)}\omega(\htau_{2j+1})-\pi_1^{(2)}\omega(t)\right| + \sup_{s,t \in [\htau_{2j+1}, \htau_{2j+2})} \left|\pi_1^{(2)}\omega(t)-\pi_1^{(2)}\omega(s)\right|.
    \end{align*}
Since $\omega$ is left-continuous, we conclude that 
    \[\sup_{t \in [\htau_k,\htau_\infty)}\pi^{(2)}_1 g(\omega)(t) \leq 3 \sup_{s,t\in [\htau_{k}, \htau_\infty)}\left|\pi^{(2)}_1 \omega(t) - \pi^{(2)}_1\omega(s)\right|+\sup_{s,t\in [\htau_k, \htau_\infty)}\left|\pi^{(2)}_2 \omega(t) - \pi^{(2)}_2\omega(s)\right| \rightarrow 0\]
as $k\rightarrow\infty$. Similarly, one can show that $\sup_{t \in [\htau_k,\htau_\infty)}\pi^{(2)}_2 g(\omega)(t)\rightarrow0$ as $k\rightarrow\infty$.

\item[iv)] According to \eqref{CB:eq:SP1d} on each interval $[\htau_{k-1},\htau_k)$ one of the components of $\bar{g}$ stays constant, while the other is the unique constraining process of a locally defined one-dimensional Skorokhod problem and hence only charges the zero set of the respective component of $g$. 

\end{enumerate}

\end{proof}

\subsection{Continuity results}\label{CB:app:contProp}

For any $z \in \R$ we define the maps
\begin{align*}
\begin{aligned}
    \tau_z:& \,D(\R_+, \R) \rightarrow \R_+, \qquad &\tau_z(\omega)&:= \inf\{t \geq 0: \omega(t) \leq z\},\\
    \tau'_z:& \, D(\R_+, \R) \rightarrow \R_+, \quad &\tau'_z(\omega) &:= \inf\{t\geq 0: \omega(t) < z\}
\end{aligned}
\end{align*}
and introduce for $k\in\N$ the function spaces 
\begin{equation}\label{CB:def:C'}
\begin{split}
    D'_z\left(\R_+, \R^k\right)&:= \left\{\omega \in D(\R_+, \R^k): \, \tau'_z\left(\pi^{(k)}_j\omega\right)=\tau_z\left(\pi^{(k)}_j\omega\right)\ \forall\ j=1,\dots,k \right\},\\
        C'_z\left(\R_+, \R^k\right)&:= \left\{\omega \in C(\R_+, \R^k): \, \tau'_z\left(\pi^{(k)}_j\omega\right)=\tau_z\left(\pi^{(k)}_j\omega\right)\ \forall\ j=1,\dots,k \right\},
        \end{split}
\end{equation}
which contain all c\`adl\`ag resp.~continuous functions whose components cross $z$ as soon as they hit $z$.  

\begin{Lem}
\label{CB:res:cont-firstExit-lastValue}
     For all $z\in \R$ the functions $\tau'_z$ and $\tau_z$ 
     are continuous at $\omega_0\in D'_z(\R_+, \R)$.
\end{Lem}

\begin{proof}
Since convergence in the Skorokhod $J_1$-topology implies convergence in the Skorokhod $M_2$-topology, we know from \cite[Theorem 13.6.4]{W02} that $\tau_z'$ 
is continuous at $\omega_0$. Now let $(\omega_n)_{n\in \N}\subset D(\R_+,\R)$ be any sequence converging towards $\omega_0$ in the Skorokhod topology and let $(a_n)_{n\in \N}\subset\R_+$ be a strictly positive null sequence. Then $\omega_n^+:=\omega_n+a_n\rightarrow \omega_0$ and $\omega_n^-:=\omega_n-a_n\rightarrow\omega_0$. 
Since $\omega_0 \in D'_z(\R_+, \R),$ we have $\tau'_z(\omega_0) = \tau_z(\omega_0)$. 
Hence, $\tau'_z(\omega_n^\pm)\rightarrow \tau'_z(\omega_0)=\tau_z(\omega_0)$. 
Since $\tau'_z(\omega_n^+)\leq\tau_z(\omega_n)\leq\tau'_z(\omega_n^-)$ for all $n\in\N$, we conclude that $\tau_z(\omega_n)\rightarrow\tau_z(\omega_0)$. 
\end{proof}

Next, we want to characterize the continuity set of the function $g$ introduced in Definition \ref{CB:def:g}. 

\begin{Lem}[Continuity of $g$]
\label{CB:res:contg}
If $\omega_0 \in D(\R_+, \R^2)$ satisfies $h_1(\omega_0) \in C'_0(\R_+,\R)$, then the function $g:D(\R_+,\R^2)\rightarrow D(\R_+,\R^2_+)$ is continuous at $\omega_0.$
\end{Lem}

\begin{proof}
Let us denote by $\Gamma:D(\R_+,\R^2)\rightarrow D(\R_+,\R_+^2)$ the extended Skorokhod map for the two-dimensional GPS ESP, cf.~\cite[Section 3.2]{R06}. As already noted in Section 4.1, Lemma \ref{CB:res:sumCompg} i) and iv) imply that $g=\Gamma$ on $[0,\htau_\infty)$. By \cite[Theorem 3.6]{R06}, $\Gamma$ is Lipschitz continuous with respect to the topology of uniform convergence on compact intervals and hence with respect to the Skorokhod topology. Thus, for any $\omega'\in D(\R_+,\R^2)$ satisfying $\htau_\infty(\omega_0)\geq\htau_\infty(\omega')$ we have for some $L>0$,
\begin{align*}
d_{J_1}\left(g(\omega_0),g(\omega')\right)&\leq d_{J_1}\left(g(\omega_0)\1_{[0, \hat{\tau}_\infty(\omega'))},g(\omega')\1_{[0, \hat{\tau}_\infty(\omega'))}\right)+d_{J_1}\left(g(\omega_0),g(\omega_0)\1_{[0,\htau_\infty(\omega'))}\right)\\
&\leq d_{J_1}\left(\Gamma(\omega_0),\Gamma(\omega')\right)+ \|g(\omega_0)\1_{[\hat{\tau}_\infty(\omega'),\htau_\infty(\omega_0))}\|_\infty\\
&\leq Ld_{J_1}(\omega_0,\omega')+\|h_1(\omega_0)\1_{[\hat{\tau}_\infty(\omega'),\htau_\infty(\omega_0))}\|_\infty,
\end{align*}
where the last inequality follows from the Lipschitz continuity of $\Gamma$ and Lemma \ref{CB:res:sumCompg} ii). 
Note that $\|h_1(\omega_0)\1_{[\hat{\tau}_\infty(\omega'),\htau_\infty(\omega_0))}\|_\infty$ becomes arbitrarily small for $\omega'$ close enough to $\omega_0$, because $\htau_\infty(\omega)=\tau_0(h_1(\omega))$ is continuous at $\omega_0$ by Lemma \ref{CB:res:cont-firstExit-lastValue} and $t\mapsto h_1(\omega_0)(t)$ is continuous with $h_1(\omega_0)(\htau_\infty(\omega_0))=0$. The case $\htau_\infty(\omega_0)<\htau_\infty(\omega')$ is treated analogously noting that 
\begin{align*}
    \|h_1(\omega')\1_{[\hat{\tau}_\infty(\omega_0),\htau_\infty(\omega'))}\|_\infty\leq  \|h_1(\omega')-h_1(\omega_0)\|_\infty+\|h_1(\omega_0)\1_{[\hat{\tau}_\infty(\omega_0),\htau_\infty(\omega'))}\|_\infty,
\end{align*}
which also becomes arbitrarily small for $\omega'$ close enough to $\omega_0$.
\end{proof}

\begin{Cor}\label{CB:res:cont-barg}
    Let $\omega_0 \in D(\R_+, \R^2)$ satisfy $h_1(\omega_0) \in C'_0(\R_+, \R)$. Then $\hat{g}: D(\R_+, \R^2) \rightarrow D(\R_+, \R)$ given by $\hat{g}:= \pi^{(2)}_2 \bar{g} - \pi^{(2)}_1 \bar{g}$ is continuous at $\omega_0.$
\end{Cor}

\begin{proof}
    Thanks to Lemma \ref{CB:res:contg} the function $g$ is continuous at $\omega_0.$ As $\hat{g}(\omega) = \pi^{(2)}_1 \omega - \pi^{(2)}_1g(\omega)$ on $[0,\htau_\infty(\omega))$ and $\hat{g}(\omega)(t)=\hat{g}(\omega)(\htau_\infty(\omega)-)$ by definition of $\bar{g}$ for all $\omega \in D(\R_+, \R^2)$, the continuity of $\hat{g}$ at $\omega_0$ follows from the continuity of $\htau_\infty(\omega)=\tau_0(h_1(\omega))$ at $\omega_0$.
\end{proof}

 The next lemma will be needed to show continuity of the reinitializations.

\begin{Lem}\label{CB:res:contH}
Let $\omega_0\in D(\R_+,\R^4)$ satisfy $h(\omega_0)\in D'_0(\R_+,\R^2)$ and $\tilde{\tau}_b(\omega_0)\neq \tilde{\tau}_a(\omega_0)$. Then $H_i(\omega):=\1_{\left\{\tilde{\tau}(\omega) = \tilde{\tau}_{i}(\omega)\right\}}$ is continuous at $\omega_0$ for $i=b,a$.
\end{Lem}

\begin{proof}
As $\tilde{\tau}_b(\omega_0)\neq \tilde{\tau}_a(\omega_0)$ there is $i\in\{b,a\}$ and $\varepsilon>0$ such that $|\tilde{\tau}_i(\omega_0)-\tilde{\tau}(\omega_0)|>\varepsilon$. Without loss of generality let $i=b$. Since $h(\omega_0)\in D'_0(\R_+,\R^2)$ and $\tilde{\tau}_b(\omega_0)=\tau_0(\pi_1^{(2)}h(\omega_0))$, Lemma \ref{CB:res:cont-firstExit-lastValue} implies that $\tilde{\tau}_b$ and $\tilde{\tau}=\tilde{\tau}_b\wedge\tilde{\tau}_a$ are continuous at $\omega_0$. Hence, there is $\delta=\delta(\varepsilon)>0$ such that $|\tilde{\tau}(\omega_0)-\tilde{\tau}(\omega')|<\frac{\varepsilon}{4}$ and $|\tilde{\tau}_b(\omega_0)-\tilde{\tau}_b(\omega')|<\frac{\varepsilon}{4}$ for any $\omega'\in D(\R_+, \R^4)$ satisfying $d_{J_1}(\omega_0,\omega') < \delta$. But then
\begin{align*}
    \left|\tilde{\tau}_b(\omega')-\tilde{\tau}(\omega')\right|> \left|\tilde{\tau}_b(\omega_0)-\tilde{\tau}(\omega_0)\right|- \left|\tilde{\tau}_b(\omega')-\tilde{\tau}_b(\omega_0)\right|- \left|\tilde{\tau}(\omega')-\tilde{\tau}(\omega_0)\right|>\frac{\varepsilon}{2},
\end{align*}
i.e.~$\tilde{\tau}_b(\omega')\neq\tilde{\tau}(\omega')$. Therefore, $H_i$ is continuous at $\omega_0$ for $i = 1,2$.
\end{proof}

\begin{Lem}[Continuity of $\tilde{\Psi}^Q_k$]
\label{CB:res:contPsik}
Let $(\omega_0, r_+, r_-) \in D(\R_+, \R^4) \times (\R^4_+)^{\N} \times (\R^4_+)^{\N}$ satisfy conditions ii) and iii) of Theorem \ref{CB:res:contPsi}.     
Then for any $k\in\N$ the map $\tilde{\Psi}^Q_k: D(\R_+, \R^4) \times (\R^4_+)^{\N} \times (\R^4_+)^{\N} \rightarrow D(\R_+, \R^4_+)$ is continuous at $(\omega_0, r_+, r_-)$.
 \end{Lem}

\begin{proof}
We only prove the claim for $k=1$ as the rest will easily follow by induction. Let $\delta>0$ and choose $(\omega', r'_{+}, r'_{-}) \in D(\R_+, \R^4) \times (\R^4_+)^\N \times (\R^4_+)^\N$ with $d((\omega_0, r_{+}, r_{-}),(\omega', r'_{+}, r'_{-})) < \delta$. Then conditions ii) and iii) together with Lemma \ref{CB:res:contH} allow us to choose $\delta>0$ small enough such that at times $\tilde{\tau}(\omega_0)$ and $\tilde{\tau}(\omega')$ both processes get either reinitialized by $r_{+,1}$ and $r'_{+,1}$ or by $r_{-,1}$ and $r'_{-,1}$. Next we choose $\lambda\in C(\R_+,\R_+)$ increasing such that $\lambda(0)=0$, $\tilde{\tau}(\omega'\circ\lambda)=\tilde{\tau}(\omega_0)$, $\lim_{T\rightarrow\infty}\lambda(T)=\infty$, and $\|\lambda-\id\|_\infty=|\tilde{\tau}(\omega_0)-\tilde{\tau}(\omega')|$. 
With this choice of $\delta$ and $\lambda$ we obtain
\begin{equation}\label{CB:eq:contPsiQ}
 \begin{split}
		&\tilde{\Psi}^{Q}_1 (\omega', r'_{+,1}, r'_{-,1}) \circ \lambda - \tilde{\Psi}^{Q}_1(\omega_0, r_{+,1}, r_{-,1})=\\
		& \left[G(\omega' \circ \lambda) - G(\omega_0)\right]  \1_{[0, \tilde{\tau}(\omega_0))}  +
  \left[G\Big(\tilde{\Psi}^Q(\omega'\circ\lambda,r'_+,r'_-)(\tilde{\tau}(\omega_0))+\omega'\circ\lambda-(\omega'\circ\lambda)(\tilde{\tau}(\omega_0))\Big)\right.\\
  &\left.\qquad\qquad \qquad\qquad\qquad\qquad\qquad-G\Big(\tilde{\Psi}^Q(\omega_0,r_+,r_-)(\tilde{\tau}(\omega_0))+\omega_0-\omega_0(\tilde{\tau}(\omega_0))\Big)  \right] \1_{[\tilde{\tau}(\omega_0), \infty)}.
\end{split}	
 \end{equation}
With conditions ii) and iii) Lemma \ref{CB:res:contg} and Lemma \ref{CB:res:contH} imply that $G$ is continuous at $\omega_0$ and at
\begin{align*}
&\tilde{\Psi}^Q(\omega_0,r_+,r_-)(\tilde{\tau}(\omega_0))+\omega_0(\cdot+\tilde{\tau}(\omega_0))-\omega_0(\tilde{\tau}(\omega_0))=\\
&\1_{\left\{\tilde{\tau}(\omega_0)=\tilde{\tau}_a(\omega_0)\right\}} r_{+,1} 
		 + \1_{\left\{\tilde{\tau}(\omega_0)  = \tilde{\tau}_b(\omega_0)\right\}} r_{-,1} + \omega_0(\cdot+\tilde{\tau}(\omega_0)) - \omega_0(\tilde{\tau}(\omega_0)).
\end{align*} 
Thus, as $\|\lambda-\id\|_\infty=|\tilde{\tau}(\omega_0)-\tilde{\tau}(\omega')|$ and $\tilde{\tau}$ is continuous at $\omega_0$, the RHS of \eqref{CB:eq:contPsiQ} will get arbitrarily small, if we choose $\delta$ small enough.
\end{proof}

\noindent 
Next we study the continuity set of $N: D(\R_+, \R^4) \times (\R^4_+)^{\N} \times (\R^4_+)^{\N} \times \R_+ \rightarrow \N_0 \cup \{+\infty\}$ given by
\begin{align*}\label{CB:def:generalN}
	N(\omega, r^+, r^-, T) := \inf\left\{k\geq 0: \tilde{\tau}(\tilde{\Psi}^Q_k(\omega, r^+, r^-)) > T\right\},
\end{align*}
where we endow $D(\R_+, \R^4) \times (\R^4_+)^{\N} \times (\R^4_+)^{\N} \times \R_+ \rightarrow \N_0 \cup \{+\infty\}$ again with the product topology.

\begin{Lem}
\label{CB:res:contNT}
	Let $(\omega_0, r_+, r_-, T) \in D(\R_+, \R^4) \times (\R^4_+)^{\N} \times (\R^4_+)^{\N} \times \big(\R_+\backslash\{\tilde{\tau}_k:\ k\in\N_0\}\big)$. If the conditions of Theorem \ref{CB:res:contPsi} are satisfied, then 
 $N$ is continuous at $(\omega_0, r_+,r_-, T)$.
\end{Lem}

\begin{proof}

Take $\delta > 0$ and let $(\omega', r'_+, r'_-, T') \in D(\R_+, \R^4) \times (\R^4_+)^{\N} \times (\R^4_+)^{\N} \times \R_+$ satisfy $$d((\omega_0, r_+, r_-, T),(\omega', r'_+, r'_-, T')) < \delta.$$ 
 Denote $N_0:=N(\omega_0 , r_+, r_-,T),\ N':=N(\omega' , r'_+, r'_-,T')$ and 
  $\tilde{\tau}'_0 := 0,\ \tilde{\tau}'_k := \tilde{\tau}(\tilde{\Psi}^Q_{k-1}(\omega', r'_+, r'_-))$ for $k\in\N.$
By assumption $N_0<\infty$ and $\tilde{\tau}_{N_0}<T<\tilde{\tau}_{N_0+1}$. Lemmata \ref{CB:res:cont-firstExit-lastValue} and \ref{CB:res:contPsik} imply that $\tilde{\tau}_k'=\tilde{\tau}\circ\tilde{\Psi}^Q_{k-1}$ is continuous at $(\omega_0,r_+,r_-)$ for all $k\in\N$. Hence, if $\delta>0$ is small enough, then also $\tilde{\tau}'_{N_0}<T'<\tilde{\tau}'_{N_0+1}$, i.e.~$N_0=N'$.
\end{proof}

\noindent
Similarly, the following continuity result can be deduced.

\begin{Lem}
\label{CB:res:contNba}
	Let $(\omega_0, r^+, r^-) \in D(\R_+, \R^4) \times (\R^4_+)^{\N} \times (\R^4_+)^{\N}$ satisfy the conditions of Theorem \ref{CB:res:contPsi}. Then the maps $N_b$ and $N_a$ introduced in Definition \ref{CB:def:PsiQ1} 
 are continuous at $(\omega_0, r^+, r^-).$
\end{Lem}
 
 \begin{proof}[Proof of Theorem \ref{CB:res:contPsi}] 
 For all $k\geq 0$ set
 \[\tilde{Y}_k(\omega_0, r^+, r^-) :=\left(h \circ \tilde{\Psi}^{Q}_k\right)(\omega_0,r^+, r^-)(\cdot + \tilde{\tau}_{k}) = h\left(\tilde{\Psi}^{Q}_{k}(\omega_0, r^+, r^-)(\tilde{\tau}_{k}) + \omega_0(\cdot + \tilde{\tau}_{k}) - \omega_0(\tilde{\tau}_{k})\right).\]
Note that the jump times of $\tilde{\Psi}^{Q}(\omega_0,r^+, r^-)$ are equal to the first hitting times of the axes by $\tilde{Y}_k(\omega_0,r^+,r^-)$, $k \geq 0$. By condition i), for every $T>0$ there is $N_T\in\N$ such that $\tilde{\Psi}^{Q}(\omega_0, r^+, r^-) = \tilde{\Psi}^{Q}_{N_T}(\omega_0, r^+, r^-)$ on $[0,T]$. By Lemma \ref{CB:res:contNT} the function $N_T$ is continuous in $(\omega_0, r^+, r^-)$ for all $T\notin\{\tilde{\tau}_k:k\in\N_0\}$. Together with Lemma \ref{CB:res:contPsik} this yields the continuity of $\tilde{\Psi}^Q$ at $(\omega_0,r^+,r^-)$.\par

Moreover, condition ii) together with Corollary \ref{CB:res:cont-barg} implies that for all $k\in\N_0$, $\bar{G}$ is continuous at $\tilde{\Psi}^Q(\omega_0, r^+, r^-)(\tilde{\tau}_k) + \omega_0(\cdot + \tilde{\tau}_k) - \omega_0(\tilde{\tau}_k)$. Hence, by condition i) on every compact interval $\tilde{\Psi}^C$ is the sum of finitely many functions that are continuous at $(\omega_0, r^+, r^-)$ and therefore $\tilde{\Psi}^C$ is itself continuous at $(\omega_0, r^+, r^-).$ 

By Lemma \ref{CB:res:contNba} the maps $N_b, N_a$ are continuous at $(\omega_0, r^+, r^-).$  Condition iii) implies that
\begin{equation*}\label{CB:cond:noSimJumps}
\Disc\left(N_a(\omega_0, r^+, r^-)\right) \cap \Disc\left(N_b(\omega_0, r^+, r^-)\right) = \emptyset,
\end{equation*}
where $\Disc(N_i(\omega_0, r^+, r^-)) := \{t\geq0: N_i(\omega_0, r^+, r^-)(t-) \neq N_i(\omega_0, r^+, r^-)(t)\}$ for $i = b,a.$ Hence, \cite[Theorem 4.1]{W80} allows us to conclude that $\tilde{\Psi}^B$ is continuous at $(\omega_0, r^+, r^-).$ 
\end{proof}

\subsection{Auxiliary results for the limiting dynamics}\label{CB:app:SRBM}

\begin{proof}[Proof of Proposition \ref{CB:res:idgW}]
By definition $(Z,l)$ is a sum-conserving SRBM with absorption associated with $W$ if and only if $(Z,l\mathcal{R})$ solves the GPS (E)SP for $W$ up to time $\tau := \inf\{t\geq 0: Z(t) = (0,0)\}$ and $(Z,l)$ is constant from time $\tau$ onwards. We have already seen above as a consequence of Lemma \ref{CB:res:sumCompg} that $\htau_\infty(W)=\inf\{t\geq0:\ g(W)(t)=(0,0)\}$ and that $(g(W),\bar{g}(W)\mathcal{R})$ is a solution to the GPS (E)SP on $[0,\htau_\infty(W))$. Hence, $g(W)$ is a sum-conserving SRBM with absorption associated with $W$. Moreover, by Definition \ref{CB:def:g} there exists for any $k\in\N$ an $i\in\{1,2\}$ such that on $[\htau_{k-1},\htau_k)$, 
\begin{align*}
    d\pi^{(2)}_ig(W)(t)&=dW_i(t)+d\pi^{(2)}_i\bar{g}(W)(t),\\
    \pi^{(2)}_ig(W)(t)d\pi^{(2)}_i\bar{g}(W)(t)&=0,\\
    \pi_j^{(2)}g(W)(t) &> 0,\\
    d\pi^{(2)}_j\bar{g}(W)(t)&=0 \quad\text{for}\quad j\neq i.
\end{align*}
As $\pi^{(2)}_i\bar{g}(W)$ is non-decreasing, we can thus apply \cite[Lemma 2.12]{Z17} to conclude that for $t \in [\htau_{k-1}, \htau_k)$,
    \begin{equation}\label{CB:eq:localtime}
        \pi^{(2)}_i\bar{g}(W)(t) - \pi^{(2)}_i\bar{g}(W)(\htau_{k-1}) = \frac{1}{2}\left(L_t(\pi^{(2)}_ig(W)) - L_{\htau_{k-1}}(\pi^{(2)}_ig(W))\right).
    \end{equation}
As neither $\pi^{(2)}_j\bar{g}(W)$ nor $L(\pi^{(2)}_jg(W))$ increases on $[\htau_{k-1},\htau_k)$, equation \eqref{CB:eq:localtime} is also valid for {$j\neq i$ instead of $i$.}
\end{proof}

\begin{Lem}\label{CB:res:cumQatSwitching}
For all $t\geq0$, on the event $\{t\neq\tilde{\tau_k^*}\,\forall\, k\in\N\}$ at most one component of $\tilde{Q}(t)$ is zero $\Pro$-a.s. Moreover, $\Pro[\tilde{\tau}^*_k=\tilde{\sigma}_1<\infty]=0$ for all $k\in\N$ and $\tilde{\sigma}_1:=\inf\{t \geq 0: \tilde{C}(t) \notin[-\kappa_-,\kappa_+]\}$.
\end{Lem}

\begin{proof}
By the strong Markov property, it suffices to prove the claim on the interval $[0,\tilde{\tau}_1^*]$. Let us denote by $(\check{\tau}_k)_{k\geq 0}$ the hitting times of $\tilde{Q}$ at alternating axes, i.e.~$\check{\tau}_0:=0$ and for $k\in\N$,
\[\check{\tau}_k:=\inf\left\{t>\check{\tau}_{k-1}:\exists\ i\in\{1,\dots,4\} \text{ s.t. } \pi_i\tilde{Q}(t)=0\neq\pi_i\tilde{Q}(\check{\tau}_{k-1})\right\}\quad \text{and}\quad\check{\tau}_\infty:=\lim_{k\rightarrow\infty}\check{\tau}_k.\]

First, we show that throughout the interval $[0,\check{\tau}_\infty)$ two components of $\tilde{Q}$ are never simultaneously zero $\Pro$-almost surely: on the interval $[0,\check{\tau}_1]$, the process $\tilde{Q}$ is a four-dimensional linear Brownian motion and hence does not hit any two-dimensional subspace $\Pro$-a.s. On the interval $[\check{\tau}_1,\check{\tau}_2]$, the process $\tilde{Q}$ behaves like a diffusion on a half-space. By \cite[Theorem 3.3]{R88}, this process does not hit any two-dimensional subspace of $\R_+^4$ $\Pro$-a.s. Hence, the claim follows by induction. 

Second, we show that $\Pro[\check{\tau}_\infty=\tilde{\tau}^*_1]=1$: on the event $\{\check{\tau}_\infty<\tilde{\tau}^*_1\}$, we may assume without loss of generality that $\tilde{Q}^{b,I}$ and $\tilde{Q}^{a,J}$ do not hit zero throughout the interval $[0,\check{\tau}_\infty]$ for some $I,J\in\{F,G\}$. In this case, on $[0,\check{\tau}_\infty)$ the queue size dynamics of the other two queues are given by ($H\neq I,K\neq J$)
\begin{align*}
d\tilde{Q}^{b,H}(t)=dX^{b,H}(t)+\frac{1}{2}dL_t(\tilde{Q}^{b,H}),\qquad 
d\tilde{Q}^{a,K}(t)=dX^{a,K}(t)+\frac{1}{2}dL_t(\tilde{Q}^{a,K}),
\end{align*}
i.e.~the other two queues behave like a 
SRBM in the positive orthant with normal reflections at the axes. By \cite[Theorem 2.2]{VW85}, this process does not hit the origin almost surely if $\mu^{b,I}=\mu^{a,I}=0$. By a Girsanov type argument, this is also true for $(\mu^{b,I}, \mu^{a,I})\neq (0,0)$. Hence, $\Pro[\check{\tau}_\infty<\tilde{\tau}^*_1]=0$.

Third, we show that $\Pro[\tilde{\tau}^*_1=\tilde{\sigma}_1<\infty]=0$: recall that $Y:=\tilde{Q}_0+X$ is a linear Brownian motion and that  $\Pro[\tilde{\tau}_{b,1}=\tilde{\tau}_{a,1}]=0$, where $\tilde{\tau}_{b,1},\tilde{\tau}_{a,1}$ were defined in \eqref{CB:eq:hitab}. On the event $\{\tilde{\tau}_1^*=\tilde{\tau}_{b,1}\}$, we may assume without loss of generality that $\tilde{Q}^{a,F}$ does not hit zero throughout the interval $[0,\tilde{\tau}_{b,1}]$. In this case, we have $\frac{1}{2}L(\tilde{Q}^{a,G})=\Gamma_1(Y^{a,G})-Y^{a,G}=\sup_{s\leq\cdot}[-Y^{a,G}_s]^+$ and  by Theorem \ref{CB:res:idQ-1} $\tilde{C}$ satisfies 
\begin{align*}
   \tilde{C}= \tilde{C}(0)+\frac{1}{2}L_{\cdot}\left(\tilde{Q}^{b,G}\right)-\frac{1}{2}L_{\cdot}\left(\tilde{Q}^{b,F}\right) - \frac{1}{2}L_{\cdot}\left(\tilde{Q}^{a,G}\right)\quad \text{on}\quad [0, \tilde{\tau}^*_{1}].
\end{align*}
By Assumption \ref{CB:ass:prob}, there are $\alpha,\beta,\gamma\in\R$ such that
\begin{align*}
    Y^{a,G}(t)=\tilde{Q}^{a,G}_0+\left(\mu^{a,G}-\alpha\mu^{b,G}-\beta\mu^{b,F}\right)t+\alpha Y^{b,G}(t)+\beta Y^{b,F}(t)+\gamma Y^\perp(t),\quad t\geq0,
\end{align*}
for a standard Brownian motion $Y^\perp$ independent of $\pi_bY$. Thus, for any finite $\sigma(\pi_b\tilde{Q}(t): t\geq0)$- measurable random time $\eta$ the conditional law of 
$\sup_{s\leq\eta}\left[-Y^{a,G}_s\right]^+$ given $\pi_b\tilde{Q}$ only has an atom at zero and hence for all $\sigma(\pi_b\tilde{Q}(t): t\geq0)$-measurable $x\neq\tilde{C}(0)+\frac{1}{2}L_\eta\left(\tilde{Q}^{b,G}\right)-\frac{1}{2}L_\eta\left(\tilde{Q}^{b,F}\right)$,
\begin{align*}   
    \Pro\left[\tilde{C}(\eta) = x\, \left|\, \pi_b \tilde{Q}\right.\right] 
    &=\Pro\left[\sup_{s\leq\eta}\left[-Y^{a,G}_s\right]^+ = \tilde{C}(0)+\frac{1}{2}L_\eta\left(\tilde{Q}^{b,G}\right)-\frac{1}{2}L_\eta\left(\tilde{Q}^{b,F}\right)-x\, \left|\, \pi_b \tilde{Q}\right.\right] = 0.
\end{align*}
On the event $\{\tilde{\tau}_{b,1}=\tilde{\tau}_1^*<\infty\}$ we have $\tilde{Q}^{b,G}(\tilde{\tau}_{b,1})=\tilde{Q}^{b,F}(\tilde{\tau}_{b,1})=0$ and thus 
$$L_{\tilde{\tau}_{b,1}}(\tilde{Q}^{b,G})-L_{\tilde{\tau}_{b,1}}(\tilde{Q}^{b,F})=Y^{b,G}(\tilde{\tau}_{b,1})-Y^{b,F}(\tilde{\tau}_{b,1})=2Y^{b,G}(\tilde{\tau}_{b,1}),$$
which is independent of $\tilde{C}(0)$ and has diffusive law as $|\rho^{(b,G),(b,F)}|\neq1$ by Assumption \ref{CB:ass:prob}. Therefore, 
$$\Pro\left[\tilde{C}(0)+\frac{1}{2}L_{\tilde{\tau}_{b,1}}\left(\tilde{Q}^{b,G}\right)-\frac{1}{2}L_{\tilde{\tau}_{b,1}}\left(\tilde{Q}^{b,F}\right)\in\{-\kappa_-,\kappa_+\}\right]=0.$$ Since $\tilde{\tau}_{b,1}$ is $\pi_b \tilde{Q}$-measurable, we may conclude that
\begin{align*}
    \Pro\left[\left.\tilde{\sigma}_1 = \tilde{\tau}_{b,1} =\tilde{\tau}^*_{1}<\infty \, \right|\, \pi_b \tilde{Q} \right]\leq\Pro\left[\left.\left\{\tilde{C}(\tilde{\tau}_{b,1}) \in \{-\kappa_-, \kappa_+\} \right\} \cap \left\{\tilde{\tau}^*_{1} = \tilde{\tau}_{b,1}<\infty\right\}\, \right|\, \pi_b \tilde{Q} \right]= 0.
\end{align*}
Thus, $\Pro[\tilde{\sigma}_1 = \tilde{\tau}_{b,1}=\tilde{\tau}^*_{1}<\infty ]=0$. Analogously it follows that also $\Pro[\tilde{\sigma}_1 = \tilde{\tau}_{a,1}=\tilde{\tau}^*_{1}<\infty ]=0$.   
\end{proof}

\begin{The}\label{CB:thm:2dBM}
    Let $W=(W_1,W_2)$ be a planar Brownian motion starting from $x=(x_1,x_2)\in\R_+^2$ with drift $\bar\mu=(\mu_1,\mu_2)\in\R^2$ and positive definite covariance matrix 
    \begin{equation*}
    \bar{\Sigma}=\begin{pmatrix}\sigma_1^2&\rho\sigma_1\sigma_2\\\rho\sigma_1\sigma_2&\sigma_2^2\end{pmatrix}\in\R^{2\times2}.
    \end{equation*}
Denote by $\tau:=\inf\{t\geq0: W_1(t)W_2(t)=0\}$ the first hitting time of the axes. Then, 
\begin{equation}\label{CB:eq:hitting}
\Pro[\tau > t]=k(x,\bar{\mu},\bar{\Sigma};t),
\end{equation}
where 
\begin{equation}\label{CB:def:condDist-Duration}
k(x,\bar{\mu},\bar{\Sigma};t):= \frac{2}{\alpha t}\exp\left(a_1 x_1 + a_2 x_2 + a_t t- \frac{U}{2t}\right)\sum_{j = 1}^{\infty} \sin\left(\frac{j\pi \theta_0}{\alpha}\right)\int_0^{\alpha} \sin\left(\frac{j\pi \theta}{\alpha}\right) k_j(\theta) d\theta,
\end{equation}
and
\begin{align*}
    k_j(\theta) &:= \int_0^{\infty} r \exp\left(-\frac{r^2}{2t}\right)\exp\left(-d_1 r \sin(\alpha-\theta) - d_2 r \cos(\alpha-\theta)\right) I_{j\pi/\alpha}\left(\frac{r \sqrt{U}}{t}\right) dr,\\
	\alpha &:= \begin{cases}
	\pi + \arctan\left(-\frac{\sqrt{1-\rho^2}}{\rho}\right) &: \rho > 0\\
	\frac{\pi}{2} &: \rho = 0\\
	\arctan\left(-\frac{\sqrt{1-\rho^2}}{\rho}\right) &: \rho < 0
	\end{cases},\quad 
	\theta_0 := \begin{cases}
	\pi + \arctan\left(\frac{x_2\sigma_1\sqrt{1-\rho^2}}{x_1\sigma_2-\rho x_2\sigma_1}\right) &: x_1\sigma_2 < \rho x_2\sigma_1\\
	\frac{\pi}{2} &: x_1\sigma_2 = \rho x_2\sigma_1\\
	\arctan\left(\frac{x_2\sigma_1\sqrt{1-\rho^2}}{x_1\sigma_2-\rho x_2\sigma_1}\right) &: x_1\sigma_2 > \rho x_2\sigma_1
	\end{cases}
\end{align*}
with parameters
\begin{align*}
	a_1 &:= \frac{\rho \mu_2 \sigma_1 - \mu_1 \sigma_2}{(1-\rho^2)\sigma^2_1 \sigma_2},  & a_2&:= \frac{\rho \mu_1 \sigma_2 - \mu_2 \sigma_1}{(1-\rho^2)\sigma^2_2 \sigma_1},\\
	d_1&:= a_1 \sigma_1 + \rho a_2 \sigma_2, & d_2&:=a_2 \sigma_2 \sqrt{1-\rho^2},\\
U&:= \frac{1}{1-\rho^2}\left(\frac{x_1^2}{\sigma^2_1} + \frac{x_2^2}{\sigma^2_2} -\frac{2 \rho x_1 x_2}{\sigma_1 \sigma_2}\right),
& a_t &:= \frac{a_1^2\sigma_1^2}{2} + \rho a_1 a_2 \sigma_1 \sigma_2 + \frac{a_2^2 \sigma_2^2}{2} + a_1 \mu_1 + a_2 \mu_2
\end{align*}
and $I_j$ denotes the $j$-th modified Bessel function of first kind.\vskip6pt

Moreover, $\tau<\infty$ $\mathbb{P}$-almost surely if and only if $a_t\leq0$ and in this case for any $z,t>0$,
    \begin{equation*}
        \Pro[W_1(\tau)\in dz,\tau\in dt]=\exp\left(-z\left(\frac{\mu_1}{\sigma_1^2}-\frac{\rho d_2}{\sigma_1}\right)-a_tt\right)p(z,t)dtdz,
    \end{equation*}
    where $p(z,t)$ is the joint density of $W_1(\tau)$ and $\tau$ on $(0,\infty)^2$  in the driftless case, cf.~\eqref{CB:def:jointdensitytau}. Especially, if $\mu_1=\mu_2=0$, then
\[\Pro[W_1(\tau)\in dz]=\frac{1}{2z\alpha}\cdot\frac{\sin\left(\frac{\pi\theta_0}{\alpha}\right)dz}{\cosh(\frac{\pi}{\alpha}\ln(\frac{\sqrt{U}}{z}))-\cos(\frac{\pi\theta_0}{\alpha})}\qquad\text{and}\qquad \Pro[W_1(\tau)>0]=\frac{\alpha-\theta_0}{\alpha}.\]
\end{The}

\begin{proof}
Formula \eqref{CB:eq:hitting} was derived in \cite[Theorem 3.5.2(ii)]{R94}. Integrating the Green function against a more general initial condition, one can even show that for any bounded Borel function $f:\R_+^2\rightarrow\R$ and for all $t\geq0$,
\begin{equation}\label{CB:generalDist}
  \E\left[f(W_t)\1_{\{\tau>t\}}\right]=\frac{2}{\alpha t}\exp\left(a_1 x_1 + a_2 x_2 + a_t t- \frac{U}{2t}\right)\sum_{j = 1}^{\infty} \sin\left(\frac{j\pi \theta_0}{\alpha}\right)\int_0^{\alpha} \sin\left(\frac{j\pi \theta}{\alpha}\right) f_j(\theta) d\theta,
\end{equation}
where
\begin{align*}
	f_j(\theta) &:= \int_0^{\infty} r \exp\left(-\frac{r^2}{2t}\right)\exp\left(-d_1 r \sin(\alpha-\theta) - d_2 r \cos(\alpha-\theta)\right)\\
	&\qquad\qquad\cdot f\big(\sigma_1r\sin(\alpha-\theta),\sigma_2r(\rho\sin(\alpha-\theta)+\sqrt{1-\rho^2}\cos(\alpha-\theta))\big) I_{j\pi/\alpha}\left(\frac{r \sqrt{U}}{t}\right) dr.
\end{align*}
   Letting $t\rightarrow\infty$ in \eqref{CB:eq:hitting}, it follows that $\tau$ is finite $\mathbb{P}$-almost surely if and only if $a_t\leq 0$.\\
    
    To derive the joint distribution of $W_1(\tau)$ and $\tau$, first suppose that $\mu=(0,0)$. In this case, for any $z>0$ the function $u(x_1,x_2):=\Pro[W_1(\tau)>z]$ satisfies the differential equation 
    $$ \frac{\sigma_1^2}{2}u_{x_1x_1}+\rho\sigma_1\sigma_2 u_{x_1x_2}+\frac{\sigma_2^2}{2}u_{x_2x_2}=0$$
    with boundary conditions $u(x,0)=\1_{(z,\infty)}(x)$ and $u(0,x)=0$ for all $x>0$. By a change of variables to polar coordinates for a standard wedge as in \cite{R94}, this equation can be transformed into a standard Laplace equation for $f(\theta,r):=u(x_1,x_2)$, i.e.
    $$ \frac{1}{2}f_{rr}+\frac{1}{2r}f_r+\frac{1}{2r^2}f_{\theta\theta}=0$$
    with boundary conditions $f(r,\alpha)=0$ and $f(r,0)=\1_{(z,\infty)}(r)$ for all $r>0$. 
    The wedge $\{(\theta,r):0<\theta<\alpha,\ r>0\}$ can be transformed into the doubly infinite strip $\{(\xi,\eta):\ 0<\xi<\pi,\eta\in\R\}$ by setting $\xi=\frac{\pi}{\alpha}\theta,\ \eta=\frac{\pi}{\alpha}\ln(r)$, 
    in which case the solution of the heat equation for steady temperature with $\xi=0$ maintained at temperature $\1_{(\frac{\pi}{\alpha}\ln(z),\infty)}(\eta)$ and $\xi=\pi$ at zero temperature is known (cf.~\cite[p.~166]{C59}) to be
    \[v(\xi,\eta)=\frac{\sin(\xi)}{2\pi}\int_\R\frac{\1_{(\frac{\pi}{\alpha}\ln(z),\infty)}(y)dy}{\cos(\pi-\xi)+\cosh(\eta-y)}.\]
    After re-substitution this gives
    \begin{equation*}
        \Pro[W_1(\tau)\in dz]=\frac{1}{2z\alpha}\cdot\frac{\sin\left(\frac{\pi\theta_0}{\alpha}\right)dz}{\cosh\big(\frac{\pi}{\alpha}\ln(\frac{\sqrt{U}}{z})\big)-\cos(\frac{\pi\theta_0}{\alpha})}
    \qquad\text{and}\qquad      
        \Pro[W_1(\tau)>0]
        =\frac{\alpha-\theta_0}{\alpha}.
    \end{equation*} 
    As $W$ is a homogeneous Markov process, one obtains for all $t,z>0$,
    \begin{equation}\label{CB:def:jointdensitytau}
        \begin{split}
        \Pro[W_1(\tau)\in dz,\tau>t]&=\E\left[\1_{\{\tau>t\}}\Pro[W_1(\tau)\in dz|W_t]\right]\\
        &=\frac{dz}{2z\alpha}\cdot\E\left[\frac{\1_{\{\tau>t\}}\sin\left(\frac{\pi\theta_0(W_t)}{\alpha}\right)}{\cosh\big(\frac{\pi}{\alpha}\ln(\frac{\sqrt{U(W_t)}}{z})\big)-\cos(\frac{\pi\theta_0(W_t)}{\alpha})}\right]=:\int_t^\infty p(z,u)dudz,
        \end{split}
    \end{equation}
    which can be explicitly computed using \eqref{CB:generalDist}. The case $\mu\neq(0,0)$ follows by Girsanov transformation.
\end{proof}

\bibliographystyle{plain}
\bibliography{lit}

\begin{thebibliography}{10}

\bibitem{SIDC}
{SIDC Stakeholder Report, November 2022}.
\newblock Available at https://www.nemo-committee.eu/
  assets/files/SIDC_stakeholder_report_1122-6ca6d8241c5028f54bc1da75d74e9b4d.pdf.

\bibitem{B99}
P.~Billingsley.
\newblock {\em {Convergence of probability measures}}.
\newblock Wiley Series in Probability and Statistics: Probability and
  Statistics. John Wiley \& Sons Inc., second edition, 1999.

\bibitem{C59}
H.S. Carslaw and J.C. Jaeger.
\newblock {\em Conduction of heat in solids.}
\newblock Clarendon Press, Oxford, second edition, 1959.

\bibitem{CFVS22}
\'{A}. Cartea, M.~Flora, T.~Vargiolu, and G.~Slavov.
\newblock Optimal cross-border electricity trading.
\newblock {\em SIAM Journal on Financial Mathematics}, 13(1):262--294, 2022.

\bibitem{CJQ19}
\'A. Cartea, S.~Jaimungal, and Z.~Qin.
\newblock Speculative trading of electricity contracts in interconnected
  locations.
\newblock {\em Energy Economics}, 79:3--20, 2019.

\bibitem{CM91}
H.~Chen and A.~Mandelbaum.
\newblock {Stochastic Discrete Flow Networks: Diffusion Approximations and
  Bottlenecks}.
\newblock {\em The Annals of Probability}, 19(4):1463 -- 1519, 1991.

\bibitem{CB20}
T.S. Christensen and F.E. Benth.
\newblock Modelling the joint behaviour of electricity prices in interconnected
  markets.
\newblock {\em Quantitative Finance}, 20(9):1441--1456, 2020.

\bibitem{CL21}
R.~Cont and A.~de~Larrard.
\newblock {Order book dynamics in liquid markets: limit theorems and diffusion
  approximations}.
\newblock arXiv e-print 1202.6412, 2012.

\bibitem{CL13}
R.~Cont and A.~de~Larrard.
\newblock {Price dynamics in a Markovian limit order market}.
\newblock {\em SIAM J. Financial Math}, 4(1):1--25, 2013.

\bibitem{ContMueller}
R.~Cont and M.~Müller.
\newblock A stochastic partial differential equation model for limit order book
  dynamics.
\newblock {\em SIAM J. Finan. Math.}, 12(2):744--787, 2021.

\bibitem{CST10}
R.~Cont, S.~Stoikov, and R.~Talreja.
\newblock A stochastic model for order book dynamics.
\newblock {\em Operations Research}, 58(3):549--563, 2010.

\bibitem{DR99}
P.~Dupuis and K.~Ramanan.
\newblock Convex duality and the skorokhod problem. ii.
\newblock {\em Probab Theory Relat Fields}, 115:197--236, 1999.

\bibitem{EFH21}
P.A. Ernst, S.~Franceschi, and D.~Huang.
\newblock {Escape and absorption probabilities for obliquely reflected Brownian
  motion in a quadrant}.
\newblock {\em Stochastic Process. Appl.}, 142:634--670, 2021.

\bibitem{FR22}
S.~Franceschi and K.~Raschel.
\newblock {A dual skew symmetry for transient reflected Brownian motion in an
  orthant}.
\newblock {\em Queueing Systems}, 102:123--141, 2022.

\bibitem{GD18}
X.~Gao and S.J. Deng.
\newblock Hydrodynamics limit of order book dynamics.
\newblock {\em Probability in the Engineering and Informational Sciences},
  32(1):96–125, 2018.

\bibitem{HPW16}
S.~Hagemann, C.~Pape, and C.~Weber.
\newblock {Are fundamentals enough? Explaining price variations in the German
  day-ahead and intraday power market}.
\newblock {\em Energy Economics}, 54:376--387, 2016.

\bibitem{HR81}
J.M. Harrison and M.I. Reiman.
\newblock {Reflected Brownian motion on an orthant}.
\newblock {\em Ann. Appl. Probab.}, 9(2):302--308, 1981.

\bibitem{HorstKreherFluid}
U.~Horst and D.~Kreher.
\newblock A weak law of large numbers for a limit order book model with fully
  state dependent order dynamics.
\newblock {\em SIAM J. Financial Math.}, 8:314--343, 2017.

\bibitem{HorstKreherDiffusion}
U.~Horst and D.~Kreher.
\newblock A diffusion approximation for limit order book models.
\newblock {\em Stochastic Processes and their Applications}, 129:4431--4479,
  2019.

\bibitem{HP17}
U.~Horst and M.~Paulsen.
\newblock {A law of large numbers for limit order books}.
\newblock {\em Mathematics of Operations Research}, 42(4):1280--1312, 2017.

\bibitem{K21}
Olav Kallenberg.
\newblock {\em Foundations of Modern Probability}, volume~99 of {\em
  Probability Theory and Stochastic Modelling}.
\newblock Cham, third edition, 2021.

\bibitem{KP17}
R.~Kiesel and F.~Paraschiv.
\newblock Econometric analysis of 15-minute intraday electricity prices.
\newblock {\em Energy Economics}, 64:77--90, 2017.

\bibitem{KM20}
D.~Kreher and C.~Milbradt.
\newblock {Jump diffusion approximation for the price dynamics of a fully state
  dependent limit order book model}.
\newblock {\em SIAM J. Financial Math.}, 14(1):1--51, 2023.

\bibitem{L03}
H.~Luckock.
\newblock {A steady state model of the continuous double auction}.
\newblock {\em Quantitative Finance}, 3(5):385--404, 2003.

\bibitem{MTW16}
P.~Markowich, J.~Teichmann, and M.~Wolfram.
\newblock {Parabolic free boundary price formation models under market size
  fluctuations}.
\newblock {\em Multiscale Modeling $\&$ Simulation}, 14(4):1211--1237, 2016.

\bibitem{PhDCM}
C.~Milbradt.
\newblock Stochastic modelling of intraday eletricity markets.
\newblock Ph.D.~Thesis -- Humboldt-Universit\"at zu Berlin, 2023.

\bibitem{MPW08}
S.~Mudchanatongsuk, J.A. Primbs, and W.~Wong.
\newblock Optimal pairs trading: A stochastic control approach.
\newblock In {\em 2008 American Control Conference}, pages 1035--1039, 2008.

\bibitem{PG93}
A.K. Parekh and R.G. Gallager.
\newblock A generalized processor sharing approach to flow control in
  integrated services networks: the single-node case.
\newblock {\em IEEE/ACM Transactions on Networking}, 1(3):344--357, 1993.

\bibitem{R06}
K.~Ramanan.
\newblock {Reflected diffusions defined via the extended Skorokhod map}.
\newblock {\em Electronic Journal of Probability}, 11(36):934--992, 2006.

\bibitem{RR03}
K.~Ramanan and M.I. Reiman.
\newblock {Fluid and heavy traffic diffusion limits for a generalized processor
  sharing model}.
\newblock {\em Ann.~Appl.~Probab.}, 13(1):100 -- 139, 2003.

\bibitem{R86}
S.~Ramasubramanian.
\newblock Hitting a boundary point by diffusions in the closed half space.
\newblock {\em Journal of Multivariate Analysis}, 20(1):143--154, 1986.

\bibitem{R88}
S.~Ramasubramanian.
\newblock Hitting of submanifolds by diffusions.
\newblock {\em Probab.~Th.~Rel.~Fields}, 78:149--163, 1988.

\bibitem{R94}
J.A. Rebholz.
\newblock Planar diffusions with applications to mathematical finance.
\newblock {\em Ph.D.~thesis, University of California, Berkeley}, 1994.

\bibitem{R84}
Martin~I. Reiman.
\newblock Open queueing networks in heavy traffic.
\newblock {\em Mathematics of Operations Research}, 9(3):441--458, 1984.

\bibitem{T90}
L.M. Taylor.
\newblock {Existence and uniqueness of semimartingale reflecting Brownian
  motions in an orthant}.
\newblock Ph.D.~Thesis -- University of California, San Diego, 1990.

\bibitem{TW93}
L.M. Taylor and R.J. Williams.
\newblock {Existence and uniqueness of semimartingale reflecting Brownian
  motions in an orthant}.
\newblock {\em Probab. Theory Relat. Fields}, 96:283--317, 1993.

\bibitem{TY13}
A.~Tourin and R.~Yan.
\newblock Dynamic pairs trading using the stochastic control approach.
\newblock {\em Journal of Economic Dynamics and Control}, 37(10):1972--1981,
  2013.

\bibitem{VW85}
S.R.S. Varadhan and R.J. Williams.
\newblock {Brownian motion in a wedge with oblique reflection}.
\newblock {\em Communications on Pure and Applied Mathematics}, 38(4):405--443,
  1985.

\bibitem{W80}
W.~Whitt.
\newblock {Some useful functions for functional limit theorems}.
\newblock {\em Math. Oper. Res.}, 5:67--85, 1980.

\bibitem{W02}
W.~Whitt.
\newblock {\em {Stochastic-Process Limits -- an introduction to
  stochastic-process limits and their applications to queues}}.
\newblock Springer Series in Operation Research, 2002.

\bibitem{W85}
R.~J. Williams.
\newblock {Reflected Brownian motion in a wedge: semimartingale property}.
\newblock {\em Z. Wahr. verw. Geb.}, 69:161--176, 1985.

\bibitem{Z17}
L.~Zambotti.
\newblock {\em {Random obstacle problems}}.
\newblock Lecture Notes in Math. 2181. Springer, Cham, 2017.

\end{thebibliography}
\end{document}